\documentclass[10pt]{amsart}
\usepackage[hidelinks]{hyperref}
\usepackage{doi,url}
\usepackage[utf8]{inputenc}
\usepackage{amsmath}
\usepackage{amssymb}
\usepackage[foot]{amsaddr}
\usepackage[numbers]{natbib}
\DeclareRobustCommand{\noopsort}[1]{}
\usepackage{prettyref}
\usepackage{tikz-cd}
\usepackage{graphicx}
\usepackage{enumitem}
\usepackage{imakeidx}
\usepackage{mathtools}

\usepackage{comment}

\numberwithin{equation}{section}

\theoremstyle{plain}
\newtheorem{theorem}{Theorem}[section]
\newtheorem{proposition}[theorem]{Proposition}
\newtheorem{corollary}[theorem]{Corollary}

\newtheorem{lemma}[theorem]{Lemma}
\newtheorem{fact}[theorem]{Fact}
\newtheorem{question}[theorem]{Question}

\newcounter{introthmcounter}
\theoremstyle{definition}

\newtheorem{definition}[theorem]{Definition}

\newtheorem{example}[theorem]{Example}

\newtheorem{notation}[theorem]{Notation}

\theoremstyle{remark}
\newtheorem{remark}[theorem]{Remark}
\newtheorem{assumption}[theorem]{Assumption}

\newcommand{\bE}{\mathbb E}

\newcommand{\bR}{\mathbb R}
\newcommand{\bN}{\mathbb N}

\newcommand{\On}{\mathbf{On}}
\newcommand{\Mod}{\mathrm{Mod}}

\DeclareMathOperator{\Fin}{\mathrm{Fin}}
\DeclareMathOperator{\Supp}{\mathrm{Supp}}
\DeclareMathOperator{\Func}{\mathrm{Func}}
\DeclareMathOperator{\Nat}{\mathrm{Nat}}

\DeclareMathOperator{\Ckr}{\mathrm{Ckr}}
\DeclareMathOperator{\ckr}{\mathrm{ckr}}
\DeclareMathOperator{\Kr}{\mathrm{Ker}}
\DeclareMathOperator{\kr}{\mathrm{ker}}
\DeclareMathOperator{\Img}{\mathrm{Im}}
\DeclareMathOperator{\img}{\mathrm{im}}
\DeclareMathOperator{\Coimg}{\mathrm{Coim}}
\DeclareMathOperator{\coimg}{\mathrm{coim}}

\DeclareMathOperator{\Yon}{\mathbb{Y}\!}
\DeclareMathOperator{\Lan}{\mathrm{Lan}}
\DeclareMathOperator{\Ann}{\mathrm{Ann}}

\DeclareMathOperator{\Hom}{\mathrm{Hom}}
\DeclareMathOperator{\Span}{\mathrm{Span}}
\DeclareMathOperator{\Sub}{\mathrm{Sub}}
\DeclareMathOperator{\Deriv}{\mathrm{Der}}
\DeclareMathOperator{\Id}{\mathrm{Id}}

\DeclareMathOperator{\Free}{\mathrm{F}}
\DeclareMathOperator{\Tors}{\mathrm{T}}

\newcommand{\hotimes}{\mathbin{\hat{\otimes}}}
\newcommand{\ort}{\mathbin{\bot}}

\newcommand{\id}{\mathrm{id}}
\newcommand{\cN}{\mathcal{N}}
\newcommand{\cA}{\mathcal{A}}
\newcommand{\cQ}{\mathcal{Q}}
\newcommand{\cS}{\mathcal{S}}
\newcommand{\ev}{\mathrm{ev}}
\newcommand{\ad}{\mathrm{ad}}

\newcommand{\Ob}{\mathrm{Ob}}

\def\bfk{\mathbf{k}}

\def\op{\mathrm{op}}

\def\cP{\mathcal{P}}
\def\cF{\mathcal{F}}
\def\cG{\mathcal{G}}
\def\cC{\mathcal{C}}
\def\cE{\mathcal{E}}
\def\cT{\mathcal{T}}
\def\cD{\mathcal{D}}
\def\cH{\mathcal{H}}
\def\cR{\mathcal{R}}
\def\cM{\mathcal{M}}
\def\cW{\mathcal{W}}
\def\-{\text{-}}
\def\bfx{\mathbf{x}}

\newcommand{\Ind}{\mathrm{Ind}}
\newcommand{\Sets}{\mathrm{Set}}
\newcommand{\Vect}{\mathrm{Vect}}
\newcommand{\TVect}{\mathrm{TVect}}
\newcommand{\Card}{\mathrm{Card}}
\newcommand{\WO}{\mathrm{WO}}

\newcommand{\Monomorphisms}{\mathrm{Mono}}
\newcommand{\Epimorphisms}{\mathrm{Epi}}

\newcommand{\two}{2}

\title{On vector spaces with formal infinite sums}
\author{Pietro Freni}



\address{Institute of Mathematics, Czech Academy of Sciences}

\email{pietro.freni@hotmail.it}

\begin{document}

\begin{abstract}
	I discuss possible definitions of categories of vector spaces enriched with a notion of formal infinite linear combination in the likes of the formal infinite linear combinations one has in the context of generalized power series, I call these \emph{reasonable categories of strong vector spaces} (r.c.s.v.s.). I show that, in a precise sense, the more general possible definition for a strong vector space is that of a small $\Vect$-enriched endofunctor of $\Vect$ that is right orthogonal for every cardinal $\lambda$, to the cokernel of the canonical inclusion of the $\lambda$-th copower in the $\lambda$-th power of the identity functor: these form the objects for a universal r.c.s.v.s.\ I call $\Sigma\Vect$. I show this is equivalent to the category of \emph{ultrafinite summability spaces} defined independently in \cite{bagayoko2024automorphisms}.
	I relate this category to what could be understood to be the obvious category of strong vector spaces $B\Sigma\Vect$ and to the r.c.s.v.s.\ $K\TVect_s$ of separated linearly topologized spaces that are generated by linearly compact spaces.
	I analyze the monoidal closed structures on various r.c.s.v.s.\  induced by the natural one on $\Ind\-(\Vect^\op)$. In particular with respect to the problem of closure under the tensor product of $\Ind\-(\Vect^\op)$.
	Most of the technical results apply to a more general class of orthogonal subcategories of $\Ind\-(\Vect^\op)$ and we work with that generality as it's cost-free.
\end{abstract}

\maketitle
\tableofcontents

\section{Introduction}

\subsection{Motivation}
    The set of generalized power series over an arbitrary field $\bfk$ in the sense of Hahn-Higman-Ribenboim (cf \cite{hahn1907nichtarchimedeischen}, \cite{higman1952ordering}, \cite{ribenboim1997semisimple}) is built up from an ordered set $(\Gamma, <)$ as the $\bfk$-vector space
	$\bfk[[\Gamma]]$ of functions $f: \Gamma \to \bfk$ whose support $\Supp f =\{\gamma \in \Gamma : f(\gamma)\neq 0\}$ is Noetherian (i.e.\ is a well partial order).
	If $\Gamma$ is endowed with an ordered monoid structure, then $\bfk[[\Gamma]]$ is naturally an algebra endowed with the straightforward generalization of the Cauchy product for power-series
	\[(f \cdot g)(\gamma)= \sum_{\substack{\gamma_0\cdot \gamma_1=\gamma\\ f(\gamma_0)g(\gamma_1)\neq 0}} f(\gamma_0)g(\gamma_1)\]
	which is well defined because of the Noetherianity constraint on the supports.
		
	These are important constructions in model theory and algebra of valued fields: for example, if $\Gamma$ is a totally ordered Abelian group, then the Hahn field $\bfk[[\Gamma]]$ is maximally valued and every $0$-equicharacteristic valued field $(K,v)$ with a residue field admitting $n$-th roots for every element is an immediate subfield of some Hahn field (Kaplansky \cite{kaplansky1942maximal}).
		
	Refining the Noetherianity constraint on the support, as imposing a bound on cardinality, results in various other algebras of series, often notably motivated by the need to define exponential operators (cf \cite{kuhlmann2005kappa}, \cite{krapp2022rayner}). Generally one can define spaces of restricted series $\bfk(\Gamma; \cF)$ consisting of functions $f: \Gamma \to \bfk$ whose support lies in the specified ideal $\cF$ of Noetherian subsets of $\Gamma$ (compare with \cite[Def.~6.2]{berarducci2021value} and \cite{krapp2022rayner}).
		
	A prominent feature of spaces of generalized power series is the notion of formal infinite sum, or rather more appropriately, of \emph{infinite $\bfk$-linear combination}: a family $(f_i)_{i \in I}\in \bfk(\Gamma; \cF)^I$ indexed by an arbitrary set $I$, is said to be \emph{summable} when $\bigcup_{i \in I}\Supp f_i\in \cF$ and for every $\gamma \in \Gamma$ the set $\{i \in I: \gamma \in \Supp f_i\}$ is finite. Given such a family it is possible to associate to any $(k_i)_{i \in I} \in \bfk^I$ the infinite linear combination $f=\sum_{i \in I} f_i k_i$ given by $f(\gamma)=\sum(f_i(\gamma): \gamma \notin \Supp(f_i))$.
	This feature is referred throughout literature as a \emph{strongly linear} structure on the space $\bfk(\Gamma; \cF)$, and maps $F:\bfk(\Gamma; \cF) \to \bfk(\Delta; \cG)$ preserving it (i.e. preserving all the $I$-ary relations $\sum_{i \in I} f_i k_i =f$, for every set $I$ and every $(k_i)_{i \in I} \in \bfk^I$) are referred to as \emph{strongly linear maps} (see e.g.\ \cite{dries2001logarithmic}, \cite{berarducci2015surreal}, \cite{hoeven2001operators}).
		
	In many applications, the strongly linear structure plays an important conceptual and algebraic role. For example: when $\Gamma$ is an ordered monoid, the product defined on $\bfk[[\Gamma]]$ is the only extension preserving infinite $\bfk$-linear combinations in both arguments of the product on $\Gamma$ along the natural inclusion $\Gamma \subseteq\ \bfk[[\Gamma]]$ via indicator functions; moreover in the case $\Gamma$ is a totally ordered group and $f \in \bfk[[\Gamma]]$ is such that $\Supp f>1$ then $(f^n)_{n \in \mathbb{N}}$ is summable (Neumann \cite{neumann1949ordered}), i.e.\ it is possible to evaluate every power series at $f$. Finally when enriching algebras of generalized power series with $\bfk$-linear operators such as derivations, it is a natural requirement to have them be strongly linear as well (the standard derivation on transseries and the Berarducci-Mantova derivation on Surreal numbers are strongly linear cf \cite{hoeven2001operators},\cite{berarducci2015surreal}).
		
	The interest in strongly linear derivations on algebras of generalized series in particular is also a good motivation for the study of an appropriate "infinite-sum-sensitive" category of modules for such algebras.
	
	Toward this goal an intermediate step would be to define an appropriate monoidal closed category $\Sigma\Vect$ of $\bfk$-vector spaces enriched with a strongly linear structure (i.e.\ a notion of infinite linear combination) in such a way that $\bfk$-algebras of generalized power series are monoids in $\Sigma\Vect$.

\subsection{Overview}
	The main goal of this article is to investigate the notion of formal infinite sum providing a properly justified definition of the wanted $\Sigma\Vect$ and relating it to the category of linearly topologized $\bfk$-vector spaces. Throughout the paper, we fix an arbitrary field $\bfk$; linear will mean $\bfk$-linear and vector space will mean $\bfk$-vector space. The letters $I,J$ will always denote arbitrary index sets.
    
	\smallskip
	In Section~\ref{sec:BasSL} a special category of strong vector spaces is introduced. This category which we call of \emph{based strong vector spaces} $B\Sigma\Vect$ is defined following what could be understood to be the ``obvious generality'': such spaces have still ``series-like'' elements, and are easier to think of.
	Some elementary results in this context are proven which will be mainly useful when working with concrete examples. In particular an appropriate notion of dual in this category allows for some pleasantly simplifying reductions. 
		
	\smallskip

	In Section~\ref{sec:GenSL} we question whether this is the right generality and address the possible generalizations. We will argue (Remark~\ref{rmk:Eur_0}) that any \emph{reasonable} ($\bfk$-additive) category of strong vector spaces the following should be required to hold 
		
	\begin{enumerate}
		\item the category whose objects are spaces of the form $\bfk^I=\bfk(I; \cP(I))$ and whose maps are the strongly linear maps among them should be a dense full subcategory;
		\item any map from some $f: \bfk^I\to X$ should be uniquely determined by the compositions
		\[\begin{tikzcd}
		    \bfk \ar[r, "\delta_i"] &\bfk^I \ar[r, "f"] &X
		\end{tikzcd}\]
		as $i$ ranges in $I$ and where $\delta_i$ sends $1$ to the function $I \to \bfk$ taking value $1$ on $i$ and $0$ everywhere else. This essentially entails that infinite linear combinations are functional and produce a unique output.
	\end{enumerate}

	So we set these conditions as the definition of \emph{reasonable category of strong vector spaces} (Definition~\ref{defn:reasonable}).
		
	It turns out that there is, up to equivalence, a unique category $\Sigma\Vect$ universal among the ones satisfying the two properties above (Theorem~\ref{thm:SigmaVectChar}).
		
	Up to equivalence such $\Sigma\Vect$ can be understood to be the category which has as objects $\bfk$-additive functors $X: \Vect \to \Vect$ such that there are no non-zero natural transformations from $\Vect(\bfk,-)^\lambda /\Vect(\bfk, -)^{\oplus \lambda}$ to $X$, and as arrows natural transformation between those functors. 

    The category of \emph{ultrafinite summability spaces} as defined independently in \cite{bagayoko2024automorphisms}, and strongly linear maps among them is also equivalent to $\Sigma\Vect$ (Theorem~\ref{thm:SigmaVect_as_mod}).

    \smallskip
    
	In Section~\ref{sec:OrtSmallRefl}, I will show that $\Sigma\Vect$ is in fact a reflective subcategory of the Abelian complete and cocomplete category $\Ind\-(\Vect^\op)$ of \emph{small} $\bfk$-additive functors from $\Vect$ to $\Vect$ (Corollary~\ref{cor:SigmaVect_is_reflective}). This entails its completeness and cocompleteness.
    In subsection~\ref{ssec:SigmaEmb} I will characterize maps which preserve and reflect infinite sums (\emph{$\Sigma$-closed embeddings}) as right orthogonals in a suitable orthogonal factorization system.
    
    Most of Section~\ref{sec:OrtSmallRefl} is presented in the more general setting of torsion theories of $\Ind\-(\Vect^\op)$ (e.g.\ the above mentioned Corollary~\ref{cor:SigmaVect_is_reflective} is an instance of the more general Theorem~\ref{thm:Q-free_are_reflective}). This generality will be useful in the last two sections.
		
	\smallskip
	Section~\ref{sec:Comp} compares $\Sigma\Vect$ with other better understood categories. After briefly making explicit the fully faithful embedding of $B\Sigma\Vect$ in $\Sigma\Vect$, attention is switched to another reasonable category of strong vector spaces which is properly topological in nature. It will be shown that a reflective subcategory of $\Sigma\Vect$ is equivalent to a specific reflective subcategory $K\TVect_{s}$ of the category $\TVect$ of linearly topologized vector spaces: in such a category infinite sums can be understood as topological limits and strongly linear maps as continuous linear maps. The category $K\TVect_{s}$ consists of the colimits in the category $\TVect_{s}$ of separated linearly topologized vector spaces, of \emph{linearly compact topological vector spaces} (in the sense of Lefschetz~\cite{lefschetz1942algebraic}). The results relating to this comparison are merely observed and are probably part of the folklore related to the theory of linearly topologized vector spaces and Lefschetz's duality.
		
	\smallskip
	It is worth mentioning here that as full subcategories of $\Ind\-(\Vect^\op)$
	$$B\Sigma\Vect \subseteq K\TVect_{s} \subsetneq \Sigma\Vect$$
	and that whether the first inclusion is an equality is a question left open (although we expect it to be strict).
	The inclusion $B\Sigma\Vect\subseteq K\TVect_{s}$ in particular implies that in the usual context of generalized series (in particular Hahn fields, or more generally, Rayner structures, cf \cite{krapp2022rayner}) preserving infinite sums can be understood in terms of continuity with respect to certain linear topologies.
	The separating example proving that the inclusion $K\TVect_{s} \subsetneq \Sigma\Vect$ is strict can already be presented here and will be used as an excuse to introduce a somewhat restricted version of various notions in the following. 
		
	If $\bfk^I=\bfk(I; \cP(I))$ is endowed with the product topology of discrete copies of $\bfk$, then $\sum_{j \in J} f_j$ as defined above coincides with the topological limit of the net $(\sum_{j \in S} f_j)_{S \in [I]^{<\omega}}$ of finite partial sums of the family $(f_j)_{j \in J}$. It turns out that in general maps which preserve infinite sums between objects of the form $\bfk^I$ for $I$ any set will coincide with continuous linear maps in the just mentioned product topology.
    However, unless $I$ is countable, subspaces which are closed under taking infinite sums of summable families need not to be topologically closed: take $I = \dim (\bfk^\bN)$, fix an isomorphism $\varphi: \bfk^{\oplus I} \to \bfk^\bN$, and and write $(e_n)_{n \in \bN}$, resp.\ $(e_i: i \in I)$ for the canonical bases of $\bfk^{\oplus \bN}$ and $\bfk^{\oplus I}$; for every $v=\sum_{n<N} k_n e_n \in \bfk^{\oplus \bN}$ let $f_v: I \to \bfk$ be defined as $f_v(i)= \sum_{n<N} \varphi(e_i)(n) k_n$, then $v \mapsto f_v$ is a linear map $\bfk^{\oplus \bN} \to \bfk^I$ and the subspace $H=\{f_v: v \in \bfk^{\oplus \bN}\}$ is topologically dense in $\bfk^I$, but all summable families $(g_j: j \in J) \in H^J$ must be such that $\{j \in J : g_j\neq 0\}$ is finite, whence $H$ is (trivially) closed under infinite sums of summable families in $\bfk^I$.
	As a consequence it is possible to define a notion of infinite sum on $\bfk^I/H$ but the quotient linear topology on $\bfk^I/H$ will be the trivial (indiscrete) topology. Example~\ref{ex:weirdo} and Proposition~\ref{prop:SigmaEmb&Quot} will give a more natural view of this.
		
	\smallskip
	In Section~\ref{sec:ClMon} a natural monoidal closed structure on the previously mentioned categories is defined. In fact the canonical monoidal closed structure on $\Ind\-(\Vect^\op)$ induces a natural monoidal closed structure on all $\Sigma\Vect$, $K\TVect_{s}$ and $B\Sigma\Vect$. Moreover the tensor product of the monoidal structure interacts as expected with infinite sums.

    I also answer basic questions as to whether $\Sigma\Vect$, $K\TVect_{s}$, and $B\Sigma\Vect$ are closed under the tensor product and internal hom of $\Ind\-(\Vect^\op)$ and give an explicit descriptions of those functors in $B\Sigma\Vect$.    
		
	The tensor product allows a natural definition of strongly linear algebras and strongly linear modules for a strongly linear algebra. 
		
	Finally I apply this result to give a notion of strongly linear Kähler differential for a strongly linear algebra.
		
	\smallskip
	The exposition is meant for a general audience who has some familiarity with commutative algebra and category theory. Most of the machinery is introduced only when needed and recalled before use. Section~\ref{sec:Comp} and Section~\ref{sec:ClMon} however make some limited use of the end-coend notation: for the purpose of this paper, the account of those notions made in \cite[Ch.~IX, Sec.s~5-6]{maclane71categories} is more than enough, a more detailed and comprehensive account on ends-coend is the book \cite{loregian2021coend}. 

\subsection{Addendum}
	The author has recently been made aware of the notion of contra-module presented for example in the exposition \cite{positselski2015contramodules}. Contra-modules over a coalgebra $C$ over a field $\bfk$, once a $\bfk$-basis $B$ of the supporting space of $C$ is chosen, can be seen as endowed with a notion of $B$-indexed sum, satisfying some compatibility constraints with the co-algebra structure. This notion of formal infinite sum seems different from the one presented in this note.
	There is a relation whose further investigation could be fruitful, but the author didn't carry it out.

\subsection{Acknowledgments}
	The author is very grateful to his supervisors Vincenzo Mantova for his support and feedback and Dugald MacPherson for his support. Thanks to Ivan Di Liberti for a chat about the topic, in particular for the suggestion to look at category-theoretical orthogonality - at the time unknown to the author - which resulted in great simplification of big chunks of the treatment. Thanks also to Vincent Bagayoko for some useful conversations about this and to the anonymous referee for several suggestion which improved the paper.
	This paper is part of a PhD project at the University of Leeds and as such is supported by a Scholarship from the University of Leeds' School of Mathematics. Its content is included in Chapter~2 of the author's doctoral dissertation \cite{freni2024structural}.
	
\section{Based strong vector spaces}\label{sec:BasSL}

In this section we will deal with spaces of generalized series of the form $\bfk(\Gamma; \cF)$ mentioned in the introduction. Since we will be interested just in strongly linear maps among these, we will drop the assumption that $\Gamma$ is ordered and let $\cF$ be any ideal of subsets containing all singletons. By doing so, we are not really dealing with a more general situation because, given an infinite set $\Gamma$, we can choose a cardinal well order on it (and even endow it with a cancellative ordered monoid structure given by Hessenberg ordinal addition); with such a choice, every ideal $\cF$ of subsets of $\Gamma$ consists of Noetherian subsets.

\subsection{Based strong vector spaces}
%
%
The definition below is essentially \cite[Def.~6.2]{berarducci2021value} or \cite[Def.~3.2]{krapp2022rayner} after forgetting about the order on $\Gamma$.
	
\begin{definition}\label{def:summable1}
	For a set $\Gamma$, let $[\Gamma]^{<\omega}$ denote the set of finite subsets of $\Gamma$ and $\cP(\Gamma)$ the Boolean algebra of all subsets of $\Gamma$.
	Call a \emph{bornology} on $\Gamma$ an ideal in the quotient Boolean algebra $\cP(X)/[\Gamma]^{<\omega}$. Bornologies are in bijection with ideals in $\cP(\Gamma)$ containing all singletons. If $\cF$ is a bornology on $\Gamma$ denote by $\widetilde{\cF}$ (or $\cF^{\sim}$) the corresponding ideal in $\cP(\Gamma)$. Call a subset $F \subseteq X$, \emph{$\cF$-bounded}, if it is in $\widetilde{\cF}$.
		
	For $\cF$ a bornology on $\Gamma$ define
	\[\bfk(\Gamma; \cF):=\{f \in \bfk^\Gamma:\: \Supp(f) \in \widetilde{\cF}\} \cong \bigcup_{F \in \cF} \bfk^F \subseteq \bfk^\Gamma\]
	to be the set of functions with $\cF$-bounded support ($\Supp(f):=\{\gamma \in \Gamma: f(\gamma) \neq 0\}$ is the \emph{support} of $f$).
    This is naturally a subspace of $\bfk^\Gamma$.
	Call a family $(f_i : i \in I)$ in $\bfk(\Gamma; \cF)$ \emph{summable} if for all $\gamma \in \Gamma$, $\{i \in I : \gamma \in \Supp(f_i)\}$ is finite and $\bigcup_{i \in I} \Supp(f_i) \in \widetilde{\cF}$. Define the sum of the family by the formula
	\[\left(\sum_{i \in I} f_i\right)(\gamma):=\sum_{\substack{i \in I\\ \Supp(f_i)\ni \gamma}} f_i(\gamma) \qquad \text{for all}\; \gamma \in \Gamma.\]
    Since $\bigcup_{i \in I} \Supp(f_i) \in \widetilde{\cF}$, $\sum_{i \in I} f_i \in \bfk(\Gamma, \cF)$.
    
	Call a $\bfk$-linear map $F : \bfk(\Gamma;\cF) \to \bfk(\Delta; \cG)$
	\begin{enumerate}
		\item \emph{$\Sigma$-preserving} (or \emph{strongly $\bfk$-linear}) if whenever $(f_i : i \in I)$ is summable in $\bfk(\Gamma;\cF)$, it happens that $(Ff_i: i \in I)$ is summable in $\bfk(\Delta; \cG)$ and $F \sum_{i \in I} f_i = \sum_{i \in I} F f_i$.
		\item \emph{$\Sigma$-reflecting} if whenever $(F f_i : i \in I)$ is summable in $\bfk(\Delta; \cG)$ it happens that $(f_i : i \in I)$ is summable in $\bfk(\Gamma; \cF)$ and $F \sum_{i \in I} f_i = \sum_{i \in I} F f_i$.
		\item \emph{weakly $\Sigma$-reflecting} if whenever $(F f_i : i \in I)$ is summable in $\bfk(\Delta; \cG)$ and $\sum_{i \in I} k_i \cdot F f_i$ is in the image of $F$ for every $(k_i)_{i \in I} \in \bfk^I$, it happens that $(f_i : i \in I)$ is summable in $\bfk(\Gamma; \cF)$ and $F \sum_{i \in I} f_i = \sum_{i \in I} F f_i$
	\end{enumerate}
	We call $\bfk(\Gamma; \cF)$ a \emph{based $\Sigma$-$\bfk$-vector space} or \emph{based strong $\bfk$-vector space}. We denote the category of based strong vector spaces and strongly linear maps by $B\Sigma\Vect$. As the field $\bfk$ will usually be fixed throughout, the subscript $\bfk$ will often be dropped (especially in diagrams).
\end{definition}

\begin{remark}\label{rmk:sum_of_summable_families}
    Note that if $k,h \in \bfk$ and $(f_i: i \in I)$ and $(g_i: i \in I)$ are two $I$-indexed summable families in $\bfk(\Gamma; \cF)$, then $(k \cdot f_i+ h \cdot g_i :i \in I)$ is summable because $\{i \in I : k \cdot f_i+ h \cdot g_i\neq 0\} \subseteq \{i \in I: f_i \neq 0\} \cup \{i \in I: g_i \neq 0\}$ and $\Supp(k \cdot f_i+h \cdot g_i)\subseteq \Supp(f_i) \cup \Supp(g_i)$. 
    Moreover it is easy to see that $\sum_{i \in I} k \cdot f_i+ h \cdot g_i= k \cdot \left(\sum_{i \in I} f_i\right) + h \cdot \left(\sum_{i \in I} g_i\right)$.

    It follows that if $F,G: \bfk (\Gamma; \cF) \to \bfk(\Delta; \cG)$ are strongly linear maps, then so is $k \cdot F + h \cdot G$. In particular the $B\Sigma\Vect(\bfk(\Gamma; \cF), \bfk(\Delta; \cG))$ is naturally a $\bfk$-vector space with the pointwise sum and scalar multiplication.
\end{remark}

\begin{remark}\label{rmk:Sigma_ref_is_weakly_Sigma_ref}
    A $\Sigma$-reflecting map is weakly $\Sigma$-reflecting.
    Note that an infinite constant non-zero family cannot be summable. In particular, a weakly $\Sigma$-reflecting map $F: \bfk(\Gamma; \cF) \to\bfk(\Delta; \cG)$ is necessarily injective.
\end{remark}

\begin{remark}\label{rmk:Sigma_emb_B}
    A map $F:\bfk(\Gamma; \cF) \to \bfk(\Delta; \cG)$ that is $\Sigma$-reflecting and $\Sigma$-preserving has an image $H:=F(\bfk(\Gamma; \cF))$ which is \emph{closed under infinite sums} (or \emph{$\Sigma$-closed}): this means that if $(g_i:i \in I)$ is a summable family of $\bfk(\Delta;\cG)$ such that $g_i \in H$ for all $i\in I$, also $\sum_{i \in I} g_i\in H$.
    We call such an $F$ a \emph{$\Sigma$-closed embedding}.
    Maps that are $\Sigma$-preserving and weakly $\Sigma$-reflecting will be called \emph{embeddings}.
\end{remark}
	
\begin{remark}\label{rmk:monomials_inclusion}
	For every $\gamma \in \Gamma$, we have that the indicator function of $\gamma$, $\delta_{\gamma}: \Gamma \to \bfk$ is an element of $\bfk(\Gamma; \cF)$. One usually denotes $\delta_{\gamma}$ by $\gamma$ and sees $\Gamma$ as a subset of $\bfk(\Gamma; \cF)$. With this identification, for every $F \in \cF$ and every $(k_\gamma)_{\gamma \in F} \in \bfk^F$ the family $(\gamma\cdot k_\gamma: \gamma \in S)$ is summable and for every $f \in \bfk(\Gamma; \cF)$ one has $f=\sum_{\gamma \in \Supp(f)}\gamma\cdot f(\gamma)$.
\end{remark}

\begin{remark}\label{rmk:density}
	Among based strong vector spaces, the ones of the form $\bfk^I=\bfk(I; \cP(I))$ play a special role, in that, by construction, every based strong vector space $\bfk(\Gamma; \cF)$ is a directed union of such spaces along a diagram of $\Sigma$-closed embeddings.
\end{remark}

We have the following correspondence between strongly linear maps $\bar{f}:\bfk(I; \cP(I)) \to \bfk(\Gamma; \cF)$ and $I$-indexed summable families $(f_i)_{i \in I}$ in $\bfk(\Gamma; \cF)$.

\begin{lemma}\label{lem:densityBSigma}
    Every summable family $(f_i)_{i \in I}$ in some $\bfk(\Gamma; \cF)$ determines a strongly linear map $\bfk(I; \cP(I)) \to \bfk(\Gamma; \cF)$ given by 
	\begin{equation}\label{ch2:eq1}
	    \bfk(I; \cP(I))=\bfk^I\ni (k_i)_{i \in I} \mapsto \sum_{i\in I} f_i \cdot k_i \in \bfk(\Gamma; \cF)
	\end{equation}
	Conversely every strongly linear map $\bar{f}: \bfk(I; \cP(I)) \to \bfk(\Gamma; \cF)$ arises this way from the summable family $(f_i)_{I \in I}:=\big(\bar{f} (\delta_i)\big)_{i \in I}$.
    \begin{proof}
        Let $\big((k_{i,j})_{i \in I}\big)_{j \in J}$ be a summable $J$-indexed family in $\bfk^{I}$, i.e.\ for all $i\in I$, $\{j\in J: k_{i,j}\neq 0\}$ is finite. To prove the first assertion we shall prove that setting $g_j:=\sum_{i \in I} k_{i,j} f_i$ defines a family $(g_j)_{j \in J}$ which is summable in $\bfk(\Gamma; \cF)$, and that $\sum_{j \in J} g_j = \sum_{i \in I} f_i \cdot \left(\sum_{j\in J} k_{i,j}\right)$. Let $\gamma \in \Gamma$ and note that
        \[\{j \in J: g_i(\gamma) \neq 0\}\subseteq \bigcup_{\substack{i\in I\\f_i(\gamma)\neq 0}}\{j \in J : k_{i,j} \neq 0\},\]
        where the right-hand side is a finite union of finite sets by the hypotheses on $(f_i)_{i \in I}$ and $\big((k_{i,j})_{i \in I}\big)_{j \in J}$. It also follows that
        \[\sum_{j} g_j(\gamma) = \sum_{i \in I} f_i (\gamma) \cdot \left(\sum_{j\in J} k_{i,j}\right).\]
        Finally note that for every $j\in J$ we have $\Supp(g_j)\subseteq \bigcup_{i \in I} \Supp(f_i)$, whence $\bigcup_{j \in J} \Supp(g_j) \subseteq \bigcup_{i \in I} \Supp(f_i) \in \cF$. Thus we have proved the first assertion of the lemma.
        
        For the second one, note that a strongly linear map $\bar{f}: \bfk(I; \cP(I)) \to \bfk(\Delta; \cG)$, since $(\delta_i)_{i \in I}$ is a summable family in $\bfk(I; \cP(I))=\bfk^I$, must be given by (\ref{ch2:eq1}) with $f_i:=\bar{f}(\delta_i)$.
    \end{proof}
\end{lemma}

\subsection{Some duality reductions}
	
The definition of $B\Sigma\Vect$ is rather unsatisfactory because objects are defined in terms of a set $\Gamma \subseteq \bfk(\Gamma; \cF)$ which is certainly not invariant under strongly linear automorphisms.
In view of Remark~\ref{rmk:density} and Lemma~\ref{lem:densityBSigma}, in order to understand $B\Sigma\Vect$ and its possible further generalizations, it is important to understand the full subcategory on the objects of the form $\bfk(I; \cP(I))=\bfk^I$. We will show, and it's the main point of this section, that this is equivalent to the category $\Vect^\op$ opposite to the category of vector spaces and linear maps (Corollary~\ref{cor:Vectop_equi}). Since the arguments of the proof easily lead to more general and somewhat useful considerations on $B\Sigma\Vect$, we chose to present them here at the cost of a slight increase in length. This will be useful in some examples and in providing another description of the category $\Sigma\Vect$ of general strong vector spaces (see Corollary~\ref{cor:fromBtogen}).

Throughout the rest of the paper we will make use of the following standard notations of linear algebra: for a subset $H$ of a $\bfk$-vector space $V$, we will denote by $\Span(H)\le V$ the $\bfk$-linear span of $H$ and by
\[\Ann(H):=\{\xi \in V^*: \forall h \in H, \;\xi(h)=0\} \le V^*:=\Vect(V, \bfk)\]
the \emph{annihilator} of $H$.
	
\begin{remark}
	The Stone space of the Boolean algebra $\cP(\Gamma)/[\Gamma]^{<\omega}$ is naturally the subspace of non-isolated points in the space $\beta\Gamma$ of ultrafilters in $\cP(\Gamma)$ (or, equivalently, of maximal ideals). $\Gamma$ is naturally embedded in $\beta\Gamma$ via principal ultrafilters (ideals).
	Filters (ideals) in $\cP(\Gamma)/[\Gamma]^{<\omega}$ correspond to closed (open) subspaces in $\beta\Gamma \setminus \Gamma$.
    We refer the reader to Stone's original paper \cite{stone1936theory} or Chapter~1 of Johnstone's book \cite{johnstone1982stone} for more on the correspondence between Boolean algebras and Stone spaces. The now standard notation $\beta\Gamma$ for the space of ultrafilters of $\cP(\Gamma)$ comes from \cite{cech1937bicompact}, as it is the Stone-Čech compactification of the discrete space $\Gamma$.
\end{remark}
	
\begin{definition}
	If $(A, \land, \lor, \top, \bot)$ is a Boolean algebra define the operator $-^\bot : I(A) \to I(A)$ on the set $I(A)$ of ideals of $A$ as
	\[\cF \mapsto \{a \in A: \forall x \in \cF, a \wedge x = \bot \}=: \cF^\bot.\]
\end{definition}
	
\begin{fact}
	For $A$ a Boolean algebra $-^\bot: I(A) \to I(A)$ is an antitone Galois connection (w.r.t. the $\subseteq$ order on ideals). It corresponds to the map sending an open set of the Stone space of $A$ to the interior of its complement.
\end{fact}
	
\begin{remark}
	If $\cF$ is a bornology on $\Gamma$, then $\cF^\bot$ is still a bornology on $\Gamma$ and
	\[(\cF^\bot)^\sim = \{S \subseteq \Gamma: \forall F \in \cF^\sim,\; |F \cap S|< \aleph_0 \}.\]
\end{remark}
	
\begin{definition}\label{def:dual_borno}
	For $\cF$ a bornology on $\Gamma$ define the $\bfk$-bilinear map
	\[\langle -, -\rangle: 
	\bfk(\Gamma; \cF) \times \bfk(\Gamma; \cF^\bot) \to \bfk
	\qquad \langle f , g \rangle = \sum_{\gamma \in \Gamma} f(\gamma)g(\gamma)\]
    Note that the sum on the right is well defined as 
	\[\{\gamma \in \Gamma: f(\gamma) \neq 0,\; g(\gamma) \neq 0 \}=\Supp(f)\cap \Supp(g) \in [\Gamma]^{<\omega}\]
	by construction of $\cF^\bot$.
\end{definition}
	
\begin{lemma}\label{lem:Dual}
	A linear functional $\xi : \bfk(\Gamma; \cF) \to \bfk$ is $\Sigma$-preserving if and only if it is of the form $\xi(-)= \langle -, g\rangle$ for some $g \in \bfk(\Gamma; \cF^\bot)$.
    \begin{proof}
		Let $(f_i : i \in I)$ be summable in $\bfk(\Gamma; \cF)$ and $g \in \bfk(\Gamma; \cF^\bot)$, then
		\[\begin{aligned}
		\{i \in I: \langle f_i, g \rangle \neq 0 \} \subseteq
		\\ \bigcup\left\{ \big\{i \in I : f_i(\gamma) \neq 0, g(\gamma) \neq 0\big\} : \gamma \in \Supp(g) \cap \bigcup_{i \in I} \Supp(f_i) \right\}
		\end{aligned}\]
		is finite as it is contained in a finite union of finite sets. Also clearly
		$\langle \sum_{i\in I} f_i, g \rangle =\sum_{i \in I} \langle f_i, g\rangle$.\\
		On the other hand assume $\xi$ is $\Sigma$-preserving, define
		$g(\gamma):=\xi(\gamma)$ ($\gamma$ ranges in $\Gamma$, $\gamma$ is a valid input for $\xi$ via the inclusion $\Gamma \subseteq \bfk(\Gamma; \cF)$ of Remark~\ref{rmk:monomials_inclusion}).
		First note that $g$ belongs to $\bfk(\Gamma; \cF^\bot)$: if not then there would be $F \in \cF^\sim$ with $|F \cap \Supp(g)|\ge \aleph_0$ so $(\xi(\gamma): \gamma \in F)$ would not be summable in $\Gamma$ whereas $(\gamma: \gamma\in F)$ is summable in $\bfk(\Gamma; \cF)$ contradicting that $\xi$ is $\Sigma$-preserving. The fact that $\langle f, g \rangle = \xi(f)$ for every $f \in \bfk(\Gamma; \cF)$ follows from the fact that $\langle -, g\rangle$ is $\Sigma$-preserving applied to the equality $f= \sum_{\gamma \in \Supp(f)} \gamma \cdot f(\gamma)$.
	\end{proof}
\end{lemma}
	
The following Proposition yields an alternative characterization of summability in $\bfk(\Gamma; \cF)$ when $\cF=\cG^{\bot}$ for some bornology $\cG$ (or equivalently when $\cF=\cF^{\bot\bot}$).
	
\begin{proposition}\label{prop:sum_imp}The following hold
	\begin{enumerate}
		\item If $(f_i: i \in I)$ in $\bfk(\Gamma; \cF)$ is summable then
		for every $g\in \bfk(\Gamma; \cF^\bot)$, $(\langle f_i, g\rangle: i \in I)$ is summable in $\bfk$.
		\item If $(\eta_i: i \in I)$ in $\bfk(\Gamma; \cF^{\bot})$ is such that for every $g\in \bfk(\Gamma; \cF)$, $(\langle \eta_i, g\rangle : i \in I)$ is summable in $\bfk$ then $(\eta_i : i \in I)$ is summable in $\bfk(\Gamma; \cF^{\bot})$.
	\end{enumerate}
    \begin{proof}
		$(1)$ is part of Lemma~\ref{lem:Dual}.
		For $(2)$, assume toward contradiction that $(\eta_i: i \in I)$ is not summable. The instance of the hypothesis given by setting $g=\gamma \in \Gamma$, says that $\{i \in I: \eta_i(\gamma)\neq 0\}$ is finite. So for $(\eta_i: i \in I)$ not to be summable we must have $\bigcup_{i \in I} \Supp(\eta_i)\notin (\cF^{\bot})^{\sim}$, which means that there is an $S \in (\cF)^{\sim}$ such that $|S\cap \bigcup_{i \in I} \Supp(\eta_i)| \ge \aleph_0$. Without loss of generality we can assume furthermore that $S\subseteq \bigcup_{i \in I} \Supp(\eta_i)$ and $|S|=\aleph_0$. To reach a contradiction we show there is $g\in \bfk(S; \cP(S))$ and an infinite $\{i_k : k <\omega\} \subseteq I$ such that $\langle g, \eta_{i_k}\rangle \neq 0$ for every $k$.
		Let $\pi: \bfk(\Gamma; \cF^{\bot}) \to \bfk(S; [S]^{<\omega})$ be the restriction map $\pi(\eta)=\eta|S$ (the image lies in $\bfk(S; [S]^{<\omega})$ because $S \in \cF$) and let $H=\{\pi(\eta_i): i \in I\}$. Note that to conclude it suffices to show that given any infinite subset $H\subseteq \bfk(S; [S]^{<\omega})$ we can always extract a further infinite countable subset $H'\subseteq H$ such that there is a linear functional $\xi : \bfk(S; [S]^{<\omega}) \to \bfk$ which does not vanish on any of the elements of $H'$.
        If $H$ spans a finite-dimensional subspace it suffices to set $H'=H \setminus \{0\}$ as there will be a linear functional $\xi$ with $\Kr (\xi) \cap \Span(H)=0$.
        If instead $H$ spans an infinite-dimensional subspace it suffices to extract an infinite linearly independent subset $H'$, complete $H'$ to a basis and consider a $\xi$ defined to be non-zero on every element of $H'$.
	\end{proof}
\end{proposition}
	
\begin{corollary}
	If $\cF=\cF^{\bot\bot}$ then $(f_i: i \in I)$ is summable if and only if $(\langle f_i, g \rangle: i \in I)$ is summable for every $g \in \bfk(\Gamma; \cF^\bot)$.
\end{corollary}
	
\begin{remark}
	The bornologies $\cF$ such that $\cF= \cF^{\bot\bot}$ are those corresponding to closed subsets of $\beta\Gamma \setminus \Gamma$ which are equal to the closure of their interior. In particular, if $\Gamma$ is infinite, maximal non trivial bornologies on $\Gamma$ (which correspond to points) never have this property: they are always such that $(\cF^{\bot\bot})^\sim=\cP(\Gamma)$. Though many bornologies used in the literature have it.
\end{remark}
	
\begin{example}\label{example:WO_bornology}
	A trivial example of a bornology with $\cF=\cF^{\bot\bot}$ is given by $\cF^{\sim}=[\Gamma]^{<\omega}$: then $(\cF^{\bot})^{\sim}=\cP(\Gamma)$ and $\cF=\cF^{\bot\bot}$. Less trivial examples are given as follows: let $\Gamma$ be a totally ordered set, define $\WO(\Gamma)$ to be the set of well ordered subsets, and $\WO_{\omega^2}(\Gamma)$ the set of subsets with order type $<\omega^2$.
	If $\cF^\sim=\WO_{\omega^2}(\Gamma)$, then $(\cF^\bot)^{\sim}=\WO(\Gamma^\op)$ and $(\cF^{\bot\bot})^\sim=\WO(\Gamma)$. Clearly any bornology $\cF$ with $\WO_{\omega^2}(\Gamma)\subseteq \cF^\sim \subset \WO(\Gamma)$ is then such that $(\cF^{\bot\bot})^\sim=\WO(\Gamma)$.
\end{example}
	
\begin{proposition}\label{prop:Sigma_pres_dual}
	If $F : \bfk(\Gamma; \cF) \to \bfk(Y; \cG)$ is $\Sigma$-preserving, then its dual $F^* : \bfk(\Delta; \cG)^* \to \bfk(\Gamma; \cF)^*$ maps $\Sigma$-preserving functionals to $\Sigma$-preserving functionals and the restriction
	$F^\bot : \bfk(\Delta; \cG^\bot) \to \bfk(\Gamma; \cF^\bot)$, characterized by $\langle Ff, g\rangle=\langle f, F^\bot g\rangle$, is $\Sigma$-preserving.
	\begin{proof}
		Clearly if $\eta \in \bfk(Y; \cG)^*$ is $\Sigma$-preserving then so is $F^*\eta = \eta \cdot F$ hence $F^*$ induces a map $F^\bot : \bfk(\Delta;\cG^\bot) \to \bfk(\Gamma;\cF^\bot)$ such that $\langle Ff, \eta\rangle=\langle f, F^\bot \eta\rangle$.
		To prove $F^\bot$ is $\Sigma$-preserving, consider a summable family $(\eta_i : i \in I)$ in $\bfk(\Delta; \cG^\bot)$. Then for every $g \in \bfk(\Delta; \cG^{\bot\bot})\supseteq \bfk(\Delta; \cG)$ one has that by Proposition~\ref{prop:sum_imp}(1) $(\langle \eta_i, g \rangle: i \in I)$ is summable in $\bfk$, so in particular for every $f \in \bfk(\Gamma; \cF)$ we have that $\langle \eta_i, Ff\rangle=\langle F^\bot \eta_i, f\rangle$ is summable over $i \in I$ in $\bfk$. And we can conclude by Proposition~\ref{prop:sum_imp}(2).
	\end{proof}
\end{proposition}
	
\begin{corollary}\label{cor:Dual_Case}
	A linear map $G: \bfk(\Delta; \cG^\bot) \to \bfk(\Gamma; \cF^{\bot})$ is $\Sigma$-preserving if and only if it is of the form $G = F^\bot$ for some $\Sigma$-preserving $F: \bfk(\Gamma; \cF^{\bot \bot}) \to \bfk(\Delta; \cG^{\bot\bot})$.
	\begin{proof}
		If $G$ is $\Sigma$-preserving then $F= G^{\bot}$ is $\Sigma$-preserving and $G=F^\bot$. Conversely if $G=F^\bot$ then it is $\Sigma$-preserving.
	\end{proof}
\end{corollary}
	
\begin{corollary}\label{cor:Sigma_ext}
	If $F : \bfk(\Gamma; \cF) \to \bfk(Y; \cG)$ is $\Sigma$-preserving, then it extends uniquely to a $\Sigma$-preserving $F^{\bot\bot} : \bfk(\Gamma; \cF^{\bot\bot}) \to \bfk(Y; \cG^{\bot\bot})$.
\end{corollary}
	
\begin{corollary}\label{cor:Vectop_equi}
	The full subcategory of $B\Sigma\Vect$ consisting of objects of the form $\bfk(I; \cP(I))$ is equivalent to $\Vect^\op$, via an extension of 
	\[\bfk^{\oplus I} \overset{f}{\longrightarrow} \bfk^{\oplus J} \qquad \mapsto \qquad \bfk(J; \cP(J)) \overset{f^*}{\longleftarrow} \bfk(I; \cP(I))\]
	\begin{proof}
		This follows from Corollary~\ref{cor:Dual_Case} observing that $\bfk^{\oplus I}= \bfk(I; [I]^{<\omega})$ and that a $\Sigma$-preserving map from $\bfk^{\oplus I}$ to $\bfk^{\oplus J}$ is just a $\bfk$-linear map.
	\end{proof}   
\end{corollary}
	
If $\bfk$ is a field and $V \in \Vect$, every basis $\Gamma$ of $V$ induces isomorphisms of vector spaces $V \simeq \bfk(\Gamma; [\Gamma]^{<\omega})$, and $V^* \simeq \bfk(\Gamma; \cP(\Gamma))$. By Corollary~\ref{cor:Dual_Case}, the notion of infinite sum we gain on $V^*$ does not depend on the choice of the basis $\Gamma$. We now see it coincides with the one we gain considering the inclusion $V^* \subseteq \bfk^V = \bfk(V,\cP(V))$.
	
\begin{lemma}\label{lem:sum_on_dual}
	If $V$ is a $\bfk$-vector space and $\Gamma$ is a basis of $V$ then the resulting isomorphism $\varphi : \bfk(\Gamma; \cP(\Gamma)) \to V^*$ when composed with the inclusion $V^* \subseteq \bfk(V; \cP(V))= \bfk^V$ is a $\Sigma$-closed embedding.
	\begin{proof}
		Let $(f_i : i \in I)$ be a summable family in $\bfk(\Gamma; \cP(\Gamma))$. For every $v\in \bfk(\Gamma; [\Gamma]^{<\omega})$, let $\tilde{v} := \sum_{\gamma \in \Supp(v)} v(\gamma)\gamma \in V$ denote its image under the natural isomorphism $\bfk(\Gamma; [\Gamma]^{<\omega})\simeq V$.
		By definition of $\varphi$, one has $\varphi(f_i)(v)= \langle v, f_i\rangle$. Now the right hand side is summable in $\bfk$ over $i \in I$ by hypothesis. This, by Proposition~\ref{prop:sum_imp}, entails that $(\varphi(f_i): i \in I)$ is summable in $\bfk^V$. The converse is clear from the same argument.
	\end{proof}
\end{lemma}		

\begin{remark}
    Given a bornology $\cF$ on some set $\Gamma$, Lemma~\ref{lem:sum_on_dual} and Proposition~\ref{prop:sum_imp} together imply that the map $\varphi: \bfk(\Gamma; \cF^{\bot}) \to \bfk(\Gamma; \cF)^{*}$ given by $\xi \mapsto \langle -, \xi\rangle$ is a $\Sigma$-closed embedding if we regard $\bfk(\Gamma; \cF)^{*}$ with its canonical $\Sigma$-structure. Thus all $\bfk(\Gamma; \cF^\bot)$ may be realized as $\Sigma$-closed subspaces of some $V^*$. Also, we know from Corollary~\ref{cor:Sigma_ext} that $\Sigma$-preserving maps $\bfk(\Gamma; \cF) \to \bfk(Y; \cG)$ are restrictions of $\Sigma$-preserving maps $\bfk(\Gamma; \cF^{\bot\bot}) \to \bfk(Y; \cG^{\bot\bot})$.
\end{remark}

\section{General strong vector spaces}\label{sec:GenSL}
	
This section is devoted to the discussion of the possible generalizations (i.e.\ fully faithful extensions) of $B\Sigma\Vect$ which are still to be regarded as categories of strong vector spaces. As premise we recall some notational conventions and standard definitions, and the foundation we will work in.

\smallskip

If $V$ is a $\bfk$-vector space and $B$ is a basis of $V$, then for every $b \in B$ we denote by $\delta_b^B: V \to \bfk$ the linear functional defined on $B$ as $\delta_b^B(b')=\delta_{b,b'}$ (Kronecker $\delta$), that is $(\delta_b^B)_{b \in B}$ is the dual family of functionals for the basis $B$ of $V$.
If $V$ is given as $V=\bfk^{\oplus I}$ then we denote the canonical basis as $e_I=\{e_i: i \in I\}$ and write $\delta_i$ for $\delta_{e_i}^{e_I}: \bfk^{\oplus I} \to \bfk$.

\smallskip
	
\noindent
{\bf Convention:} we write vertical compositions of natural transformations with $\cdot$ (this includes composition of arrows), we write horizontal composition of natural transformations (including composition of functors or of functors and natural transformations) as $\circ$ or by juxtaposition. For $V \in \Vect$ we identify $v \in V$ with the linear map $v: \bfk \to V$ mapping $1 \in \bfk$ to $v$.

\smallskip
	
\noindent
\textbf{Foundations:} we work in standard ZFC + existence of strongly inaccessible cardinals. We fix a non-trivial Grothendieck universe $\mathcal{U}$ and call every set in $\mathcal{U}$ \emph{small}. A category will be called \emph{locally small} if each of its hom-sets is in $\mathcal{U}$ and \emph{small} if it is locally small and its class of objects is small. The word \emph{large} will just mean non-small. We will write $\Sets$ to denote the category of small sets. Also recall that we use $\Vect$ to denote the category of \emph{small} $\bfk$-vector spaces. The word \emph{class} will just be a synonym for set (not necessarily small). 

\smallskip

\noindent
\textbf{$\bfk$-Enriched categories and functors:} recall that a \emph{$\bfk$-enriched category} is a category $\cC$ where all the hom-sets $\cC(c,c')$ are (possibly large) $\bfk$-vector spaces and the composition is bilinear. A functor between $\bfk$-enriched categories $F: \cC \to \cD$ is said to be \emph{$\bfk$-additive} if each component $\cC(c, c') \to \cD(Fc,Fc')$ is $\bfk$-linear.
When $\cC$ and $\cD$ are $\bfk$-enriched, we will denote the category of $\bfk$-additive functors from $\cC$ to $\cD$ and natural transformations among them by $[\cC, \cD]$.

\smallskip

\noindent
\textbf{$\bfk$-Enriched structure on additive categories:} Recall that if a $\bfk$-enriched category $\cC$ has biproducts when seen as a classic category, then it is said to be \emph{$\bfk$-additive}. For such a category, the $\bfk$-enriched structure can be recovered from the biproducts and the ring homomorphism $L_{\cC}:\bfk \to \Nat (\id_\cC, \id_\cC)$. Moreover any functor $F: \cC \to \cD$ between $\bfk$-additive categories is $\bfk$-additive if and only if it preserves biproducts and $F(L_{\cC}k)=(L_{\cD}k)F$ for every $k\in \bfk$ (see e.g.\ \cite[Sec.~8.5]{kashiwara2005categories}). So $[\cC, \cD]$ can be identified with the full subcategory of functors which preserve biproducts and respect the ring homomorphisms $L_{\cC}, L_{\cD}$ in the sense above. If $\cD$ is complete (resp.\ cocomplete), then limits (resp.\ colimits) in $[\cC, \cD]$ are computed pointwise, that is: the inclusion $[\cC,\cD] \subseteq \Func(\cC, \cD)$ of the category of $\bfk$-additive functors in the category of \emph{all} functors from $\cC$ to $\cD$ creates limits (resp.\ colimits).
	
\subsection{Reasonable categories of strong vector spaces}\label{ssec:RCSV}
We now discuss some properties which a reasonable category $\cC$ of strong vector spaces extending the category of based strong vector spaces of the form $\bfk(I; \cP(I))$ and $\Sigma$-preserving maps among them, should have, and postulate them as the definition of \emph{reasonable category of strong vector spaces}.

\begin{remark}[Heuristics]\label{rmk:Eur_0}
	The first assumption we make on a reasonable category $\cC$ of strong vector spaces is that each of its objects $X$ should be regarded as a small vector space $\mathfrak{U}X$ with some additional structure $\mathfrak{S}_X$, and maps $\cC(X, Y)$ should be maps in $\Vect(\mathfrak{U}X, \mathfrak{U}Y)$ which preserve such additional structure. Moreover we require that a finite linear combination of strongly linear maps is still strongly linear. It is easily understood that this means that we are requiring $\cC$ to be a $\bfk$-additive category admitting a $\bfk$-additive forgetful functor to $\mathfrak{U}: \cC \to\Vect$. The additional structure $\mathfrak{S}_X$ should carry the following information: for each set $I$ a subspace of $(\mathfrak{U}X)^{I}$ of families $(x_i)_{i \in I} \in (\mathfrak{U}X)^I$ which will be declared to be \emph{summable} and  for each such family, for every $(k_i)_{i \in I}\in \bfk^I$ the value of the infinite linear combination $\sum_{i \in I} x_i \cdot k_i$.
		
	We saw that the subcategory of $B\Sigma\Vect$ consisting of spaces of the form $\bfk(I; \cP(I))$ and strongly linear maps is equivalent to $\Vect^\op$, via a functor sending $\bfk(I; \cP(I))$ to $\bfk^{\oplus I}=\bfk(I; [I]^{<\omega})$ and $f$ to $f^\bot$ (Corollary~\ref{cor:Dual_Case}).
	Therefore it makes sense that a category of strong vector spaces comes as full extension of $\Vect^\op$ i.e.\ as a $\bfk$-additive fully faithful functor $\iota: \Vect^\op \to \cC$.
	Now we assume that:
	\begin{enumerate}[label=(H\arabic*), ref=(H\arabic*)]
		\item\label{SV:H1} the faithful underlying-vector-space functor $\mathfrak{U}$ is the representable $\cC(\iota\bfk, -)$, i.e.\ $\iota\bfk$ is a \emph{separator}\index[ind]{separator};
		\item\label{SV:H2} a family $(x_i)_{i \in I} \in (\mathfrak{U}X)^I$ is \emph{summable} if and only if there is an arrow $\overline{x}: \iota(\bfk^{\oplus I}) \to X$ such that $(\mathfrak{U} \overline{x}) (\delta_i^I) = x_i$ for every $i$, where $\delta_i^I \in \cC(\iota\bfk, \iota(\bfk^{\oplus I}))\cong \Vect(\bfk^{\oplus I}, \bfk)\cong \bfk^I$; when this is the case the infinite linear combinations are given by
		\[\sum_{i \in I} k_i x_i = (\mathfrak{U}\overline{x}) (\xi) \quad \text{where}\quad \forall i \in I,\, \xi \cdot e_i = k_i\]
		\item\label{SV:H3} every arrow $f: \iota(\bfk^{\oplus I}) \to X$ is uniquely determined by the family $(x_i)_{i \in I} \in (\mathfrak{U}X)^I$ given by $(\mathfrak{U} f)(\delta_i^I) = x_i$, i.e.\ by the family of compositions
		\[\begin{tikzcd}
            \iota\bfk \arrow[r,"\iota\delta_i"] &\iota(\bfk^{\oplus I}) \arrow[r, "f"] &X
        \end{tikzcd}\]
        as $i$ ranges through $I$;
		\item\label{SV:H4} a map $\hat{f} \in \Vect(\mathfrak{U} X, \mathfrak{U}Y)$ has the form $\hat{f}=\mathfrak{U}f$ for a (unique by faithfulness of $\mathfrak{U}$) $f \in \cC(X, Y)$ if and only if it is \emph{strongly linear} in the sense that for every set $I$, for every summable $(x_i)_{i \in I} \in (\mathfrak{U}X)^I$ and for every $(k_i)_{i \in I}\in \bfk^I$ it happens that $(\hat{f}x_i)_{i \in I}$ is summable and $\hat{f} \sum_{i \in I} k_i x_i= \sum_{i \in I} k_i \hat{f} x_i$.
	\end{enumerate}
	In view of \ref{SV:H2} and \ref{SV:H3} the information for $\mathfrak{S}_X$ can then be redundantly encoded in the \emph{canonical diagram} (cf \cite[p.~2]{adamek1994locally}) 
	\[\mathfrak{S}_X: (\iota \downarrow X) \to \cC \qquad (f_V: \iota V \to X) \mapsto \iota V.\]
	Given the assumptions \ref{SV:H2} and \ref{SV:H3} on the meaning of the notions of \emph{summable} family and infinite linear combination we get that for $\hat{f} \in \Vect(\mathfrak{U} X, \mathfrak{U}Y)$ to be \emph{strongly linear} is equivalent to the fact that for every $V \in \Vect^\op$ and every $h \in \cC(\iota V, X)$ there is a (unique by faithfulness of $\mathfrak{U}$) $h' \in \cC(\iota V ,Y)$ such that
	\[\label{eq1}
	\begin{tikzcd}
	    &\mathfrak{U}(\iota V)\ar[dl, "\mathfrak{U}h"] \ar[dr,"\mathfrak{U}h'"']\\
		\mathfrak{U}X \ar[rr,"\hat{f}"] &&\mathfrak{U}Y
	\end{tikzcd}\]
	commutes.
	Now we see that \ref{SV:H4} implies that each $X \in \cC$ is the colimit of the diagram of $\mathfrak{S}_X : (\iota \downarrow X) \to \cC$ with the obvious cone $h: \mathfrak{S}_X(V, h) \to X$. Indeed if $\varphi_{\bullet}: \mathfrak{S}_X \to \Delta Y$ is a cone, then it defines a map $\hat{f}:\mathfrak{U}X \to \mathfrak{U}Y$ by $\hat{f} (x) = \varphi_{(\bfk, x)} \in \cC(\iota\bfk, Y)$. We see that $\hat{f}$ is strongly linear: for any $(V, h) \in (\iota \downarrow X)$ we indeed have
	$\hat{f} \cdot (\mathfrak{U} h) = \mathfrak{U}\varphi_{V,h}$. So by assumption \ref{SV:H4}, $\hat{f}=\mathfrak{U}f$ for $f \in \cC(X, Y)$ and by faithfulness of $\mathfrak{U}$ we have $f\cdot h =\varphi_{(V, h)}$ for every $(V, h) \in (\iota \downarrow X)$.
	Uniqueness easily follows because any $f$ such that $f\cdot h =\varphi_{(V, h)}$ for every $h$ should satisfy $\mathfrak{U}f=\hat{f}$.
\end{remark}
	
Motivated by the above Remark we give the following definition:
	
\begin{definition}\label{defn:reasonable}
    A locally small $\bfk$-additive fully faithful extension $\iota: \Vect^\op \to \cC$ of $\Vect^\op$ is a \emph{reasonable category of strong vector spaces}\index[ind]{strong vector space!reasonable category of ---s|textbf} if the following hold:
	\begin{enumerate}[label=(A\arabic*), ref=(A\arabic*), align = left, leftmargin =*, itemindent = !]
		\item\label{SV:A1} the functor $\iota: \Vect^\op \to \cC$ is \emph{dense} (cf \cite[Ch.~X, Sec.~6]{maclane71categories}), that is, every object of $X \in \cC$ should be the colimit of the canonical (possibly large) diagram
	    \[(\iota \downarrow X) \to \cC \qquad (f_V: \iota V \to X) \mapsto \iota V;\]
		\item\label{SV:A2} $\iota \bfk$, where $\bfk \in \Vect^{\op}$ is the one dimensional space, is a separator (i.e.\ $\cC(\iota\bfk, -): \cC \to \Vect$ is faithful);
		\item\label{SV:A3} every map $f: \iota H \to X$ for $H \in \Vect^{\op}$ is uniquely determined by its value on the points $\delta^{B}_b \in \Vect^{\op}(\bfk , V)$ as $b$ ranges through a basis $B$ of $V$.
	\end{enumerate}
\end{definition}
	
\begin{remark}
	Remark~\ref{rmk:Eur_0} shows that \ref{SV:A1}-\ref{SV:A3} of Definition~\ref{defn:reasonable} are a consequence of the assumption \ref{SV:H1}-\ref{SV:H4} we made. It will be shown that they are equivalent in Lemma~\ref{lem:sum-pres_characterization}.
\end{remark}
	
\begin{remark}
	It will turn out (Remarks~\ref{rmk:Q_2,3,4_rel} and \ref{rmk:A2_redundant}) that \ref{SV:A2} is actually redundant as it follows from \ref{SV:A1} and \ref{SV:A3}. 
\end{remark}

\smallskip
	
We shall prove (Theorem~\ref{thm:SigmaVectChar}) that there is a universal reasonable category of strong vector spaces: that is, a unique $\bfk$-additive fully faithful extension $\iota: \Vect^\op \to \Sigma\Vect$ satisfying all the above constraints and moreover
\begin{enumerate}[label=(U), ref=(U)]
	\item\label{SV:U} for every reasonable category $\iota': \Vect^\op \to \cC$, there is a unique up to a unique isomorphism fully faithful functor $\eta: \cC \to \Sigma\Vect$ such that $\iota \cong \eta \circ \iota'$.
\end{enumerate}

\subsection{Consequences of \texorpdfstring{\ref{SV:A1}}{(A1)} and \texorpdfstring{\ref{SV:A2}}{(A2)}}\label{ssec:A12}
In this subsection we recall some facts in category theory and draw some consequences of \ref{SV:A1} and \ref{SV:A2}.
	
\begin{remark}\label{rmk:A1^*}
	It follows from the general theory of dense functors that if \ref{SV:A1} holds, then $\cC$ has to be equivalent to a full subcategory of $[\Vect, \Vect]$ containing the essential image of the Yoneda embedding
	\[\Yon: \Vect^{\op} \to [\Vect, \Vect], \qquad \Yon V \cong \Vect^\op(-, V) = \Vect(V,-).\]
	More precisely: by \cite[Ch.~X Sec.~6 Prop.~2]{maclane71categories} if $\iota : \Vect^\op \to \cC$ is dense, then $U_* \circ \cC(\iota, -) : \cC \to \Func(\Vect, \Sets)$ is fully faithful where $U_* : [\Vect, \Vect] \to \Func(\Vect, \Sets)$ is the postcomposition with the forgetful $U: \Vect \to \Sets$. Since $U_*$ is faithful, this entails that $\cC(\iota, -)$ must be fully faithful as well.
	Thus, since we are interested in determining $\cC$ just up to equivalence, in view of \ref{SV:A1} we may replace $\cC$ with the essential image of the fully faithful $\cC(\iota, -)$ and make the following stronger assumption
	\begin{enumerate}[label=(A1$^*$), ref=(A1$^*$), align = left, leftmargin = *]
		\item\label{SV:A1*} $\cC$ is a full subcategory of $[\Vect, \Vect]$ and 
		\[\iota =\Yon : \Vect^{\op} \to [\Vect, \Vect].\]
	\end{enumerate}
	In particular $\iota \bfk = \Vect(\bfk, -)$ is isomorphic to the identity functor on $\Vect$ and $\iota (\bfk^{\oplus I}) = \Vect(\bfk^{\oplus I},-)\cong \Vect(\bfk, -)^I$ is isomorphic to the $I$-th power functor.
	Moreover by Yoneda's Lemma
	\[[\Vect, \Vect] (\Vect(\bfk^{\oplus I},-), X) \cong X(\bfk^{\oplus I}),\]
    so, under assumption \ref{SV:A1*}, condition \ref{SV:A2} instantiates to
	\begin{enumerate}[label=(A2$^*$), ref=(A2$^*$), resume, align = left, leftmargin = *]
		\item\label{SV:A2*} $-(\bfk): \cC \to \Vect$ is faithful.
	\end{enumerate}
\end{remark}

Recall the following folklore result entailing that the functor $-(\bfk): [\Vect, \Vect] \to \Vect$ has a left adjoint.
	
\begin{lemma}\label{lem:Fin_adj}
	The functor
	\[-(W) \cong \Nat (\Yon W, -): [\Vect, \Vect] \to \Vect\] has a left adjoint
	$\Fin^W: \Vect \to [\Vect,\Vect]$ given by 
	\[\Fin^W V \cong  V \otimes \Vect(W,-).\]
	\begin{proof}
		By Yoneda's Lemma there is an isomorphism 
		\[\Vect(V, XW) 
        \cong \Nat (V \otimes \Vect(W,-), X)\]
		sending $f: V \to XW$ to the natural transformation
		\[V \otimes \Vect(W, H) \to XH \qquad (v \otimes h)\mapsto Xh \cdot f \cdot v.\]
		The map $f$ can be recovered from such natural transformation considering $H=W$ and $h=1_W: W \to W$.
	\end{proof}
\end{lemma}
	
\noindent
{\bf Notation}: we will denote the left adjoint $\Fin^\bfk$ of $-(\bfk)$ just by $\Fin$.
\smallskip
	
When $\cC$ satisfies \ref{SV:A1*}, it has to contain all representables. If it furthermore satisfies \ref{SV:A2*}, then  every $X \in \cC$ must have the property that $\Nat(\Yon \bfk, -)$ acts as a monomorphism on $\Nat(\Yon V, X)=\cC(\Yon V, X)$. This is readily observed to have some equivalent formulations which will be useful.
	
\begin{lemma}\label{lem:simple_char}
	Let $X \in [\Vect, \Vect]$. The following are equivalent
	\begin{enumerate}
		\item for every pair of maps $f, g : \Vect(V, -) \to X$ we have $f=g$ if and only if $f\bfk= g\bfk: \Vect(V, \bfk) \to X\bfk$;
		\item $\Nat(-, X)$ sends the canonical map $V^* \otimes - \to \Vect(V,-)$ to a monomorphism;
		\item $\Nat(\Yon V/ \Fin (V^*), X)=0$.
	\end{enumerate}
	\begin{proof}
        (1) states that $-(\bfk) : \Nat(\Yon V, X) \to \Vect(V^* , X\bfk)$ is mono. By Lemma~\ref{lem:Fin_adj}
        \[\Vect(V^*, X\bfk) \cong \Nat(\Fin(V^*), X)=\Nat(V^* \otimes -, X).
		\]
        The resulting composite map $\Nat(\Yon V, X) \to \Nat(V^* \otimes -, X)$ is precisely $\Nat(i,X)$ where $i: \Fin (V^*) \to \Yon V$ is the natural inclusion. Thus (1) is equivalent to (2).

		Since $\Nat(-,-): [\Vect, \Vect]\times [\Vect,\Vect]\to \Vect$ is a left exact $\bfk$-additive bifunctor and $[\Vect, \Vect]$ has cokernels, $(2)$ is in turn equivalent to $(3)$.
	\end{proof}
\end{lemma}

\begin{definition}
	A $\bfk$-additive functor $X: \Vect \to \Vect$ will be said to be \emph{$\bfk$-concrete} if one of the equivalent conditions above holds for every small vector space $V$.
\end{definition}

\begin{lemma}\label{lem:simple_prop}
    If $X$ is $\bfk$-concrete then for every $Z \in [\Vect,\Vect]$ the map $-(\bfk) : \Nat(Z,X) \to \Vect(Z\bfk, X\bfk)$ is one-to-one.
    \begin{proof}
        Suppose $f,g: Z \to X$ are such that $f \neq g$. Then there is $V \in \Vect$ and $z \in ZV$ such that $(fV)\cdot z \neq (gV) \cdot z$. But then considering $\dot{z}:=(Z-)\cdot z : \Yon V \to Z$, we have $f \cdot \dot{z} \neq g \cdot \dot{z}$ and thus by simplicity of $X$, $(f\cdot \dot{z})(\bfk) \neq (g\cdot \dot{z})(\bfk) : V^* \to X\bfk$, whence $f\bfk\neq g \bfk$.
    \end{proof}
\end{lemma}

\begin{example}
	All representables $\Yon V$ and all objects of the form $\Fin V$ are simple.
\end{example}
	
\begin{remark}\label{rmk:cN_2}
	The expression $\Nat(\Yon V / \Fin(V^*), X)$ is functorial in $V$ and defines a subfunctor of $X \cong \Nat( \Yon -, X)$ via Yoneda. Such subfunctor will be denoted by $\cN_{2} X$ and has the following explicit description
	\[\Nat(\Yon V / \Fin (V^*), X)\cong (\cN_2 X) V \coloneqq \bigcap_{\xi \in V^*} \Kr X \xi \subseteq XV.\]
\end{remark}
	
\begin{theorem}\label{thm:simpcat}
	The full subcategory on the class of objects $\{X: \cN_2 X=0\}$ of $\bfk$-concrete functors is the largest full subcategory of $[\Vect, \Vect]$ containing all representables and such that $\Yon \bfk=\Vect(\bfk,-)$ is a separator. 
	\begin{proof}
		By Lemma~\ref{lem:simple_char}, if a full subcategory of $[\Vect, \Vect]$ contains all representables and has $\Yon \bfk$ as a separator, then all of its objects $X$ must be such that $\cN_2 X=0$. On the other hand the representables are a family of separators for $[\Vect, \Vect]$ and hence, a fortiori, they are a family of separators for any subcategory containing them, in particular for full subcategory on $\{X: \cN_2 X=0\}$. Since in turn $\Vect(\bfk, -)$ will be a separator for (the full subcategory on) the representables, it will be a separator for the whole full subcategory on $\{X: \cN_2 X=0\}$. 
	\end{proof}
\end{theorem}
	
\subsection{\texorpdfstring{\ref{SV:A3}}{(A3)} and uniqueness of infinite sums}\label{ssec:A3}

In this subsection we will prove an analogue of Theorem~\ref{thm:simpcat} for condition \ref{SV:A3}.
In Remark~\ref{rmk:Eur_0}, assumption \ref{SV:H2} stipulated that a family $(x_i)_{i \in I}$ in $\cC(\iota\bfk, X)$ has to be deemed summable if and only if there is an arrow $\emph{x}: \iota (\bfk^{\oplus I}) \to \bfk$ such that $x \cdot \iota(\delta_{i}^I) = x_i$. In view of assumption \ref{SV:A1*} that $\iota : \Vect^\op \to \cC$ is a full subcategory of $\Yon : \Vect^\op \to [\Vect, \Vect]$ we turn it into the following definition.
	
\begin{definition}\label{defn:summable2}
	For $X \in [\Vect, \Vect]$, a family $(x_i)_{i \in I}\in (X\bfk)^I$ is \emph{summable}\index[ind]{summable|textbf} if there is an $x \in X(\bfk^{\oplus I})$ such that $x_i= (X\delta_i^I) \cdot x$ for all $i \in I$, or equivalently if it is in the image of the map
	\[\begin{tikzcd}
	    X(\bfk^{\oplus I}) \ar[rr, "(X\delta_i)_{i \in I}"] &&(X\bfk)^I.
	\end{tikzcd}\]
	Denoting by $\sigma : \bfk^{\oplus I} \to \bfk$ the linear functional defined by $\sigma \cdot e_i = 1$ for every $i \in I$, we call a \emph{sum}\index[ind]{sum|textbf} for a summable $(x_i)_{i \in I}$ any element of the form $(X\sigma) \cdot x$ for some $x\in X(\bfk^{\oplus I})$ such that $x_i= (X\delta_i^I) \cdot x$ for every $i \in I$.
\end{definition}

\begin{example}
    Let $X= \Yon(\bfk^{\oplus J})\cong (-)^J$, so $X(\bfk)=\Vect(\bfk^{\oplus J}, \bfk) \cong \bfk^J$, and let $I \subseteq J$. For every $I$-tuple $(k_i)_{i \in I}$ of elements of $\bfk$ the family $(\delta_i^J\cdot k_i)_{i \in I}$ is indeed summable in $X\bfk$, as it has the form $\big((X \delta_i^I )\cdot x\big)_{i \in I}$ where $x \in X(\bfk^{\oplus I})=\Vect(\bfk^{\oplus J}, \bfk^{\oplus I})$ is given by $x(e_j^J) = e_j^I \cdot k_i$ if $j \in I$ and $x(e_j^J)=0$ otherwise. Now the element $(X\sigma) \cdot x \in X(\bfk) = \Vect(\bfk^{\oplus J}, \bfk)$ is such that $((X\sigma) \cdot x) \cdot e_j=k_j$ if $j \in I$ and $0$ otherwise, that is it corresponds via $X(\bfk)\cong \bfk^{J}$ to the tuple $(k_j')_{j \in J}$ where $k_j=k_j'$ if $j' \in I$ and $k_j'=0$ otherwise.
\end{example}
    
As defined above sums are not necessarily unique, but they are if we restrict to certain objects.
	
\begin{remark}
    Under assumption \ref{SV:A1*}, assumption \ref{SV:A3} becomes
	\begin{enumerate}[label=(A3$^*$), ref=(A3$^*$), align = left, leftmargin = *, itemindent = !]
		\item\label{SV:A3*} for every $X$ in $\cC$ and every $(x_i)_{i \in I} \in (X\bfk)^I$ there is at most one $x \in X\bfk^{\oplus I}$ such that every $i \in I$, $X\delta_i^I \cdot x = x_i$. 
	\end{enumerate}
	So in particular if $\cC$ satisfies \ref{SV:A1*} and \ref{SV:A3*}, then for every $X\in \cC$ every summable family has a unique sum.
\end{remark}	
	
\begin{lemma}\label{lem:k_sums}
	Let $X: \Vect \to \Vect$ be $\bfk$-additive and $\kappa$ be a cardinal. Then the following statements are equivalent
	\begin{enumerate}
		\item every $\kappa$-indexed summable family in $X(\bfk)$ has a unique sum;
		\item one has $\bigcap_{i \in I}\Kr(X\delta_i) \subseteq \Kr(X\sigma)$ in $X(\bfk^{\oplus I})$ with $\delta_i$ and $\sigma$ as in Definition~\ref{defn:summable2} and $|I|= \kappa$;
		\item for every set $I$ of size $|I|=\kappa$ one has
		\[\bigcap_{i \in I} \Kr(X\delta_i) = \bigcap_{\xi \in (\bfk^{\oplus I})^*} \Kr (X\xi)= \cN_2 X(\bfk^{\oplus I}).\]
	\end{enumerate}
	\begin{proof}
		The only non trivial implication is $(2) \Rightarrow (3)$. To show it, note that if $\xi \in \Vect(\bfk^{\oplus I}, \bfk)$, and one defines $f: \bfk^{\oplus I} \to \bfk^{\oplus I}$ by $f \cdot e_i =e_i \cdot (\xi \cdot e_i)$, then $\Kr(\delta_i) \subseteq \Kr((\xi \cdot e_i)\cdot \delta_i)=\Kr(\delta_i\cdot f)$. Thus since $X$ is left exact $\Kr(X\delta_i) \subseteq \Kr (X(\delta_i \cdot f))$. On the other hand $\sigma \cdot f = \xi$ so again by left exactness $\Kr(X\xi)=\Kr (X(\sigma \cdot f)) \supseteq \bigcap_{i \in I} \Kr (X (\delta_i \cdot f)) \supseteq \bigcap_{i \in I} \Kr (X\delta_i)$.
	\end{proof}
\end{lemma}
	
\begin{definition}\label{def:strongly_lin}
	Let $X \in \Ind\-(\Vect^{\op})$. Say that
	\begin{enumerate}
		\item $X$ \emph{has unique $\kappa$-sums} if the equivalent conditions of Lemma~\ref{lem:k_sums} hold;
		\item $X$ \emph{has unique infinite sums} if it has unique $\kappa$-sums for every small cardinal $\kappa$;
		\item $X$ is a \emph{strong $\bfk$-vector space} if it is $\bfk$-concrete and has unique infinite sums.
	\end{enumerate}
	If $X$ has $\kappa$-sums, for a summable family $(x_i: i \in I) \in (X\bfk)^I$, $|I|\le \kappa$, we denote its unique sum as $\sum_{i \in I} x_i$.\\
	If $X$ and $Y$ have $\kappa$-sums we say that a linear map $\tilde{f}: X\bfk \to Y\bfk$ \emph{preserves $\kappa$-sums} if and only if, for every summable family $(x_i: i \in I) \in X(\bfk)^I$ with $|I|\le \kappa$, we have that $(\tilde{f} \cdot x_i: i \in I) \in Y(\bfk)^I$ is summable and $\tilde{f} \cdot \sum_{i \in I} x_i= \sum_{i \in I} \tilde{f} \cdot x_i$.\\
	If $X$ and $Y$ have infinite sums and $\tilde{f}: X\bfk \to Y\bfk$ preserves $\kappa$-sums for every $\kappa$, then we say that $\tilde{f}$ \emph{preserves infinite sums}.
\end{definition}
	
Similarly to Lemma~\ref{lem:simple_char} we note
	
\begin{lemma}\label{lem:SigmaVect_char}
	For $X \in [\Vect, \Vect]$ the following are equivalent:
	\begin{enumerate}
		\item $X$ is a strong $\bfk$-vector space;
		\item for every set $I$ one has
		$\bigcap_{i \in I} \Kr (X \delta^I_i)= 0$;
		\item every $x \in X(\bfk^{\oplus I})$ is uniquely determined by the family $(X\delta_i^I \cdot x)_{i \in I} \in (X\bfk)^I$, that is the map $(X\delta_i)_{i \in I} : X(\bfk^{\oplus I})  \to (X\bfk)^I$ is one-to-one;			\item for every set $I$, $\Nat(-, X)$ sends the natural inclusion $\delta_I: (\Yon \bfk)^{\oplus I} \to (\Yon \bfk)^I$ to a monomorphism in $\Vect$;
		\item $\Nat((\Yon \bfk)^I/(\Yon \bfk)^{\oplus I}, X)=0$ for every set $I$.
	\end{enumerate}
	\begin{proof}
		$(1)\Leftrightarrow(2)$ follows from the equivalence $(1)\Leftrightarrow (3)$ of Lemma~\ref{lem:k_sums}. $(2)\Leftrightarrow (3)$ is clear. The equivalence of $(3)$ with $(4)$ follows from the observation that the map $(X\delta_i)_{i \in I}$ is precisely the composition
		\[\begin{tikzcd}
		    X(\bfk^{\oplus I}) \cong \Nat((\Yon \bfk)^I, X)
            \arrow[rr, "{\Nat(\delta_I,X)}"]
            &&\Nat((\Yon \bfk)^{\oplus I}, X)\cong (X\bfk)^I,
	    \end{tikzcd}\]
		where the first and last isomorphisms are the obvious ones. The equivalence of $(4)$ and $(5)$ is proven like the equivalence of $(2)$ and $(3)$ in Lemma~\ref{lem:simple_char}.
	\end{proof}
\end{lemma}
	
\begin{remark}\label{rmk:A2_redundant}
	From the equivalence $(1)\Leftrightarrow (3)$ of the previous Lemma it follows that if $\cC$ satisfies \ref{SV:A1*} and \ref{SV:A3*}, then every $X\in \cC$ has to be a strong vector space, so in particular it has to be simple, therefore \ref{SV:A1*} and \ref{SV:A3*} together imply \ref{SV:A2*}. Thus \ref{SV:A1} and \ref{SV:A3} imply \ref{SV:A2}. It follows that the full subcategory of $[\Vect, \Vect]$ of strong vector spaces is the largest (i.e.\ the maximum) full subcategory of $[\Vect, \Vect]$ satisfying \ref{SV:A2*} and \ref{SV:A3*} and by Remark~\ref{rmk:A1^*} it therefore satisfies \ref{SV:U}.
\end{remark}
	
For future reference we summarize the discussion above in a definition and a theorem.
	
\begin{definition}
	Denote the full subcategory of $[\Vect, \Vect]$ with objects the strong vector spaces by $\Sigma\Vect$.
\end{definition}
	
\begin{theorem}\label{thm:SigmaVectChar}
	The following hold
	\begin{enumerate}
		\item $\Sigma\Vect$ is the biggest full subcategory of $[\Vect, \Vect]$ satisfying \ref{SV:A3*} and containing all representables. In particular every $\bfk$-additive locally small category satisfying \ref{SV:A1}, \ref{SV:A3} (and hence also \ref{SV:A2}) is equivalent to a full subcategory of $\Sigma\Vect$.
		\item $\Sigma\Vect$ contains all of the $\Fin V$ for $V \in \Vect$ and is the largest full subcategory containing $\Yon (\Vect^\op)$ and $\Fin(\Vect)$ for which all the arrows of the form $\delta_I: (\Yon\bfk)^{\oplus I} \to \Vect(\bfk^{\oplus I},-)\cong (\Yon\bfk)^{I}$ are epimorphisms.
	\end{enumerate} 
	\begin{proof}
		We already proved (1). So let us prove (2). Note that for each $V$, $\Fin V$ is a strong $\bfk$-vector space and that $(\Yon\bfk)^{\oplus I}\cong \Fin (\bfk^{\oplus I})$, whence in particular $(\Yon\bfk)^{\oplus I}$ for every set $I$.
        If the maps $\delta^I$ are epimorphisms in a full subcategory $\cC\subseteq [\Vect, \Vect]$ containing $\Yon (\Vect^\op)$ and $\Fin(\Vect)$, then every object $X\in \cC$ must be in $\Sigma\Vect$ by Lemma~\ref{lem:SigmaVect_char}~(4).
	\end{proof}
\end{theorem}
	
Let $\iota : \Vect^\op \to \cC$ be a $\bfk$-additive functor. In Remark~\ref{rmk:Eur_0} we saw that under assumptions \ref{SV:H1}, \ref{SV:H2}, and \ref{SV:H3} on the meaning of summability, assumption \ref{SV:H4} that strongly linear maps must coincide with arrows in $\cC$, implies that $\iota$ has to be a dense functor, hence leading to condition \ref{SV:A1} (and, up to equivalence, \ref{SV:A1*}). Now we are about to check a converse: i.e.\ that \ref{SV:A1} and \ref{SV:A3} imply that every $\tilde{f}: \cC(\iota\bfk,X) \to \cC(\iota\bfk, Y)$ preserving infinite sums  has the form $\cC(\iota\bfk, f)$ for a unique $f\in \cC(X, Y)$. In other words (instantiating to \ref{SV:A1*} and using Yoneda's Lemma): if $X,Y$ are strong $\bfk$-vector spaces, the infinite-sum-preserving maps $\tilde{f}: X\bfk \to Y\bfk$ are precisely those of the form $f\bfk$ for some $f: X\to Y$, and such $f$ is unique. 
	
\begin{lemma}\label{lem:sum-pres_characterization}
	For a map $\tilde{f}: X(\bfk) \to Y(\bfk)$ where $X$ and $Y$ are strong $\bfk$-vector spaces, and for $I$ a set with $|I|=\kappa$, the following are equivalent:
	\begin{enumerate}
		\item $\tilde{f}$ preserves $\kappa$-sums;
		\item there is a unique map $f_{I}$ making the following diagram commute
		\[\begin{tikzcd}
			X(\bfk^{\oplus I}) \ar[rr,"{(X-)\cdot -}"] \ar[d, "f_I"']
			&&\Vect (\Vect(\bfk^{\oplus I}, \bfk), X(\bfk)) \ar[d, "{\Vect(\Vect(\bfk^{\oplus I}, \bfk), \tilde{f})}"]\\
			Y(\bfk^{\oplus I}) \ar[rr,"{(Y-)\cdot -}"]
			&&\Vect (\Vect(\bfk^{\oplus I}, \bfk), Y(\bfk))
		\end{tikzcd}\]
		where $(X-)\cdot -$ maps $x$ to $(X-)\cdot x: \Vect(\bfk^{\oplus I}, \bfk) \to \bfk$.
		\item for every $V\in \Vect$ of dimension $\le \kappa$ there is a unique $f_V: XV \to YV$ such that the following diagram commutes
		\[\begin{tikzcd}
			V^* \otimes_\bfk XV \ar[rr, "{-\cdot (X-)}"] \ar[d,"{V^* \otimes f_V}" ']
			&&X\bfk \ar[d, "{\tilde{f}}"]\\
			V^* \otimes_\bfk YV \ar[rr, "{-\cdot (Y-)}"]
			&&Y\bfk
		\end{tikzcd}\]
	\end{enumerate}
	\begin{proof}
		Assume $(2)$ holds and $(x_i)_{i \in I}$ is a summable family. Since $X \in \Sigma\Vect$ there is a unique $x \in X(\bfk^{\oplus I})$ such that $(X\delta_i) \cdot x = x_i$ and for such $x$ we have $(X \sigma) \cdot x = \sum_{i \in I} x_i$. By $(2)$ for every $i\in I$ we have $\tilde{f}\cdot x_i = (Y \delta_i) \cdot (f_I \cdot x)$ so
		\[\tilde{f} \cdot \sum_{i \in I} x_i = \tilde{f} \cdot (X \sigma) \cdot x = (Y\sigma) \cdot (f_I \cdot x) = \sum_{i \in I} \tilde{f} \cdot  x_i,\]
		where the middle equality follows again from $(2)$.
		Conversely, to prove $(1)\Rightarrow (2)$ observe that the kernels of the horizontal arrows in the diagram are respectively $(\cN_2 X )(\bfk^{\oplus I})$ and $(\cN_2 Y)(\bfk^{\oplus I})$, so since $X$ and $Y$ are simple, the horizontal arrows in the diagram are mono. It follows that if a map $f_I$ making the diagram commute exists, then it is unique. Now assume $x \in X(\bfk^{\oplus I})$ and let us show $\tilde{f} \cdot (X-)\cdot x$ is of the form $(Y-) \cdot y$ for some $y \in Y(\bfk^{\oplus I})$. Since $X$ and $Y$ have unique $\kappa$-sums the $I$-tuple $(x_i)_{i \in I}:=((X\delta_i) \cdot h)_{i \in I}$ suffices to determine $(X-) \cdot x$; in particular if $\xi \in \Vect(\bfk^{\oplus I}, \bfk)$ we have $(X \xi) \cdot x = \sum_{i \in I} (\xi \cdot e_i) \cdot x_i$. Given that $\tilde{f}$ preserves summability, there is $y$ such that $(\tilde{f} \cdot x_i)_{i \in I} = ((Y\delta_i) \cdot y)_{i \in I}$. On the other hand, given any $\xi$, since $\tilde{f}$ preserves infinite sums, we have
		\[\begin{aligned}
		\tilde{f} \cdot ((X \xi) \cdot x) &= \tilde{f} \cdot \sum_{i \in I} x_i \cdot (\xi \cdot e_i) = \sum_{i \in I} (\tilde{f} \cdot x_i) \cdot (\xi \cdot e_i) =\\
		&= \sum_{i \in I} ((Y \delta_i) \cdot y) \cdot (\xi \cdot e_i) =(Y \xi) \cdot y.
		\end{aligned}\]
		So actually $(Y-) \cdot y = \tilde{f} \cdot ((X -)\cdot x)$ and we proved $(1) \Rightarrow (2)$. Finally observe that $(3)$ for $\dim(V)=\kappa$ is just a rewriting of $(2)$ and that maps preserving $\kappa$-indexed infinite sums also preserve $\kappa'$-indexed infinite sums for $\kappa'\le \kappa$.
	\end{proof}
\end{lemma}

\subsection{Strong vector spaces as structures in a signature}\label{ssec:AltSigmaVect}
	
We give an equivalent description of $\Sigma\Vect$ in terms of models of a theory in a (large) signature. We do this by defining $\Sigma\Vect'$ to be the category of models of such a theory, defining how to interpret every object of $\Sigma\Vect$ as a model, and finally proving that every object in $\Sigma\Vect'$ is isomorphic to some model coming from $\Sigma\Vect$. In the following we denote by $\Card$ the (large) set of small cardinals.
	
\begin{definition}
	Define the category $\Sigma\Vect'_\bfk$ as follows. The objects are tuples $(V, (\mathfrak{S}_\lambda)_{\lambda \in \Card}, (\Sigma_\lambda)_{\lambda \in \Card})$ such that 
	\begin{enumerate}[label=(B\arabic*), ref=(B\arabic*)]
    \setcounter{enumi}{-1}
		\item\label{SV:B0} $V$ is a small $\bfk$-vector space;
		\item\label{SV:B1} for every $\lambda\in \Card$, $\mathfrak{S}_\lambda$ is a subspace of $V^\lambda$ containing $V^{\oplus \lambda}$ and $\Sigma_\lambda: \mathfrak{S}_\lambda \to V$ is a $\bfk$-linear function;
		\item\label{SV:B2} if $\lambda, \mu, \nu \in \Card$, $(v_i)_{i<\lambda} \in \mathfrak{S}_\lambda$, $(s,t): \nu \to \mu \times \lambda$ is injective, and $(k_{i,j})_{i<\lambda, j<\mu}$ is such that for every $i$ the set $\{j<\mu: k_{i,j}\neq 0\}$ is finite, then $(k_{s(i), t(i)} \cdot  v_{s(i)})_{i<\nu} \in \mathfrak{S}_{\nu}$ and 
		\begin{equation*}
			\Sigma_{\nu} (k_{s(i), t(i)} \cdot v_{s(i)})_{i<\nu} = \Sigma_{\lambda}\left( c_i \cdot v_i\right)_{i<\lambda} \; \text{where}\; c_i:=\sum_{\substack{j\in t(s^{-1}(i))}} k_{i,j}.
		\end{equation*}
	\end{enumerate}
	
	The arrows $f: (V, (\mathfrak{S}_\lambda)_{\lambda \in \Card}, (\Sigma_\lambda)_{\lambda \in \Card}) \to (V', (\mathfrak{S}'_\lambda)_{\lambda \in \Card}, (\Sigma'_\lambda)_{\lambda \in \Card})$ of $\Sigma\Vect'$ are linear functions $f: V \to V'$ such that for all $\lambda \in \Card$, 
    \[f^\lambda \mathfrak{S}_\lambda \subseteq \mathfrak{S}'_\lambda \quad \text{and} \quad
	\Sigma_\lambda(f(v_i))_{i<\lambda} = f\Sigma_{\lambda}(v_i)_{i<\lambda}.\] 
    Here $f^\lambda: V^\lambda \to V^{\lambda'}$ is the function acting component-wise as $f$.
\end{definition}
	
\begin{remark}
	The axioms above are equivalent to the axiomatization of the \emph{ultrafinite summability spaces} as defined in \cite[Def.~1.1]{bagayoko2024automorphisms}.
\end{remark}
	
\begin{remark}
	For finite $\lambda$, condition \ref{SV:B1} implies that $\mathfrak{S}_\lambda=V^\lambda$.
\end{remark}
	
\begin{remark}\label{rmk:setsums_in_SigmaV'}
	We note the following consequences of \ref{SV:B0}, \ref{SV:B1}, and \ref{SV:B2}:
	\begin{itemize}
		\item[(a)] for finite $\lambda$, $\Sigma_\lambda (v_{i})_{i<\lambda}$ is the sum of the $v_i$s in the sense of $V$, $\Sigma_\lambda (v_{i})_{i<\lambda}=\sum_{i<\lambda} v_i$;
		\item[(b)] for every $\lambda\in \Card$, for every bijection $r: \lambda \to \lambda$, $\Sigma_\lambda(v_i)_{i< \lambda}=\Sigma_\lambda(v_{r(i)})_{i< \lambda}$;
		\item[(c)] for any $\mu < \lambda$, and every $(v_i)_{i< \mu} \in \mathfrak{S}_\mu$, setting $v_i'=v_i$ for $i<\mu$ and $v_i'=0$ otherwise, we have $(v_i')_{i< \lambda} \in \mathfrak{S}_\lambda$ and $\Sigma_{\mu}(v_i)_{i<\mu} = \Sigma_{\lambda} (v_i')_{i<\lambda}$;
		\item[(d)] if $(v_i)_{i<\lambda} \in\mathfrak{S}_\lambda$ and $(k_i)_{i<\lambda} \in \bfk^\lambda$, then $(v_i \cdot k_i)_{i <\lambda} \in \mathfrak{S}_\lambda$.
	\end{itemize}
	In particular for an object $X= (V, (\mathfrak{S}_\lambda)_{\lambda \in \Card}, (\Sigma_{\lambda})_{\lambda \in \Card})$ of $\Sigma\Vect'$ we can give the following unambiguous definitions for any small set $I$:
	\begin{itemize}
		\item a family $(v_i)_{i \in I}$ is \emph{summable} if for some (or equivalently every, by (b)) bijection $r: |I| \to I$, $(v_{r(i)})_{i<|I|} \in\mathfrak{S}_{|I|}$;
		\item if $(v_i)_{i \in I}$ is a summable family, its \emph{sum} is defined as
		\[\sum_{i \in I} v_i:=\Sigma_{|I|}(v_{r(j)})_{j<|I|}\]
		for any bijection $r: |I| \to I$ (again this does not depend on the choice of $r$, by (b)).
	\end{itemize}
\end{remark}

\begin{theorem}\label{thm:SigmaVect_as_mod}
	$\Sigma\Vect$ is equivalent to $\Sigma\Vect'$.
	\begin{proof}
		For $X \in \Sigma\Vect$, set $FX:= (V, (\mathfrak{S}_{\lambda})_{\lambda \in \Card}, (\Sigma_\lambda)_{\lambda \in \Card})$ where $V=X\bfk$, $\mathfrak{S}_\lambda(V):= (X\delta_i)_{i<\lambda} \cdot (X (\bfk^{\oplus \lambda}))$ and for $v \in X (\bfk^{\oplus \lambda})$
		\[\Sigma_\lambda ((X\delta_i) \cdot f)_{i < \lambda} := (X \sigma) \cdot v.\]
		Clearly $\mathfrak{S}_\lambda (V)$ is a subspace of $V^\lambda$ and $\Sigma_\lambda$ is linear.
			
		Assume $\mu, \lambda, \nu \in \Card$ and $((k_{i,j})_{i<\lambda})_{j<\mu}$ is such that for every $i$, $\{j: k_{i,j} \neq 0\}$ is finite. Let $(v_i)_{i< \lambda} \in \mathfrak{S}_\lambda(V)$ so $v_i=(X\delta_i)\cdot v$ for some $v \in X(\bfk^{\oplus \lambda})$. Consider the map $\varphi: \bfk^{\oplus \lambda} \to \bfk^{\oplus (\lambda \times \mu)}$ given by $\varphi(e_i)=\sum_{j} e_{i,j} \cdot k_{i,j}$, so that $\delta_{(i,j)} \cdot \varphi = k_{i,j} \cdot \delta_i$. Now 
		\[X \delta_{(i,j)} \cdot X \varphi \cdot v = X(k_{i,j} \cdot \delta_i) \cdot v = v_i \cdot k_{i,j}\]
		so $(v_i \cdot k_{i,j})_{i<\lambda, j<\mu}$ is summable, hence for every injection $(s,t): \nu \to \lambda \times \mu$,
		$(v_{s(i)}\cdot k_{s(i), t(j)})_{i<\nu} \in \mathfrak{S}_\nu$ and the identity in \ref{SV:B2} is easily checked as well.
			
		By Lemma~\ref{lem:sum-pres_characterization}, the construction extends to a unique fully faithful functor $F: \Sigma\Vect \to \Sigma\Vect'$, such that the map $V \to V'$ underlying an arrow $F\alpha: FX \to FX'$ is $\alpha\bfk$.
			
			
		Observe that $F\Yon(\bfk^{\oplus I}) \cong (\bfk^I, (\mathfrak{S}_\lambda)_{\lambda \in \Card}, (\Sigma_\lambda)_{\lambda \in \Card})$ where 
		\[\mathfrak{S}_\lambda=\big\{((k_{i,j})_{i\in I})_{j<\lambda}: \forall i \in I, \; |\{j<\lambda: k_{i,j}\}|<\aleph_0 \big\}\quad \text{and}\]
		\[\big(\Sigma_\lambda((k_{i,j})_{i \in I})_{j<\lambda}\big)_i=\sum_{j<\lambda} k_{i,j}.\]
		Thus $\Sigma\Vect'_\bfk$ extends the category of strong vector spaces of the form $\bfk^I$ and strongly linear maps among them. In fact it also satisfies conditions \ref{SV:H1} to \ref{SV:H4} in Remark~\ref{rmk:Eur_0}. Among the four, the only one that requires a comment is \ref{SV:H2}: i.e.\ that any summable $(v_i)_{i\in I}$ induces a strongly linear function 
        \[F\Yon(\bfk^{\oplus I}) \to (V,(\mathfrak{S}_\lambda)_{\lambda \in \Card}, (\Sigma_\lambda)_{\lambda \in \Card}),\quad \text{via} \quad v: (k_i)_{i\in I} \mapsto \sum_{i \in I} k_i \cdot v_i.\]
		A summable $J$-indexed family in $F\Yon (\bfk^{\oplus I})$ (that is, $\bfk^I$ with the usual summability notion) is precisely a family $((k_{i,j})_{i \in I})_{j \in J}$ such that $|\{j \in J: k_{i,j}\neq 0\}|<\aleph_0$ for every $i\in I$. The image of such a family by $v$ is $(\sum_{i\in I}k_{i,j}\cdot v_{i})_{j \in J}$, but this is summable and its sum agrees with $\sum_{i \in I}\sum_{j\in J}k_{i,j}\cdot v_{i}$ because \ref{SV:B2} entails that $(k_{i,j}\cdot v_{i})_{(i,j) \in I \times J}$ is summable.

        It follows that $F \Yon : \Vect^\op \to \Sigma\Vect'$ is a reasonable category of strong vector spaces and therefore $\Sigma\Vect'(F\Yon, -): \Sigma\Vect' \to [\Vect, \Vect]$ is fully faithful and factors through the inclusion $\Sigma\Vect \subseteq [\Vect, \Vect]$.

        Thus to prove the statement it suffices to show that the fully faithful functor $\Sigma\Vect'(F\Yon, -)$ is essentially surjective, i.e.\ that for every $X \in \Sigma\Vect$, there is $Y \in \Sigma\Vect'$ such that $X \cong \Sigma\Vect'(F\Yon, Y)$. But by the fully faithfulness of $F$ we have natural isomorphisms
        \[\Sigma\Vect'(F\Yon (V), F X) \cong \Sigma\Vect(\Yon (V), X) \cong X(V),\]
        whence $\Sigma\Vect'(F\Yon, -)$ is essentially surjective.
	\end{proof}
\end{theorem}

\section{Orthogonality, Smallness, and Reflectivity}\label{sec:OrtSmallRefl}
	
The main goal of this section is to show that $\Sigma\Vect$ only consists of functors that are \emph{small} colimits of representables (Subsection~\ref{ssec:small}) so consequently $\Sigma\Vect$ is a full \emph{reflective} subcategory of $\Ind\-(\Vect^\op)$. Most of the results apply to the more general set-up of some specific (pre-)torsion theories of $[\Vect, \Vect]$ (described in Subsection~\ref{ssec:setup}) and their restriction to $\Ind\-(\Vect^\op)$: since the generality is cost-free and will be used in Section~\ref{sec:Comp}, we will work in such a context.
To this end, we start with a very brief introduction to the basic machinery of factorization systems and orthogonality which will be used throughout the chapter: this is the subsection below and the reader who has some familiarity with these notions may just skim through the definitions to get used to the notations. This is motivated by the fact that $\Sigma\Vect$ and also some other categories considered in Section~\ref{sec:Comp} in fact consist of the torsion free part for a certain torsion theory on $[\Vect,\Vect]$.

Given an arrow $f: X \to Y$ in a pointed category $\cC$ with kernels and cokernels, we will use the standard notations $\Kr(f)$, $\Ckr(f)$, $\Img(f)$, $\Coimg(f)$ for the objects respectively of the kernel, cokernel, image, and cokernel, and we will use their lower case version $\kr(f)$, $\ckr(f)$, $\img(f)$, $\coimg(f)$ for the corresponding connecting arrows. 
	
\subsection{Factorization systems and torsion theories}\label{ssec:OrtRecap}
Throughout this subsection $\cC$ will always denote a category.

As was shown in Lemma~\ref{lem:simple_char}, under condition \ref{SV:A1}, condition \ref{SV:A2} says that all objects $c$ of the  category $\cC$ are such that $\cC(-,c)$ maps a certain class of maps $\cA$ to monomorphisms of sets (or of vector spaces, in the $\bfk$-enriched setting).
Similarly, once we assume \ref{SV:A1*} and that $\cC$ contains all objects of the form $\Fin V$ for $V \in \Vect$, \ref{SV:A3*} can also be phrased in such form (see Lemma~\ref{lem:SigmaVect_char}~(4)).

This is indeed a naturally occurring construction: it amounts to require that objects $c\in \cC$ ``see'' the arrows in the class $\cA$ as epimorphisms. A natural example is: in the category of topological spaces the separated ones are those which ``see'' dense maps as epimorphisms.
	
If $\cC$ has cokernel pairs (i.e.\ push-outs of the spans fo the form $b \overset{f}{\longleftarrow} a \overset{f}{\longrightarrow} b$) and coequalizers, then the construction is actually an orthogonality construction: given an arrow $f$, if $(g, g')$ is the cokernel pair of $f$, then to ask that $\cC(f, c)$ is one-to-one is equivalent to require that $\cC(g, c)=\cC(g',c)$.
If $h$ is the coequalizer of $(g,g')$, the equality $\cC(g, c)=\cC(g',c)$, in turn, is equivalent to $\cC(h,c)$ being an isomorphism, i.e.\ to the fact that $c$ is \emph{right-orthogonal} to $h$ (see Definition~\ref{def:orthogonal} below).
	
In the case of $\bfk$-additive categories, as observed in the equivalences $(3)\Leftrightarrow (2)$ of Lemma~\ref{lem:simple_char} and $(4)\Leftrightarrow (5)$ in Lemma~\ref{lem:SigmaVect_char} this becomes even easier: if $\cC$ is $\bfk$-additive and has cokernels, $c$ is an object and $f$ an arrow, then $\cC(f, c)$ is a one-to-one if and only if $\cC(\Ckr f, c)=0$.
	
The main definitions recalled in this subsection are that of orthogonality in category theory (Definition~\ref{def:orthogonal}) and that of a factorization system (Definition~\ref{def:factorization}) first introduced in a primitive version in \cite{maclane1948groups} and \cite{isbell1957remarks} under the name of bicategorical structure or bicategory and later developed in \cite{ringel1970diagonalisierungspaare} and \cite{freyd1972categories} (cf \cite[1.4]{freyd1972categories} for a more detailed yet brief account of the earliest contributions, \cite{riehl2008factorization} for a more modern exposition).
We will be interested in the particular case of this given by a torsion theory (Definition~\ref{def:torsion_theory}), whose definition in the case of Abelian categories is from \cite{dickson1966torsion}. None of the definitions or results presented in this subsection are original.
	
\begin{definition}\label{def:orthogonal}
	Recall the notion of orthogonality from from \cite[2.1]{freyd1972categories} and \cite[Sec.~3]{ringel1970diagonalisierungspaare}: in a category $\cC$, one says that a map $f \in \cC(a,b)$ is \emph{left orthogonal} to $g\in \cC(c,d)$ (or that $g$ is \emph{right orthogonal} to $f$) and writes $f \ort g$, if and only if the square
	\[\begin{tikzcd}
		\cC(b,c) \ar[rr, "{\cC(b,g)}"] \ar[d, "{\cC(f, c)}"'] &&\cC(b,d)\ar[d,"{\cC(f,d)}"]\\
		\cC(a,c) \ar[rr, "{\cC(a,g)}"] &&\cC(a,d)
	\end{tikzcd}\]
	is a \emph{pullback}, or equivalently if for every $(l,k) : f \to g$ in $(\id_\cC \downarrow \id_\cC)\cong \Func(\two, \cC)$ there is a unique $h \in \cC(b, d)$ with $l=h\cdot f$ and $k = g \cdot h$.
    
	This orthogonality relation on arrows is usually extended to a relation defined on arrows and objects as follows.
		
	An object $d$ is \emph{right orthogonal} to an arrow $f : c \to c'$, written $f \ort d$, if and only if the map
	\[\cC(f, d)= - \cdot f : \cC(c', d) \to \cC(c,d)\]
	is an isomorphism. It is \emph{left orthogonal}, written $d \ort f$, if it is right orthogonal in $\cC^\op$.
		
	Finally an object $d$ is \emph{right orthogonal} to an object $c$, written $c \ort d$ if $|\cC(c, d)|=1$ (i.e.\ if there is only one arrow from $c$ to $d$).
\end{definition}

\begin{remark}\label{rmk:orthogonality_with_extreme_objects}
	In the definition above if $\cC$ has a terminal object $\top$, $c$ is an object, and $f$ is an arrow, then $f \ort c$ if and only if $f \ort\, (c \to \top)$. Similarly if $\cC$ has an initial object $\bot$, then $c \ort f$ if and only if $(\bot \to c) \ort f$. If $\cC$ has both a terminal and an initial object and $c,d$ are objects, then $c \ort d$ if and only if $(\bot \to c) \ort \,(d \to \top)$.
\end{remark}

\begin{definition}
	If $\cA$ is a class of arrows or objects in $\cC$, then $\cA^{\bot}=\{m: \forall a \in \cA,\; a\ort m\}$ will denote the class of \emph{arrows} that are right orthogonal to everything in $\cA$ and similarly $^{\bot}\cA=\{e: \forall a \in \cA,\; e\ort a\}$ will denote the set of arrows that are left orthogonal to everything in $\cA$. If $\cA$ consists of arrows, then a class of the form $\cA^{\bot}$ will be called a \emph{right orthogonal class of arrows} whereas a class of the form ${}^{\bot}\cA$ will be called a \emph{left orthogonal class of arrows}. 
\end{definition}

\begin{remark}
	If $\cC$ has an initial object, then $\cA^\bot$ is always a right orthogonal class of arrows, and similarly if $\cC$ has a terminal object, then ${}^\bot\cA$ is always a left orthogonal class of arrows. From now on we will assume that the category $\cC$ in consideration always has an initial and a terminal object.
\end{remark}


\begin{remark}
	The operators $-^{\bot}$ and $^{\bot}-$ on classes of arrows define an antitone Galois connection, in the sense that for classes of arrows $\cA$ and $\cA'$, one has $\cA' \subseteq \cA^\bot$ if and only if $\cA \subseteq {}^\bot\cA'$. In particular right orthogonal classes of arrows are equivalently characterized as the classes $\cM$ such that $(^{\bot}\cM)^{\bot}=\cM$ and dually the left orthogonal class of arrows are those $\cE$ such that $^{\bot}(\cE^\bot)= \cE$.
\end{remark}

\begin{remark}\label{rmk:Facts_On_Right_Orthogonal_Classes}
	If $\cM$ is a right orthogonal class of arrows in a category $\cC$ then it always has the following properties (cf \cite[Proposition~2.1.1]{freyd1972categories} plus the observation that $\cM$ is closed under pointwise equalizers in $\Func(2, \cC)$):
	\begin{enumerate}
		\item it contains all isomorphisms and is closed under composition;
		\item the pullback of an $m \in \cM$ along any arrow $f$ in $\cC$ is in $\cM$;
		\item the class $\cM$ in $\Func(2, \cC)$ is closed under pointwise limits;
		\item for every $c \in \cC$ the class $\cM_c$ of objects of $(\cC\downarrow c)$ consisting of arrows $(m:x \to c) \in \cM$ is closed under limits in $(\cC\downarrow c)$;
		\item if $f \cdot m \in \cM$ and $f \in \cM \cup \Monomorphisms$ then $m \in \cM$.   
	\end{enumerate}
	Another fact we are going to use is 
	\begin{enumerate}[resume]
		\item if $\cM \subseteq \Monomorphisms$ and $\cC$ has pullbacks, then $e \ort \cM$ if and only if whenever $e= f \cdot m$ with $m \in \cM$, $m$ has to be an isomorphism \cite[1.4, Prop.~1.1, i]{kelly1980unified}.
	\end{enumerate}
	Formally dual statements hold for left orthogonal classes of arrows. 
\end{remark}

\begin{definition}\label{def:factorization}
	Following \cite[2.1]{freyd1972categories} we will call a pair $(\cE, \cM)$ of classes of arrows an \emph{orthogonal prefactorization system} or just \emph{prefactorization system} if $\cE$ is a left orthogonal class and $\cM=\cE^\bot$. A prefactorization system $(\cE, \cM)$ is called a \emph{orthogonal factorization system} or just \emph{factorization system} when every arrow $f$ factors as $f= m \cdot e$ for some $m \in \cM$ and $e \in \cE$. Associated with a prefactorization $(\cE, \cM)$ are the classes of objects
	\[\begin{aligned}
	\Free(\cM)&\coloneqq \{c : \cE \ort c\} = \{c : (c \to \top) \in \cM\},\\
	\Tors(\cE)&\coloneqq \{c : c \ort \cM\} = \{c : (\bot \to c) \in \cE\}.
	\end{aligned}\]
	A factorization system is \emph{proper} if $\cE \subseteq \Epimorphisms$ and dually $\cM \subseteq \Monomorphisms$.
	Prefactorization systems are naturally ordered, the usual order is the one of inclusions of their right part, i.e.\ $(\cE, \cM) \le (\cE', \cM')$ if and only if $\cM \subseteq \cM'$.
\end{definition}

\begin{example}
	In an Abelian category $\cC$, there is a unique proper factorization system $(\cE, \cM)$ given by
	\[\cE=\Epimorphisms \quad \text{and} \quad \cM=\Monomorphisms.\]
	For such a factorization system the classes $\Free(\cM)$ and $\Tors(\cE)$ are given just by the zero objects and the factorization of a map $f: c_0 \to c_1$ has as $\cE$-component $\coimg f := \ckr \kr (f): c_0 \to \Coimg f\cong \Img f$ and as $\cM$-component $\img f := \kr \ckr f : \Img f \to c_1$.
\end{example}

\begin{remark}\label{rmk:uniqueness_of_orth_fact}
	If $f = m \cdot e = m' \cdot e'$, with $e' \ort m$ and $e \ort m'$, then the factorizations are isomorphic in the sense that there is an isomorphism $h$ such that $e' = h \cdot e$ and $m'= m \cdot h^{-1}$. More generally, if $(h,k): f \to g$ is an arrow in $\Func(2, \cC)$ (so $k\cdot f =g \cdot h$) and 
	$f=m \cdot e$, $g=m' \cdot e'$ with $e \ort m'$, then there is a unique $l$ such that $e' \cdot h = l \cdot e$ and $m' \cdot l = k \cdot m$.
	In particular if $(\cE, \cM)$ is an orthogonal factorization system, for every arrow $f \in \cC(a,b)$ the $\cM$-component $m=\cM(f) \in \cM$ of a factorization $f = m \cdot e$  with $m \in \cM$ and $e \ort \cM$ is determined up to isomorphism. Moreover (all choices of) the assignments $f\mapsto \cM(f)$ for $f \in \Func(2, \cC)$ lift to a functor $\cM: \Func(2, \cC) \to  \cM$, which is uniquely determined by the factorization system up to isomorphism.
	Such functor $\cM: \Func(2, \cC) \to  \cM$ is a reflector, the unit of the adjunction being given for an object $(f:a \to b)$ of $\Func(2, \cC)$ by the pair $(e, 1_b): f \to \cM(f)$, where $e$ is s.t.\ $e \ort \cM$ and $f= \cM(f) \cdot e$ (see \cite[Thm.~3.2]{grandis2006natural}).
\end{remark}
	
\begin{remark}\label{rmk:Fact&Refl}
	It follows in particular that if $(\cE, \cM)$ is a factorization system, then $\Psi(\cE, \cM):=\Free(\cM)$ is a \emph{$\cE$-reflective subcategory} of $\cC$, that is, the reflection morphisms are all in $\cE$, the reflection being given by $\rho_c=\cE(c \to \top)$. More generally this happens whenever $(\cE, \cM)$ is a pre-factorization system satisfying
	\begin{enumerate}[label=(Sf), ref=(Sf), align = left, leftmargin =*, itemindent = !]
		\item\label{Sf} for every $c\in \Ob(\cC)$, the arrow $c \to \bot$ has a $(\cE, \cM)$-factorization.
	\end{enumerate}
	Conversely, given a replete reflective subcategory $\cD$, one can form a pre-factorization system $\Phi(\cD)=(\cE, \cM)$ satisfying \ref{Sf} by setting $\cE = {}^\bot\cD$. Such $\cE$ is then the biggest class such that $\Free(\cE^\bot)=\cD$, thus $(\cE, \cM)$ is the minimum prefactorization system with $\Free(\cE^\bot)=\cD$.
	These correspondences define a monotone \emph{Galois co-insertion} (i.e.\ $\Phi(\cD)\le (\cE, \cM) \Leftrightarrow \cD \subseteq \Psi(\cE, \cM)$ and $\Psi \circ \Phi = \mathrm{Id}$) from the poset of replete reflective subcategories into the one of pre-factorization systems satisfying \ref{Sf}(cf \cite[p.~4 and Prop.~2.2]{cassidy1985reflective}). In particular the map $\Phi \circ \Psi$ defines an interior operator on the class of pre-factorization systems satisfying \ref{Sf} called \emph{reflective interior}. All pre-factorzation systems $(\cE, \cM)$ satisfying \ref{Sf} which equal their interior satisfy the extra property that if $e \cdot f \in \cE$ and $e \in \cE$, then $f \in \cE$ as well. This characterizes the pre-factorization systems equaling their interior within the ones satisfying \ref{Sf} (cf \cite[Thm.~2.3]{cassidy1985reflective}) and motivated the following
\end{remark}
	
\begin{definition}\label{def:ReflFact}
	A prefactorization system $(\cE, \cM)$ with the property that if $e \cdot f \in \cE$ and $e \in \cE$, then $f \in \cE$ as well, is called \emph{reflective} (cf \cite[3.2]{rosicky2008factorization} and \cite[p.291]{cassidy1985reflective}).
\end{definition}
	
From the dual of Remark~\ref{rmk:Facts_On_Right_Orthogonal_Classes}(5), it follows that if $(\cE, \cM)$ is reflective then $\cE$ has the \emph{3-for-2 property}, i.e.\ whenever $|\cE \cap \{f,g, g\cdot f\}|\ge 2$ actually $\{f,g, g\cdot f\}\subseteq \cE$.

\begin{remark}
	Remark~\ref{rmk:Fact&Refl} and the consequent Definition~\ref{def:ReflFact} dualize in the obvious way, giving rise to the concept of \emph{coreflective} pre-factorization system and \emph{coreflective closure} of a pre-factorization system.
\end{remark}
	
\begin{definition}\label{def:torsion_theory}
	A (pre-)factorization system $(\cE, \cM)$ in which both $\cE$ and $\cM$ have the 3-for-2 property is called a \emph{(pre-)torsion theory} (see \cite[4.4]{rosicky2008factorization}; in \cite[p.~21]{cassidy1985reflective} it has different definition which is essentially equivalent in Abelian categories, \cite[Thm~8.20]{cassidy1985reflective}).
	For a (pre-)torsion theory $(\cE, \cM)$, the classes $\Free(\cM)$ and $\Tors(\cE)$ are called respectively \emph{class of the torsion-free objects} and \emph{class of the torsion objects}.
\end{definition}

\begin{remark}\label{rmk:GenTorsionTheories}
	A pre-torsion theory $(\cE, \cM)$ is uniquely determined by the full subcategories $\cF=\Free(\cM)$ and $\cT=\Tors(\cE)$, through
	\[\cE={}^\bot\cF=\{e : e \ort \cF\}, \qquad \cM = \cT^\bot = \{m : \cT \ort m\}.\]
	Moreover any pair $(\cT, \cF)$ such that $\cT = \Tors({}^\bot\cF)$ and $\cF= \Free(\cT^\bot)$ defines a left orthogonal class of arrows $\cE$ and a right orthogonal one $\cM$ by the above formulae; however the reflective prefactorization system $(\cE, \cE^\bot)$ and coreflective one $({}^\bot\cM, \cM)$ need not to agree in general. When they do, or equivalently when $\cE \ort \cM$, we say the pair $(\cT, \cF)$ \emph{induces the pretorsion theory} $(\cE, \cM)$; if, moreover, the induced pretorsion theory is a torsion theory, we say it \emph{induces a torsion theory}.
\end{remark}

\begin{remark}
	Recall that a category $\cC$ is \emph{finitely well complete} in the sense of \cite[p.~292]{cassidy1985reflective} if it is finitely complete (so all has finite limits) and \emph{all} (even large) intersections of \emph{strong} monomorphisms (i.e.\ mononrphisms that are right orthogonal to all epimorphisms).
    By \cite[Cor.~3.4]{cassidy1985reflective}, if $\cC$ is finitely well complete then the pre-torsion theory induced by a pair $(\cT, \cF)$ as in Remark~\ref{rmk:GenTorsionTheories}, is a torsion theory if and only if $\cF$ is $\cE$-reflective, and when this is the case $\cT$ is $\cM$-coreflective.
\end{remark}

\begin{definition}
    Following \cite[5.1]{rosicky2008factorization}, given a finitely bicomplete category $\cC$ with $0$ object, we call \emph{standard torsion theory} on $\cC$ a pair $(\cT, \cF)$ of classes of objects such that $\cT=\Tors({}^\bot\cF)$, $\cF=\Free(\cT^\bot)$ and for all objects $c \in \cC$, there are $\tau : t \to c$ and $\varphi: c \to f$ such that $t \in \cT$, $f \in \cF$, $\kr (\varphi) \cong \tau$, and $\ckr(\tau)\cong\varphi$.
\end{definition}

\begin{remark}
    If for all objects $c \in \cC$, there are $\tau : t \to c$ and $\varphi: c \to f$ such that $t \in \cT$, $f \in \cF$, $\kr (\varphi) \cong \tau$, and $\ckr(\tau)\cong\varphi$, then provided $\cT$ and $\cF$ are closed under isomorphisms, $(\cT, \cF)$ is a standard torsion theory if and only if ${}^\bot\cF \ort \cT^\bot$. 
\end{remark}

The following Proposition is well known and follows from much more general results \cite[Thm.s~8.18 and 8.20]{cassidy1985reflective}, or \cite[Thm.~5.2]{rosicky2008factorization}: it entails the equivalence of the notions of torsion theory and standard torsion theory in an Abelian category. A short standard proof is included for the convenience of the reader.

\begin{proposition}\label{prop:TorsionTheoriesAbelian}
    If $\cC$ is an Abelian category then the following are equivalent for a pair $(\cF, \cT)$ such that $\cT=\Tors({}^\bot \cF)$ and $\cF=\Free(\cT^\bot)$:
    \begin{enumerate}
        \item $(\cT, \cF)$ induces a torsion theory;
        \item $(\cT, \cF)$ is a standard torsion theory;
        \item $\cF$ is reflective;
        \item $\cT$ is coreflective.
    \end{enumerate}
    \begin{proof}
        Note that $(1) \Rightarrow (3)$ and $(1) \Rightarrow (4)$ are trivial.
        Also observe that since by construction $\cF=\Free(\cT^\bot)$ is closed under subobjects if $\cF$ is reflective, then it is epi-reflective: indeed let $\rho : Id_\cC \to F : \cC \to \cC$ be the reflection of $\cC$ onto $\cF$. The image $\Img(\rho c)$ of each $\rho c: c \to Rc$ is in $\cF$, whence if $\rho c = m \cdot e$ is the epi-mono factorization of $\rho c$, we have for all $f \in \cF$, $\cC(\rho c, f) = \cC(e, f) \cdot \cC(m, f)$. Thus since $\cC(\rho c, f)$ is bijective, $\cC(e,f)$ must be bijective for all $f \in \cF$, but given that $\Img(\rho c) \in \cF$, this implies that also $e: c \to \Img(\rho c)$ is initial in $(c \downarrow \cF)$ and thus $m$ is an isomorphism.

        $(3) \Rightarrow (2)$. Let $\rho$ denote the reflection as above.
        First observe that $\Kr(\rho c)$ is in $\cT$. In fact, for any map from $h: \Kr(\rho c) \to f\in \cF$, the push-forward $h': c \to f'$ of $h$ along $\kr(\rho c)$ has as codomain $f'$ an extension of $f$ by $F c$, whence it is in $\cF$.
        But then $h'$ factors through $\rho c$ and thus $0 = h' \cdot \kr(\rho c) = k' \cdot h$ where $k'$ is a monomorphism (because monomorphisms are stable under pushouts in Abelian categories) whence $h=0$. Thus we proved $\Kr(\rho c) \in \cT$ and (2) holds.

        $(2) \Rightarrow (3)$ readily follows from the fact that if $0 \to t \to c \to f \to 0$ is exact with $t \in \cT$ and $f \in \cF$, then for all $f' \in \cF$, $\cC(c, f')\cong \cC(f, f')$ whence $f$ would be the reflection of $c$ onto $\cF$.

        $(2) \Rightarrow (1)$. By $(2)\Leftrightarrow (3)$ and its dual $(2) \Leftrightarrow (4)$ we have that if $(2)$ holds, then we have a reflection $\rho : Id_\cC \to F : \cC \to \cC$ and a coreflection $\sigma: S \to Id_\cC : \cC \to \cC$ of $\cC$ onto $\cF$ and $\cT$ respectively. Thus to conclude it suffices to show that ${}^\bot \cF \ort \cT^\bot$.
        Consider the morphism of exact rows
	    \[\begin{tikzcd}
	    0 \ar[r] &Sc_0 \ar[r]\ar[d,"{Sf}"] &c_0 \ar[r]\ar[d, "f"]& Rc_0 \ar[r]\ar[d, "{Rf}"] &0\\
		  0 \ar[r] &Sc_1 \ar[r] &c_1 \ar[r]& Rc_1 \ar[r] &0 
	    \end{tikzcd}\]
	    and observe that by definition of $\cT^\bot$, if $f\in\cT^\bot$ then $S(f)$ is an isomorphism. Since the diagram has exact rows, this happens if and only if the rightmost square is a pullback. Now it is easy to see that since $g \in {}^\bot\cF$ if and only if $Rg$ is an isomorphism, we must have ${}^\bot \cF \ort h$ for every $h$ between objects in $\cF$, in particular ${}^\bot \cF \ort Rf$, so since $({}^\bot\cF)^{\bot}$ is closed under pullbacks we have $f \in ({}^\bot\cF)^{\bot}$. Thus (1) is proven.
    \end{proof}
\end{proposition}

\subsection{Some special torsion theories in \texorpdfstring{$[\Vect, \Vect]$}{[Vect,Vect]ₖ}}\label{ssec:setup}
In the remainder of this paper we are going to be interested in torsion theories on $[\Vect, \Vect]$ generated by quotients of representables.

\begin{remark}
	As remarked before, the class of those $X$ such that $\Nat(\alpha,X)$ is one-to-one for every $\alpha$ in a given class of arrows $\mathcal{A}$ is actually the class
	\[\Free(\cQ^\bot)=\Free(\{0 \to Q: Q \in\cQ\}^\bot)=\{X : \forall Q \in \cQ,\,\Nat(Q,X)=0\}\]
	where $\cQ = \{\Ckr\, \alpha : \alpha \in \mathcal{A}\}$. For simplicity we will call an $X \in\Free(\cQ^\bot)$ a \emph{$\cQ$-free object}.
	By Remark~\ref{rmk:GenTorsionTheories} and Proposition~\ref{prop:TorsionTheoriesAbelian}, the pair $\big({}^\bot \Free(\cQ^\bot), \Tors({}^\bot\Free(\cQ^\bot))^\bot\big)$ will then be a torsion theory as soon as $\Free(\cQ^\bot)$ is reflective.
\end{remark}
	
In the cases we are interested in, $\cA$ will always be a class of monomorphisms into representables.
In this regard it is useful to know what subobjects of a representable in $[\Vect, \Vect]$ look like.
	
\begin{lemma}\label{lem:subobjects_of_representables}
	Every subobject $Z \subseteq \Vect(V, -)$ which is small (i.e.\ a small colimit of representables) is up to isomorphism of the form
	\[Z\cong \bigcup_{H \in \cF} \Vect(V/H, -)\]
	for some filter $\cF$ of subspaces of $V$.
	\begin{proof}
		If $Z$ is a small subobject of $\Vect(V,-)$, it is the image of some 
	    \[A=[\Vect(\alpha_i, -)]_{i \in I}: \bigoplus_{i \in I} \Vect(W_i, -) \to \Vect(V,-)\]
        for some small family $\{W_i: i \in I\}$ of vector spaces and linear maps $\alpha_i \in \Vect(V, W_i)$.
		Define
		\[\cF= \{H < V : \exists F \in [I]^{<\omega},\; H\supseteq \bigcap_{i \in F} \Kr \alpha_i \}.\]
		It is then easily verified that the image of $A$ is $\bigcup_{H \in \cF} \Vect(V/H, -)$.
	\end{proof}
\end{lemma}

\begin{remark}
	If the large cardinal $|\mathcal{U}|=\sup\{|x|: x \in \mathcal{U}\}$ of the chosen Grothendieck universe $\mathcal{U}$ is a Vopenka cardinal, then every $Z \subseteq \Vect(V,-)$ has to be small. See \cite[Cor.~6.31]{adamek1994locally}.
\end{remark}
	
\begin{definition}\label{defn:cQ()}
	For $V \in \Vect$ and $\cF$ a filter of subspaces of $V$, denote by $a_{V, \cF}$ the inclusion $\bigcup_{H \in \cF} \Yon (V/H) \subseteq \Yon V$.
	For $\bE$ a (possibly large) class of pairs $(V, \cF)$ with $V \in \Vect$ and $\cF$ a filter of subspaces of $V$ define the classes 
	\[\mathcal{A}(\bE)=\{a_{V, \cF}: (V, \cF) \in \bE\}, \qquad \cQ(\bE)=\{\Ckr\ a : a \in \mathcal{A}(\bE)\}.\]
	Call $\bE$ \emph{degenerate}, if there is $(V, \cF) \in \bE$ such that $\bigcap \cF\neq 0$ and \emph{non-degenerate} otherwise. 
\end{definition}
	
\begin{lemma}
	The following are equivalent for $\bE$ as above:
	\begin{enumerate}
		\item $\bE$ is non-degenerate;
		\item all representables are $\cQ(\bE)$-free.
		\item there is a non-zero $\cQ(\bE)$-free representable.
	\end{enumerate}
	\begin{proof}
		To prove $(1) \Rightarrow (2)$, let $W \in \Vect$, $(V, \cF) \in \bE$ and $Q:=\Ckr a_{V,\cF}$. Given any map $f: Q \to \Yon W$, its composite with $f \cdot \ckr(a_{V,\cF}) : \Yon V \to \Yon W$, by Yoneda must have the form $\Yon h$ for $h\in \Vect(W, V)$ and thus must have a representable kernel. Thus $\Kr (f\cdot a_{V, \cF}) := \Yon (V/ H) \le \Yon V$ with $H \le \bigcap \cF$, but by $(1)$, $\bigcap \cF=0$ so $H=0$ and $\Kr( f \cdot a_{V, \cF} ) = \Yon V$. Since this holds for every $f$ and every $(V, \cF) \in \bE$, $\Yon W$ is $\cQ$-free.

		$(2) \Rightarrow (3)$ is trivial.
			
		$(3) \Rightarrow (1)$. Assume that $(1)$ does not hold and that there is $(V, \cF)\in \bE$ with $H=\bigcap \cF \neq 0$. Then $Q= \Ckr a_{V, \cF}$ has a non zero representable quotient $h : Q \to \Yon H$. Since $\Yon H$ admits an epimorphism to every $\Yon H'$ with $\dim H' \le \dim H$, and a monomorphism to every $\Yon H''$ with $\dim H'' \ge \dim H$ it follows that no representable is $Q$-free.
	\end{proof}
\end{lemma}
	
\begin{example}\label{ex:mainex}
	Regarding $\Sigma\Vect$ we are mainly interested in the following two classes
	\begin{itemize}
		\item $\bE_2:=\big\{(V,\{\text{finite codimensional subspaces of V}\}): V \in \Vect\big\}$, to which corresponds
		\[\cQ_2:=\cQ(\bE_2)=\{\Yon V/ \Fin V^*: V \in \Vect\}\]
        so $\cQ_2$-free objects are precisely \emph{simple} objects (cf Lemma~\ref{lem:simple_char});
		\item $\bE_3:=\big\{(\bfk^{\oplus \kappa}, \{\text{filter generated by the $\Kr(\delta_i)$, $i \in \kappa$}\}): \kappa \in \Card\big\}$, to which corresponds
        \[\cQ_3:=\cQ(\bE_3)=\big\{\Ckr \big(\Vect(\bfk,-)^{\oplus \kappa} \to \Vect(\bfk,-)^\kappa\big): \kappa \in \Card\big\}\]
        so an $X \in [\Vect, \Vect]$ is $\cQ_3$-free if and only if $X\in \Sigma\Vect$ (see Lemma~\ref{lem:SigmaVect_char}).
	\end{itemize}
\end{example}
	
In Subsection~\ref{ssec:Comp_Top} we will also be interested in the following special class	
	
\begin{itemize}
	\item $\bE_4:=\big\{(V, \cF): V \in \Vect, \; \bigcap\cF=0\big\}$ so $\cQ_4:=\cQ(\bE_4)$ consists of all $\Yon V/ S$ where $V \in \Vect$ and \[S = \bigcup_{H \in \cF} \Yon (V/H)\] with $\cF$ any proper filter of subspaces of $V$ such that $\bigcap \cF =0$.
\end{itemize}
	
The class of $X \in [\Vect,\Vect]$ which are small colimits of representables and are $\cQ_4$-free, in fact, coincides with the class of small directed unions of representables.
	
\begin{proposition}\label{Q4-free_char}
	Let $X \in [\Vect, \Vect]_{\bfk}$ be a small filtered colimit of representables. Then the following are equivalent
	\begin{enumerate}
		\item for every $V \in \Vect$, every arrow $\Yon V \to X$ has a representable kernel;
		\item $X$ is $\cQ_4$-free;
		\item $X$ is a directed union of representables, i.e.
        \[X \cong \varinjlim_{d \in D} \Yon Gd\]
		for some directed diagram $G: D \to \Vect^\op$ with $G(d<d')$ an epimorphism of $\Vect$. 
	\end{enumerate}
	\begin{proof}
		$(1)\Rightarrow (2)$ Let $\cF$ be a filter of subspaces on $V$ with $\bigcap \cF=0$. Given a map $f:\Ckr a_{V, \cF} \to X$, consider the kernel $\Kr(f\cdot \pi)$ of its composition with $\pi:=\ckr(a_{V, \cF}) :\Yon V \to \Ckr a_{V, \cF}$: by construction it must contain every $\Yon (V/H) \subseteq \Yon V$ for $H \in \cF$, but by $(1)$ such a kernel should be representable, so it must be $\Yon(V/\bigcap \cF)=\Yon V$ and therefore $f\cdot \pi =0$. Since $\pi$ is epi, it follows $f=0$.
			
		$(2)\Rightarrow (1)$ Let $f:\Yon V \to X$. Since $X$ is a small colimit of representables, so is $\Kr(f)$. Thus by Lemma~\ref{lem:subobjects_of_representables}, $\Kr (f) = \Img a_{\cF, V}$ for some $\cF$. Note that we have the inclusions 
		\[\Kr(f)\subseteq \Yon (V /\bigcap \cF) \subseteq \Yon V \]
		This means that there is a monomorphism $\Yon (V /\bigcap \cF) / \Kr(f) \to X$, but $(2)$ then implies that this is $0$, so $\Kr(f)=\Yon (V/\bigcap \cF)$ is representable.
			
		$(1)\Rightarrow (3)$ Recall that each filtered colimit is isomorphic to some directed colimit (cf \cite[Thm.~1.5]{adamek1994locally}), so we may assume that $X$ is a directed colimit of representables. Let $H:D \to \Vect^\op$ be a directed diagram and let $X \cong \varinjlim_{d \in D} \Yon (Hd)$. Consider the colimiting cone $\eta: \Yon H \to X$. By $(1)$ each component $\eta d$ has a representable image $\Img(\eta d)$. Since the image factorization is functorial, the map $H': d \mapsto \Img (\eta d)$ extends to a functor $H': D \to [\Vect, \Vect]$ and the natural transformation $\eta$ factors as $\eta = \iota \cdot \pi$ where $\iota d: \Img (\eta d) \to X$ and $\pi d: \Yon Hd \to \Img (\eta d)$ are the natural maps. Note that $\varinjlim H' = X$ and that by $(1)$ each $\Img (\eta d)\cong \Ckr(\ker(\eta d))$ is representable so by Yoneda $H' \cong \Yon G$ for some $G: D \to \Vect^\op$. Finally we observe that all the arrows $G(d\le d')$ are monomorphisms of $\Vect$ because $\Img (\eta d)\le \Img(\eta d')$ as subobjects of $X$.
		
		$(3)\Rightarrow (1)$ Since representables are finitely presented in $[\Vect, \Vect]$ every map $f:\Yon V \to X$ factors as $f = (\eta d) \cdot f'$ for some $\eta d: \Yon Gd \to X$ where $\eta$ is the colimit cone $\eta: \Yon G \to X$. Since $(3)$ implies\footnote{recall directed limits are exact in $\Vect$}
        that such $\eta d$ is a monomorphism, we have $\Kr f= \Kr f'$ and $f'$ has a representable kernel.
	\end{proof}
\end{proposition}

\begin{remark}
    We will see that every colimit of representables in $[\Vect, \Vect]$ is in fact a filtered colimit of representables (Remark~\ref{rmk:IndVect_and_(co)limits}). Thus in Proposition~\ref{Q4-free_char} we can drop the hypothesis that $X$ is a small filtered colimit, and just assume that it is a small colimit.
\end{remark}
	
\begin{remark}\label{rmk:Q_2,3,4_rel}
	It is easy to see that every object in $\cQ_2$ is a quotient of an object in $\cQ_3$ and that $\cQ_4 \supseteq \cQ_3 \cup \cQ_2$. It follows that every $\cQ_4$-free object is $\cQ_3$-free and every $\cQ_3$-free object is $\cQ_2$-free:
	\[\cQ_{4}\text{-free}\Rightarrow \cQ_{3}\text{-free}\Rightarrow\cQ_{2}\text{-free}.\]
	Moreover each of the implications is strict as shown by the examples below.
\end{remark}

\begin{example}
	Let $\dim W = \aleph_0$ and consider $Z \in [\Vect \Vect]$ defined by
    \[Z V = \Vect(W, V ) / (W^* \otimes_\bfk V),\]
    where we identify $W^* \otimes_\bfk V$ with the set of $f \in \Vect(W, V)$ with finite rank. We see that $Z \bfk = 0$ so $(\cN_{2} Z) V = ZV$ for every $V$ (cf Remark~\ref{rmk:cN_2}) and that $ZV \neq 0$ for $\dim V \ge \aleph_0$, so $Z$ is not $\cQ_2$-free.
\end{example}

\begin{remark}\label{rmk:concreteNat(Q,X)}
	If $Q= \Yon V/(\bigcup_{H \in \cF} \Yon (V/H))=\Ckr\, a_{\cF, V}$ for some filter $\cF$ of subspaces of $V$, then
	\[\begin{aligned}
		&\Nat(Q, X) \cong \Kr(\Nat(\Yon V,X) \to \Nat(\bigcup_{H \in \cF} \Yon (V/H), X)) \cong\\
		&\cong \bigcap_{H \in \cF} \Kr(X \pi_H)\subseteq XV,
	\end{aligned}\]
	where $\pi_H: V \to V/H$ denotes the canonical projection and the arrow inside the kernel is the restriction map. More precisely the map 
	\[(X-)\cdot - :  \bigcap_{H \in \cF} \Kr(X \pi_H) \to \Nat(Q, X),\quad x \mapsto (X-) \cdot x\]
	is a bijection.
	In particular, for the case of the canonical inclusion $\delta_I : (\Yon \bfk)^{\oplus I} \to (\Yon \bfk)^I$ relative to $\cQ_3$,
	\[\Nat\left( \frac{\Yon (\bfk^{\oplus I})}{\Fin (\bfk^{\oplus I})}, X\right)\cong \bigcap_{i \in I}\Kr(X \delta_i^I) \subseteq X(\bfk^{\oplus I}).\]
\end{remark}
	
\begin{example}\label{ex:doubledual}
	If $X$ is the double dualization functor $XV = V^{**}$ one has $\cN_{2} X = 0$ but for $\dim V = \aleph_0$ and $B$ a basis of $V$
	\[\bigcap_{b \in B} \Kr(X \delta_b^{B}) \neq 0\]
	To see this observe that $\bigcap_{b \in B} \Kr(X \delta_b^{B})$ is the annihilator of $\Span\{\delta_b^B: b \in B\}\le V^*$ and that $\Span\{\delta_b^B : b \in B\}\neq V^*$ so its annihilator is not $0$.
\end{example}

Although $X=(-)^{**}$ is interesting in its own right, it will not be satisfying as an example separating $\cQ_2$-free and $\cQ_3$-free, because it will be important to restrict to objects that are small colimits of representable functors (see Definition~\ref{defn:small}) and $X=(-)^{**}$ is not such. Therefore we present below a more down to earth (and small) example.
	
\begin{example}
	Consider the small functor $X: V \mapsto V^{\bN}/V^{\oplus \bN}$. To show $(\cN_{2} X) V=0$ first observe that $(v_n)_{n \in \bN} + V^{\oplus \bN} \in (\cN_{2} X) V$ if and only if
	\[\bigcup_{N \in \bN} \Ann\{v_n : n \ge N\}= V^*.\]
	In fact, $(v_n)_{n \in \bN} + V^{\oplus \bN} \in \Kr (X \xi)$ if and only if $(\xi \cdot v_n)_{n \in \bN} \in \bfk^\bN$ is eventually null, that is: if there is $N \in \bN$ such that $\xi \in \Ann\{v_n: n \ge N\}$. Now if $V$ is finite dimensional then the only way $(v_n)_{n \in \bN} + V^{\oplus \bN} \in (\cN_2 X) V$ can happen, is that $H_N=:\Ann\{v_n : n \ge N\}$ is eventually $V^*$, but then $v_n$ is eventually $0$. On the other hand if $V$ is infinite dimensional and $H_N$ is not eventually $V^*$, we can then extract a cofinal linearly independent subsequence of $(v_n: n \in \bN)$ and there is a $\xi$ not annihilating any element of such a subsequence.\\
	On the other hand $X \in \cQ_3$ and $X \neq 0$ so it is obvious that $\Nat(X,X) \neq 0$.
\end{example}
	
\begin{example}\label{ex:weirdo}
	Let $V$ be an infinite dimensional space. There is a canonical embedding $V \otimes_\bfk - \subseteq \Vect(V^*, -)$. Its cokernel, the quotient $X= \Vect(V^*,-)/(V \otimes_\bfk -)$ is $\cQ_3$-free, but not $\cQ_4$-free. To see that $X$ is $\cQ_3$-free: assume $f\in \Vect(V^*, \bfk^{\oplus \kappa})$ is such that for every $i\in \kappa$ there is $v_i \in V$ such that
	\[\forall \xi \in V^*, \qquad \delta_i \cdot f \cdot \xi = \xi \cdot v_i.\]
	Then the set $S=\{i \in \kappa : v_i \neq 0\}$ is finite, for otherwise there would be a $\xi\in V^*$ such that $\xi \cdot v_i \neq 0$ for infinitely many $i<\kappa$ contradicting the fact that $\{i \in I: \delta_i (f\cdot \xi) \neq 0\}$ is finite for every $\xi$. From this it follows that for every $\xi \in V^*$ one has $f(\xi) = \sum_{i \in S} e_{i} \cdot \xi(v_i)$.
		
	That $X$ is not $\cQ_4$-free, follows from Proposition~\ref{Q4-free_char} and the fact that $V \otimes_\bfk -$ is not representable when $V$ is infinite-dimensional.
\end{example}

\subsection{Strong vector spaces are small functors}\label{ssec:small}
	
The category $\Sigma\Vect$ is by construction the torsion free part of the pre-torsion theory $\big({}^\bot \Free(\cQ_3^\bot), \Tors({}^\bot\Free(\cQ_3^\bot))^\bot\big)$ on the (non-locally small) Abelian category $[\Vect, \Vect]$.
We can however show that $\Sigma\Vect$ is a torsion free part for a torsion theory on a locally small Abelian category (Theorem~\ref{cor:SigmaVect_is_reflective}): in this subsection we show that $\Sigma\Vect$ is contained in the (locally small) full subcategory of \emph{small} (see Definition~\ref{defn:small}) $\bfk$-additive functors $\Vect\to\Vect$ and natural transformations among them. In the following subsection $\kappa$ will always be a small infinite cardinal $\kappa$. We will denote by $\Vect_{\kappa}$ the category of vector spaces of dimension strictly smaller than $\kappa$.

\begin{lemma}\label{lem:Lan(X|<k)}
    Let $\iota: \Vect_{\kappa} \subseteq \Vect$ be the inclusion functor and let $\tilde{X}: \Vect_{\kappa} \to \Vect$ be $\bfk$-additive. Then the left Kan extension $\Lan_{\iota}\tilde{X}$ is pointwise and given by
    \[\big(\Lan_{\iota}(\tilde{X})\big)(V) \cong \varinjlim_{\substack{H \le V\\\dim_{\bfk}(H)<\kappa}}X(V).\]
    Moreover $\Lan_\iota(\tilde{X}) \circ \iota \cong \tilde{X}$.
    \begin{proof}
        Since $\Vect_{\kappa}$ is essentially small and $\Vect$ is small-cocomplete, by \cite[Ch.~X, Sec.~5, Thm.~3]{maclane71categories}, for every $X:\Vect \to \Vect$
        \[\Lan_\iota(\tilde{X})(V)\cong \varinjlim_{\substack{h: H \to V\\ \dim_{\bfk}(H)<\kappa}} \tilde{X}(H)\]
        is a pointwise Kan extension. Thus to prove the statement it suffices to show that the category $\Sub_\kappa(V)$ of subspaces of $V$ of dimension $<\kappa$ is final in the category $(\iota \downarrow V)$. For this it suffices to show that for every object $h\in \Vect(H, V)$ of $(\iota \downarrow V)$, the category $(h \downarrow \Sub_\kappa(V))$ is connected. In fact we will show it has a terminal object. For this note that we can factor $h$ as $\bar{h}\cdot i$ where $\bar{h}=\coimg(h)$ and $i=\img(h): \Img(h) \le V$ is the inclusion of the image of $h$. Then the arrow $\bar{h} : (H, h) \to (\Img(h), i)$ gives a terminal object in $(h \downarrow \Sub_\kappa(V))$.

        The ``moreover'' follows from the fact that for each $V$ of dimension $< \kappa$, $\Sub_\kappa(V)$ has a terminal object.
    \end{proof}
\end{lemma}
	
\begin{lemma}\label{lem:smallness_crit}
	For $X: \Vect \to \Vect$ $\bfk$-additive and $\kappa$ a regular cardinal, the following are equivalent:
	\begin{enumerate}
		\item $XV \cong \varinjlim_{d \in D} \Vect(Fd, V)$ for some small directed poset $(D, \le)$ and some diagram $F: D \to \Vect^\op$ with $\dim Fd <\kappa$ for every $d \in D$;
		\item for every $V \in \Vect$, one has
		\[XV \cong \varinjlim_{\substack{H \le V\\ \dim H < \kappa}} XH;\]
		\item said $(\Yon \downarrow X)_{\kappa} \subseteq (\Yon \downarrow X)$ the full subcategory with objects $(V, x : \Yon V \to X)$ with $\dim V <\kappa$, the inclusion $(\Yon \downarrow X)_{\kappa} \subseteq (\Yon \downarrow X)$ is \emph{final}\index[ind]{functor!final} (cf \cite[Ch.~ IX, Sec.~3]{maclane71categories});
		\item $X$ is the left Kan extension of its restriction $X|_{\Vect_{\kappa}}$ to the category of spaces of dimension $<\kappa$. 
        \item $X$ is a left Kan extension of some $\bfk$-additive functor $\tilde{X}: \Vect_{\kappa} \to \Vect$.
	\end{enumerate}
	\begin{proof}
		$(1)\Rightarrow (2)$ Let $X \cong \varinjlim_{d \in D} \Vect(Fd, V)$ for some $F: D \to \Vect^\op$ with $\dim Fd <\kappa$ for every $d \in D$. Then
		\[\begin{aligned}
		\varinjlim_{\substack{H < V\\ \dim H < \kappa}} XH &\cong \varinjlim_{\substack{H < V\\ \dim H < \kappa}} \varinjlim_{d \in D} \Vect(Fd, H) \cong\\
		&\cong \varinjlim_{d \in D} \varinjlim_{\substack{H < V\\ \dim H < \kappa}} \Vect(Fd, H) \cong \varinjlim_{d \in D} \Vect(Fd, V)
		\end{aligned}\]
		where the last isomorphism holds because every map from $Fd$ factors through some $H<V$ with $\dim H < \kappa$.
			
		$(2) \Rightarrow (3)$ It suffices to show that for every object $x$ in $(\Yon \downarrow X)$ the slice category $(x \downarrow (\Yon \downarrow X)_{<\kappa})$ is connected.
        In the proof we will use Yoneda to identify the category $(\Yon \downarrow X)$ with the category of elements of $X$, that is we will identify the objects of $(\Yon \downarrow x)$ with pairs $(V, \tilde{x}\in XV)$. The arrows $f:(V_0, \tilde{x}_0 \in X(V_0)) \to (V_1, \tilde{x}_1 \in X(V_1))$ will be maps $f: V_0 \to V_1$ such that $Xf \cdot x_0= x_1$.
        
        First observe that $(2)$ tells us that $(x \downarrow (\Yon \downarrow X)_{<\kappa})$ is nonempty. Let $x:=(V, \tilde{x} \in X(V))$ be an element of $X$, by $(2)$ there is some $H < V$ of dimension $< \kappa$ with $\tilde{x} \in XH \le XV$, set $z=(H, \tilde{x} \in X(H))$ then the inclusion $H < V$ gives us a map in $(\Yon \downarrow X)(x, z)$. Finally note that $(x \downarrow (\Yon \downarrow X)_{<\kappa})$ has binary coproducts and coequalizers, so being non empty, it is connected.
			
		$(3) \Rightarrow (1)$ take $D=(\Yon \downarrow X)_{<\kappa}$ and $F$ the projection on the first factor.
			
		$(4) \Leftrightarrow (2)$ and $(5) \Leftrightarrow (4)$ follow from Lemma~\ref{lem:Lan(X|<k)}. 
	\end{proof}
\end{lemma}
	
\begin{definition}\label{defn:small}
	A $\bfk$-additive functor $X$ is \emph{$\kappa$-ary} if one of the above conditions holds. If $X$ is $\kappa$-ary for some small $\kappa$, then it is \emph{small}. The category of small $\bfk$-additive functors and $\bfk$-additive natural transformations is the full subcategory $\Ind\-(\Vect^\op)$ of $[\Vect, \Vect]$. With a slight abuse of notation we will denote by $\Ind\-(\Vect_{\kappa}^\op)$ the category of $\kappa$-ary $\bfk$-additive endofunctors of $\Vect$.
\end{definition}
 
\begin{remark}
	Strictly speaking, $\Ind\-\cC$ is the full sub-category of the category $\Func(\cC^\op, \Sets)$ consisting only of the $\Sets$-valued functors $\cC^{\op}\to \Sets$ which are small filtered colimits of representables. However, a representable $\Vect(H,-) : \Vect \to \Sets$ factors uniquely through the forgetful functor $U: \Vect  \to \Sets$ as a $\bfk$-enriched $\Vect(H,-) : \Vect \to \Vect$, and since the forgetful functor $U$ creates filtered colimits (see e.g.\ \cite[Ch.~IX, Sec.~2, Thm.~2]{maclane71categories}), this holds as well for every small filtered colimit of representables. Thus the postcomposition by $U$
	\[U_*: [\Vect, \Vect] \to \Func(\Vect, \Sets),\] restricts to an equivalence of categories from the full subcategory of small $\bfk$-additive functors in $[\Vect, \Vect]$ to 
	$\Ind\-(\Vect^{\op})$. We will regard $\Ind\-(\Vect^\op)$ as genuinely included in $[\Vect, \Vect]$.
\end{remark}
	
\begin{remark}\label{rmk:IndVect_and_(co)limits}
	The category $\Ind\-(\Vect^{\op})$ is locally small and Abelian (see \cite[Thm.~8.6.5~i]{kashiwara2005categories}). For a more detailed presentation of the properties of the indization and its enriched structure in the case of Abelian categories such as $\Vect$ we refer to \cite{kashiwara2005categories}.
	In particular we recall that
	\begin{enumerate}
		\item $\Ind\-(\Vect^\op)$ has all small limits and the inclusions
		\[\Ind\-(\Vect^\op) \subseteq [\Vect, \Vect] \subseteq \Func(\Vect, \Sets)\]
		are fully faithful and create limits (\cite[Prop.~6.1.16]{kashiwara2005categories} and the fact that the forgetful functor $U: \Vect \to \Sets$ creates limits).
		\item $\Ind\-(\Vect^\op)$ has all small colimits and the inclusion
		\[\Ind\-(\Vect^\op)\subseteq [\Vect, \Vect]\]
		preserves them (follows from items ii and v of \cite[Proposition~8.6.5]{kashiwara2005categories}).
	\end{enumerate}
\end{remark}

\begin{remark}\label{rmk:IndVect_well_powered}
    $\Ind\-(\Vect^\op)$ is \emph{well powered} (and thus also \emph{well-copowered} as it is Abelian). This is because a subobject of a $\kappa$-ary object is $\kappa$-ary as well: to see this use for example \cite[Lem.~8.6.7]{kashiwara2005categories}, to write a monomorphism $f: X \to Y$ as a colimit of monomorphisms of representables and recall that if $\Yon g : \Yon V \to \Yon W$, then $g:W \to V$ is epi and $\dim (V) \le \dim W$.
\end{remark}

\begin{remark}\label{rmk:k-ary_form_coreflective_subcat}
    For every infinite cardinal $\kappa$, the subcategory $\Ind\-(\Vect_{\kappa}^\op)\subseteq [\Vect, \Vect]$ is coreflective. In fact if $\iota: \Vect_{\kappa} \subseteq \Vect$ is the inclusion functor and $X: \Vect \to \Vect$ is $\bfk$-additive, then the left Kan extension
    \[\Lan_{\iota} (X\iota) (V) \cong \varinjlim_{\substack{H \le V\\\dim(H)<\kappa}} X(H)\]
    is $\kappa$-ary and for every $\kappa$-ary $Y\cong \varinjlim_{d \in D} \Yon G (d)$ with $\dim(G(d))<\kappa$, we have a chain of natural isomorphisms
    \[\begin{aligned}\Nat \big(Y , \Lan_{\iota} (X\iota)\big) \cong 
    \varprojlim_{d \in D} \Nat\big(\Yon Gd, \Lan_{\iota} (X\iota)\big)\cong \\
    \cong \varprojlim_{d \in D}\Lan_{\iota} (X\iota)(Gd) \cong 
    \varprojlim_{d \in D} X(G(d)) \cong\\
    \cong \varprojlim_{d \in D} \Nat(\Yon G d, X) \cong \Nat(Y,X).
    \end{aligned}\]
    Here the first and last isomorphism are given by the continuity of $\Nat(-,-)$, the second is by Yoneda, and the third is because $\Lan_\iota(X\iota)\circ \iota \cong X\iota$ (Lemma~\ref{lem:Lan(X|<k)}).
\end{remark}
    
\begin{remark}
	The category $\{X: \cN_{2} X=0\}$ of $\bfk$-concrete (equiv.\ $\cQ_2$-free) $X$ is necessarily locally small (every hom-set has an injection into a small vector space), but is not a subcategory of $\Ind\-(\Vect^{\op})$: for example the double-dualization functor
	\[(-)^{**}:=\Vect(\Vect(-,\bfk), \bfk) : \Vect \to \Vect\]
	is $\cQ_2$-free (Example~\ref{ex:doubledual}), but is not small.
\end{remark}
	
The fact that $\Sigma\Vect\subseteq \Ind\-(\Vect^\op)$ will be a corollary of the following
	
\begin{lemma}
	If $X \in \Sigma\Vect$ then $X$ is $\kappa$-ary if and only if for every set $I$, every summable $I$-family $(f_i)_{i \in I}$ in $X\bfk$ is such that $|\{i \in I: f_i \neq 0\}|<\kappa$.
	\begin{proof}
		Let $I, J_0, J_1$ be sets, $I = J_0 \sqcup J_1$.
		There is a natural monomorphism between the two split exact sequences (in fact between biproducts).
		\[\begin{tikzcd}
             0 \ar[r] &(-)^{\oplus J_0}\ar[d]\ar[r] &(-)^{\oplus I} \ar[d]\ar[r] &(-)^{\oplus J_1}\ar[d]\ar[r]&0\\
			0 \ar[r] &(-)^{J_0}\ar[r] &(-)^{I}\ar[r] &(-)^{J_1}\ar[r] & 0
		\end{tikzcd}\]
		Applying $\Nat(-,X)$ to the above diagram results in the following diagram in $\Vect$ with exact rows (the rows are in fact split-exact so exactness is preserved)
		\[\begin{tikzcd}
			0 \ar[r] &(X\bfk)^{J_0}\ar[r] &(X\bfk)^{I} \ar[r] &(X\bfk)^{J_1} \ar[r] &0\\
			0 \ar[r] &X(\bfk^{\oplus J_0})\ar[r]\ar[u] &X(\bfk^{\oplus I})\ar[r]\ar[u] &X(\bfk^{\oplus J_1})\ar[u] \ar[r] &0
		\end{tikzcd}\]
		If $X\in \Sigma\Vect$ the vertical arrows are monomorphisms. It follows that the leftmost square is a pullback. 
			
		Thus, given any $X \in \Sigma\Vect$ and sets $J \subseteq I$ we have that the diagram
		\[\begin{tikzcd}
			X(\bfk^{\oplus J}) \ar[r, hook]\ar[d, hook]
			&(X\bfk)^J\ar[d, hook]\\
			X(\bfk^{\oplus I}) \ar[r,hook] &(X\bfk)^I,
        \end{tikzcd}\]
		is an intersection square.\\
		Since in $\Vect$ filtered colimits preserve pullbacks, we get that for any $\kappa$, the following is a pullback as well
		\[\begin{tikzcd}
			\displaystyle{\bigcup_{J \in [I]^{<\kappa}}X(\bfk^{\oplus J})} \ar[r]\ar[d]
			&\displaystyle{\bigcup_{J \in [I]^{<\kappa}}(X\bfk)^J}\ar[d]\\
			X(\bfk^{\oplus I}) \ar[r] &(X\bfk)^I.
		\end{tikzcd}\]
		But if every summable $I$-family $(f_i)_{i \in I}$ in $X\bfk$ is such that $|\{i \in I: f_i \neq 0\}|<\kappa$ then the image of the lower horizontal monomorphism lies in the image of the rightmost vertical monomorphism, from which it follows that the leftmost vertical arrow is in fact an isomorphism. This in turn implies by Lemma~\ref{lem:smallness_crit} that $X$ is $\kappa$-ary. The converse is trivial. 
	\end{proof}
\end{lemma}
	
\begin{corollary}\label{cor:strlin_are_small}
	If $X \in \Sigma\Vect$ then $X$ is small.
	\begin{proof}
		Suppose $X \in \Sigma\Vect$, $I$ is a set, $(f_i)_{i \in I}$ is summable in $X\bfk$ and $f_i= c \in X\bfk\setminus \{0\}$ for every $i\in I$; then $I$ is finite.
		In fact if $I$ was infinite then we would have
		\[\sum_{i \in I} f_i = \sum_{i \in I} f_i + \sum_{i \in I} f_i,\]
		so $\sum_{i \in I} f_i = 0$, but this would hold also for any $I \setminus \{i\} $ so it would have to be that $f_i = 0$.
		It follows that if $\kappa = |X\bfk|^{++}$, for every summable $I$-family $(f_i)_{i \in I}$ one has that $|\{i \in I: f_i \neq 0\}|<\kappa$; for otherwise, there would be $c \in X\bfk\setminus \{0\}$ such that $J_c=\{i \in I: f_i=c\}$ is infinite and since $(f_i)_{i \in J_c}$ is summable, this is impossible. 
	\end{proof}
\end{corollary}

$\Ind\-(\Vect^\op_\bfk)$ is a small-bicomplete locally small Abelian category. It is a general phenomenon that under enough completeness assumptions on the ambient category, a right-orthogonal class of objects forms a reflective subcategory: this is the case for example if the category is locally presentable \cite[Thm.~1.39]{adamek1994locally} (see also \cite[Thm.~4.1.3]{freyd1972categories} for another earlier results along this lines). 
Unfortunately $\Ind\-(\Vect^\op)$ is not locally-presentable, because it lacks a \emph{small} dense subcategory. 
We can nevertheless derive reflectiveness of $\cQ$-free objects from \cite[Thm.~1.39]{adamek1994locally} by reducing to reflections defined on the finitely presentable $\Ind\-(\Vect_{\kappa}^\op)$, when the class $\cQ=\cQ(\bE)$ consists of quotients of representables as in Definition~\ref{defn:cQ()}.

\begin{theorem}\label{thm:Q-free_are_reflective}
    If $\cQ=\cQ(\bE)$ is a class of quotients in representables as in Definition~\ref{defn:cQ()}, then $\Free(\cQ^\bot)$ is an epi-reflective subcategory of $\Ind\-(\Vect^\op)$. Consequently it is the torsion-free part for a torsion theory on $\Ind\-(\Vect^\op)$.
    \begin{proof}
        Let $X$ be $\kappa$-ary, so $X \in \Ind\-(\Vect_{\kappa}^\op)\subseteq \Ind\-(\Vect^\op)$.
        Without loss of generality we can assume that $\cQ$ is closed under taking quotients: this is possible because if $\cQ'$ is the class of quotients of objects in $\cQ$, then $\Free(\cQ^\bot)=\Free(\cQ'^\bot)$.
        With this assumption we see that $X$ is $\cQ$-free if and only if it is $\cQ^\kappa$ free where $\cQ^\kappa$ consists of the $\kappa$-ary objects in $\cQ^\kappa$.
        In fact, let $X$ be $\cQ^\kappa$-free and $\kappa$-ary and let $Q$ be in $\cQ$, so in particular $Q=\Ckr(a)$ for some map with representable codomain $a:F \to \Yon V$. Then, by $\kappa$-arity of $X$, Lemma~\ref{lem:smallness_crit}~(2), and Yoneda, every $x=\tilde{x} \cdot \ckr(a) : \Yon V \to X$ factors as $x=y \cdot (\Yon \pi_H)$ for some $H \le V$ with $\dim(H) <\kappa$ and some $y: \Yon H \to X$. But then since $\Img(y)\cong\Coimg(y)$ is naturally a quotient of $Q$ and since $\cQ$ is closed under taking quotients we have $\Img(y) \in \cQ^\kappa$. Thus, by $\cQ^\kappa$-freeness of $X$, $y=0$ and hence $\tilde{x}=0$.
        This implies that if $X$ has a an epi-reflection in the category $\Free\big((\cQ^{\kappa})^{\bot}\big)\subseteq \Ind\-(\Vect_{\kappa}^\op)$, then this reflection is an epi-reflection in $\Free(\cQ^\bot)$. In order to conclude, it is therefore enough to show that for all small cardinals $\kappa$, $\Free(\cQ^{\kappa,\bot})$ is an epi-reflective subcategory of $\Ind\-(\Vect_{\kappa}^\op)$, but this follows from the familiar fact that orthogonal subcategory are epi-reflective in a locally finitely presentable category as $\Ind\-(\Vect_{\kappa})$ \cite[Thm.~1.39]{adamek1994locally}.

        As for $\big({}^{\bot}(\cQ^\bot), \cQ^\bot\big)$ being a torsion theory, by \cite[Cor.~3.4]{cassidy1985reflective}, it suffices to observe that $\Ind\-(\Vect^\op)$ is finitely well complete and $\Free(\cQ^\bot)$ is ${}^{\bot}(\cQ^\bot)$-reflective.
        
        Finite well-completeness of $\Ind\-(\Vect^\op)$ readily follows from Remarks~\ref{rmk:IndVect_and_(co)limits} and \ref{rmk:IndVect_well_powered}.
        As for ${}^{\bot}(\cQ^\bot)$-reflectivivity of $\Free(\cQ^\bot)$, note that ${}^{\bot}(\cQ^\bot)$ contains all epimorphisms and $\Free(\cQ^\bot)$ is epi-reflective.

        The last assertion follows from Proposition~\ref{prop:TorsionTheoriesAbelian}.
    \end{proof}
\end{theorem}

\begin{corollary}\label{cor:SigmaVect_is_reflective}
    $\Sigma\Vect \subseteq \Ind\-(\Vect^\op)$ is reflective and the torsion free part for a torsion theory on $\Ind\-(\Vect^\op)$.
\end{corollary}

\begin{notation}\label{not:reflections}
    We will use $\cR_\cQ^\infty: \Ind\-(\Vect^\op) \to \Ind\-(\Vect^\op)$ to denote the endofunctor of $\Ind\-(\Vect^\op)$ given by the composition of the reflection onto $\Free(\cQ^\bot)$ with the inclusion $\Free(\cQ^\bot)\subseteq \Ind\-(\Vect^\op)$, and $\rho_\cQ^\infty: \Id \to \cR_\cQ^\infty$ to denote the reflection transformation.
    The well pointed endofunctor $(\cR_\cQ^\infty, \rho_\cQ^\infty)$ can be explicitely constructed as a transfinite iteration (in the sense of \cite[Sec.~5]{kelly1980unified}) of the well pointed endofunctor $\rho_\cQ: \Id \to \cR_\cQ: \Ind\-(\Vect^\op) \to \Ind\-(\Vect^\op)$ given by
    the natural quotients
    \[\rho_\cQ (X):  X \to \cR_\cQ(X):= X/\cN_\cQ(X)\]
    where $\cN_\cQ(X)$ is the smallest subobject of $X$ containing the image of all maps $Q \to X$ for all $Q \in \cQ$.
    When $\cQ \in \{\cQ_2, \cQ_3, \cQ_4\}$ we will use $\cR_i^\infty$ and $\rho_i^\infty$ instead of $\cR_{\cQ_i}^\infty$, and $\rho_{\cQ_i}^\infty$.
\end{notation}

\subsection{\texorpdfstring{$\Sigma$}{Σ}-embeddings}\label{ssec:SigmaEmb}
In this subsection we study monomorphisms of $\Sigma\Vect$ which reflect infinite sums. 
	
It makes sense to distinguish a more general class of monomorphism between objects in $\Ind\-(\Vect^\op)$ with the property that their cokernel in $\Ind\-(\Vect^\op)$ is in $\Sigma\Vect$. Proposition~\ref{prop:SigmaEmb&Quot} will show that when restricted to arrows in $\Sigma\Vect$ such a class consists precisely of the monomorphisms which reflect infinite sums.
	
Since the generality comes cost-free and it will have some applications, in the remainder of this chapter we develop a slightly more general theory for $\cQ$-free objects for $\cQ=\cQ(\bE)$ with $\bE$ as in Definition~\ref{defn:cQ()}.

\begin{definition}\label{def:AQ} 
	The arrows of $\Ind\-(\Vect^\op)$ that are in $\cM_3:=(\cA_3 \cup \Epimorphisms)^\bot$ or, equivalently, the monomorphisms which are right orthogonal to the arrows $\delta_\kappa : (\Yon\bfk)^{\oplus \kappa} \to (\Yon \bfk)^{\kappa}$ for every $\kappa \in \Card$, will be called \emph{$\Sigma$-closed embeddings}.
	The arrows of $\Ind\-(\Vect^\op)$ that are in $\cM_2:=(\cA_2 \cup \Epimorphisms)^\bot$ or, equivalently, the monomorphisms which are right orthogonal to the natural inclusion $(V^*\otimes_\bfk -) \to \Vect(V,-)$ for every $V \in \Ob(\Vect)$, will be called \emph{embeddings}.
\end{definition}

\begin{remark}
	Propositions~\ref{prop:SigmaEmb&Quot} and \ref{prop:emb_quot} below will show that the definitions of \emph{$\Sigma$-closed embedding} and \emph{embedding} given in Definition~\ref{def:AQ} agree with and extend the ones of Remark~\ref{rmk:Sigma_emb_B}.
\end{remark}
	
\begin{lemma}\label{lem:toname1}
	Assume $\mathcal{A}$ is a class of monomorphisms in an Abelian category $\cC$, each with projective codomain, $\cQ=\{\Ckr\, a: a \in \mathcal{A}\}$, and 
	\[\begin{tikzcd}
		0 \ar[r]&X \ar[r,"f"] \ar[r] &Y \ar[r, "p"] &C \ar[r] &0
	\end{tikzcd}\]
	is exact in $\cC$. Then $f$ is right-orthogonal to every $a \in \mathcal{A}$ if and only if $C$ is $\cQ$-free.
	\begin{proof}
        Let $a \in \cA$ with domain $F$ and projective codomain $P$ and let $\pi = \ckr(a): P \to Q=\Ckr(a)$.
        Note that any square $(\alpha, \beta) : a \to f$, induces a unique arrow $\gamma: Q \to C$, making $(\alpha, \beta, \gamma)$ a morphism of short exact sequences of $\cC$
        \[\begin{tikzcd}
		0 \ar[r]&F \ar[r,"a"]\ar[d,"\alpha"] \ar[r] &P \ar[r, "\pi"] \ar[d, "\beta"] &Q \ar[r, "p"] \ar[d, "\gamma"] &0\\
        0 \ar[r]&X \ar[r,"f"] \ar[r] &Y \ar[r, "p"] &C \ar[r] &0.
	    \end{tikzcd}\]
        
        Suppose that $Q\ort C$ so $\cC(Q, C)=0$: then $\gamma=0$ and the image of $\beta$ is contained in the image of $f$, whence $a \ort f$.
        
        Conversely, suppose $a \ort f$. If $\gamma \in \cC(Q, C)$ is any arrow, then, since $P$ is projective, $\gamma \cdot \pi: P \to C$ lifts to an arrow $\beta: P \to Y$ such that $p \cdot \beta = \gamma \cdot \pi$. Since then $p\cdot \beta \cdot a=0$, it follows that $\beta \cdot a$ is in the image of $f$ and thus induces an arrow $\alpha: F \to X$.
        By $a \ort f$, it then follows that actually the image of $\beta$ is included in the image of $f$ and thus it must be $\gamma=0$.
	\end{proof}
\end{lemma}

\begin{proposition}\label{prop:SigmaEmb&Quot}
	The following are equivalent for a monomorphism $X \le Y$ in $\Sigma\Vect$:
	\begin{enumerate}
		\item $X\le Y$ is a $\Sigma$-closed embedding;
		\item for every set $I$, every family $(f_i)_{i \in I}\in (X\bfk)^I$ which is summable in $Y\bfk$ is summable in $X\bfk$ (in which case the sums agree);
		\item the quotient $Y/X$ in $\Ind\-(\Vect^{\op})$ is again in $\Sigma\Vect$.
	\end{enumerate}
	\begin{proof}
		That $j:X \le Y$ is right-orthogonal to the inclusion $\delta_I : (\Yon\bfk)^{\oplus I} \le (\Yon \bfk)^I$, means that for every arrow $f: (\Yon \bfk)^I \to Y$ such that $\Img(f\cdot \delta_\lambda) \le X$ we also have $\Img(f) \le X$. Remembering that arrows from $(\Yon \bfk)^I\cong \Yon(\bfk^{\oplus I})$ correspond to summable families we see the equivalence $(1) \Leftrightarrow (2)$.

		The equivalence $(2)\Leftrightarrow (3)$ follows from Lemma~\ref{lem:toname1}.
	\end{proof}
\end{proposition}

\begin{remark}[Summable families in quotients]
	If $X \le Y \in \Sigma\Vect$ is a $\Sigma$-closed embedding, then it follows that summability on $Y/X$ behaves as follows: summable families in $(Y/X)\bfk \cong Y\bfk/X\bfk$ are given by families $(y_i+ X\bfk)_{i \in I}$ for which there is a family $(x_i)_{i\in I} \in (X\bfk)^I$ such that $(y_i+x_i)_{i \in I}$ is summable in $Y$, in which case the sum is just 
	\[\sum_{i \in I} (y_i + X\bfk) =\big(\sum_{i \in I} y_i + x_i\big) + X\bfk.\]
\end{remark}
		
\begin{proposition}\label{prop:emb_quot}
	Let $X \le Y$ be a monomorphism and $Y \in \Sigma\Vect$. Then the following are equivalent:
	\begin{enumerate}
		\item for every set $I$, a family $(x_i)_{i \in I} \in (X\bfk)^I$ is summable if and only if it is summable in $Y\bfk$ and, if this happens, then for every $(k_i)_{i \in I} \in \bfk^I$ one has $\sum_{i \in I} x_i \cdot k_i \in X\bfk$;
		\item $X\le Y$ is an embedding;
		\item for every $V \in \Vect$ the following square is Cartesian
		\[\begin{tikzcd}
			XV \ar[rr, hook] \ar[d, "{(X-)\cdot -}"'] &&YV\ar[d, "{(Y-)\cdot -}"]\\
			\Vect_(V^*, X\bfk) \ar[rr] && \Vect(V^*, Y\bfk);
		\end{tikzcd}\]
		\item the quotient $Y/X$ is $\bfk$-concrete (i.e.\ $\cQ_2$-free). 
	\end{enumerate}
	\begin{proof}
		The condition $(2)$ that the inclusion $j:X \le Y$ is right orthogonal to every inclusion $i:\Fin (V^*) \le \Yon(V)$ means that for every arrow $f:\Yon(V) \to Y$ such that $\Img(f \cdot i)\le X$, there is a unique $\tilde{f}: \Yon (V) \to X$ such that $\tilde{f} \cdot j = f \cdot i$.
        Note that $\Img(f \cdot i)\le X$ is equivalent to $\Img(f\bfk) \le X\bfk$. Thus $(2)$ is equivalent to the requirement that every map $f: \Yon(V) \to Y$ such that the image of $f\bfk$ is in $X\bfk$, lifts to a map $\tilde{f}: \Yon(V) \to X$. Recalling that every $\Yon(V)$ is isomorphic to $(\Yon\bfk)^I$ for some set $I$ and that arrows from $(\Yon \bfk)^I$ correspond to summable families we have the equivalence $(1)\Leftrightarrow (2)$.
  
        $(2)\Leftrightarrow (3)$. Recall that orthogonality of $j: X \le Y$ to $i: \Fin(V^*) \le \Yon V$ is equivalent to the fact that the square
        \[\begin{tikzcd}
            \Nat(\Yon V, X) \ar[d, "{\Nat(i, X)}"'] \ar[rr,"{\Nat(\Yon V, j)}"] &&\Nat(\Yon V, Y)\ar[d, "{\Nat(i, Y)}"]\\
            \Nat(\Fin V^*, X) \ar[rr, "{\Nat(\Fin V^*, j)}"] &&\Nat(\Fin V^*, Y)
        \end{tikzcd}\]
        is Cartesian. The square in $(3)$ is precisely this square up to the natural isomorphisms $\Nat(\Yon V, X) \cong X(V)$ and $\Vect(V^*, X\bfk) \cong \Nat(V^* \otimes- , X)$.
        
		$(2)\Leftrightarrow (4)$ readily follows from Lemma~\ref{lem:toname1}.
	\end{proof}
\end{proposition}

\section{Comparisons}\label{sec:Comp}
	
In this section we compare $\Sigma\Vect$ with other more easily understood categories.
It is easy to see that $B\Sigma \Vect$ is a reasonable category of strong vector spaces: the first subsection explicitly describes the unique fully faithful functor $B\Sigma\Vect \to \Sigma\Vect$ given by \ref{SV:U}.
	
The second subsection relates $\Sigma\Vect$ to the category $\TVect$ of linearly topologized vector spaces. The relation of $[\Vect, \Vect]$ and $\TVect$ should be a well known consequence of a result of Lefschetz and the application of the nerve construction. It consists of two main parts:

\begin{enumerate}
	\item There is a pair of adjoint functors between $\TVect$ and $[\Vect, \Vect]$ whose induced comonad on $\TVect$ is idempotent and whose coalgebras are the $K$-spaces, i.e.\ the linearly topologized vector spaces which are ``generated'' by \emph{linearly compact spaces}.
	\item The coreflective category of separated $K$-spaces is equivalent to the category of $\cQ_4$-free objects in $[\Vect, \Vect]$; these are automatically small and form a reflective subcategory of $\Ind\-(\Vect^\op)$ (in fact of $\Sigma\Vect$).
\end{enumerate}
	
\subsection{General versus based strong vector spaces}\label{ssec:Comp_B} This subsection explicitely describes the fully faithful embedding of $B\Sigma\Vect$ as described in Section~\ref{sec:BasSL} into $\Sigma\Vect$.

	
\begin{definition}\label{def:BasedToUnbased}
	Let $(\Gamma, \cF)$ be a bornological set. Inside of $\bfk^{\oplus \Gamma}$ consider the filter of subspaces $U(\Gamma,\cF)$ generated by the family $\{\bfk^{\oplus(\Gamma \setminus F)}: F \in \cF\}$. This induces a filtered diagram $F_{(\Gamma, \cF)}: U(\Gamma, \cF)^{\op} \to \Vect^{\op}$ given by $F_{(\Gamma, \cF)}(H) = \bfk^{\oplus \Gamma}/H$ (with the obvious projection maps). Define 
	\[\bfk[\Gamma; \cF] = \varinjlim_{H \in U(\Gamma, \cF)} \Yon F_{(\Gamma, \cF)} H \cong \varinjlim_{F\in \cF} \Vect(\bfk^{\oplus \Gamma}/\bfk^{\oplus (\Gamma \setminus F)}, -)\in [\Vect, \Vect]\]
	and note that canonically $(\bfk[\Gamma; \cF])(\bfk) \cong \bfk(\Gamma; \cF)$ as a $\bfk$-vector space.
	More concretely $\bfk[\Gamma; \cF]$ is canonically isomorphic to
	\[\bfk[\Gamma; \cF](V)\cong \bigcup_{F \in \cF} V^F \subseteq V^\Gamma\]
	where $V^F$ sits in $V^\Gamma$ as $V^F \times \{0\}^{\Gamma\setminus F}$.
\end{definition}
	
\begin{proposition}
	If $(\Gamma, \cF)$ is a bornological set, then $\bfk[\Gamma; \cF]$ is a strong vector space. Moreover a family $(f_i : i \in I)$ in $\bfk(\Gamma; \cF)=\bfk[\Gamma;\cF](\bfk)$ is summable with respect to Definition~\ref{def:summable1} if and only if it is summable with respect to Definition~\ref{defn:summable2}, and the sums coincide.
	\begin{proof}
		Any $x \in \bfk[\Gamma; \cF](\bfk^{\oplus I})\cong\bigcup_{F \in \cF} (\bfk^{\oplus I})^F$ is inside some $(\bfk^{\oplus I})^F$. Now 
		\[\big(\bfk[\Gamma; \cF]\big)\delta_i \cong \bigcup_{F \in \cF} \delta_i^F: \bigcup_{F \in \cF} (\bfk^{\oplus I})^F \to \bigcup_{F \in \cF} \bfk^{F},\]
		so if $\big(\bfk[\Gamma; \cF]\delta_i\big)\cdot x=0$ for every $\delta_i$, we must have $x=0$. This suffices to deduce $\bfk[\Gamma, \cF]$ is $\cQ_3$-free and hence in $\Sigma\Vect$.
		
		Now, a family $(f_i: i \in I)\in \bfk(\Gamma; \cF)$ is summable with respect to Definition~\ref{def:summable1}, if and only if there is $F \in \cF$ such that for every $i \in I$, $f_i\in \bfk^F$ and $(f_i:i \in I)$ is summable as a subfamily of $\bfk^F$. This is in turn equivalent to the existence of a (unique) linear function $f:\bfk^{\oplus F} \to \bfk^{\oplus I}$ such that $f^*: \bfk^{I} \to \bfk^{F}$ has the property that $f^*\cdot \delta_i=f_i$.
		
		Now set $x$ to be the image of such $f$ for the natural map
		\[\Vect(\bfk^{\oplus F}, \bfk^{\oplus I})\cong (\bfk^{\oplus I})^F\subseteq \bigcup_{F \in \cF}(\bfk^{\oplus I})^{F}\cong \bfk[\Gamma; \cF](\bfk^{\oplus I}).\]
		With such $x$ we have $f_i = (\bfk[\Gamma;\cF]  \delta_i) \cdot x$ and $\sum_{i \in I} f_i$ as defined in Definition~\ref{def:summable1} equals $(\bfk[\Gamma;\cF]  \sigma) \cdot x$.
		
		A similar argument shows that whenever $(f_i: i \in I)\in \bfk[\Gamma; \cF](\bfk)$ is summable, the corresponding family in $\bfk(\Gamma; \cF)$ is summable.
	\end{proof}
\end{proposition}
	
\begin{corollary}\label{cor:fromBtogen}
	Every $X \in \Sigma\Vect$ is $X\simeq \bfk[\Gamma; \cF]/H$ for some set $\Gamma$, some bornology $\cF$ and some $\Sigma$-closed embedding $H\subseteq \bfk[\Gamma; \cF]$. The choice of the bornological set $(\Gamma,\cF)$ and of the $\Sigma$-closed subobject $H$ is not canonical.
	\begin{proof}
		Every object of $\Sigma\Vect$ is a quotient in $\Ind\-(\Vect^\op)$ of some
		$P=\bigoplus_{i \in I} \Yon V_i$
		for some $\Sigma$-closed embedding $H\subseteq P$. Now choose bases $B_i$ of $V_i$ for each $i \in I$. Let $B= \bigsqcup_{i\in I} B_i$ and let $\cF$ be the ideal of subsets of $B$ whose intersection with each $B_i$ is finite. Then $\bfk[B, \cF]\simeq P$.
	\end{proof}
\end{corollary}

\subsection{Topological vector spaces}\label{ssec:Comp_Top}

$\Vect^{\op}$ embeds fully faithfully in the category $\TVect$ of linearly topologized vector spaces and continuous linear maps (see Proposition~\ref{prop:Lefschetz-Duality} below, essentially due to Lefschetz \cite{lefschetz1942algebraic}). Recall that a \emph{linearly topologized} (short \emph{l.t.}) $\bfk$-vector space is a vector space with a vector topology generated by subspaces. The topology on a l.t.\ $\bfk$-vector space is uniquely determined by the filter of its open subspaces. So for us a l.t. $\bfk$-vector space will be a pair $(V, \cF)$ with $V \in \Vect$ and $\cF$ a filter of subspaces. Note that the topology is T0 (and hence T2 by uniformity) if and only if $\bigcap \cF=0$: that is if the filter is not concentrated on any proper quotient (i.e.\ is not a subfilter of a principal filter).
	
The following notion is from \cite[27.1]{lefschetz1942algebraic}.
	
\begin{definition}\label{def:linearly_compact}
	A topological vector space is said to be \emph{linearly compact} if every family of closed affine subspaces having the finite intersection property has a nonempty intersection.
\end{definition}
	
\begin{proposition}[Lefschetz Duality]\label{prop:Lefschetz-Duality}
	Let $V$ be a vector space and $B$ be a basis of $V$. The subspace topology on $V^*\subseteq \bfk^V$ induced by the product topology on $\bfk^V$ is the same as the topology induced by the isomorphism $V^*\simeq \bfk^B$ and is a linearly compact topology.
	A linear map $V^* \to W^*$ is continuous w.r.t.\ the above topologies if and only if it has the form $f^*$ for some linear map $f: W \to V$.
	\begin{proof}
		This follows from \cite[27.2 and 28.2]{lefschetz1942algebraic}.
	\end{proof}
\end{proposition}
	
		
\begin{definition}
	Define $L: \Vect^{\op} \to \TVect$ to be the faithful functor defined by the assignments in the proposition above, i.e.\ $LV$ is $V^*$ with the linearly compact topology and $Lf$ for $f \in \Vect(V, W)$ is $f^* : W^* \to V^*$.
\end{definition}
	
We will also need
	
\begin{lemma}
	If $V$ is a vector space, then the canonical inclusion $V \subseteq L(V^*)$ has dense image. Moreover the topology induced on $V$ by the inclusion has as filter of open subspaces the filter of finite-codimensional subspaces.
	\begin{proof}
		Open subspaces of $L(V^*)$ are given by finite intersections of subspaces of the form $H_\xi:=\{x \in V^{**}: x(\xi)=0\}$ for $\xi \in V^*$. To show the image is dense, observe that given any finite tuple $(\xi_i: i <n)$ in $V^*$ and any $x \in V^{**}$, there is $v \in V$ such that for all $i$, $x(\xi_i)=\xi_i (v)$. Thus for every $x \in V^{**}$ every neighborhood of $x$ in $L(V^*)$ contains an element of $V\subseteq V^{**}$.
		The claim on the induced topology on $V$ follows from the fact that for every subspace $W$  of $V$, $\Ann(\Ann(W)) \cap V=W$.
	\end{proof}
\end{lemma}
		
A consequence of Proposition~\ref{prop:Lefschetz-Duality} then is
	
\begin{theorem}[Lefschetz]
	The functor $L$ is fully faithful and the essential image of $L$ consists of the linearly compact spaces.
\end{theorem}
	
The nerve construction (see e.g.\ \cite[Rmk.~2.3 and App.~A]{diliberti2020codensity} or \cite[Sec.~2]{leinster2012codensity}) now provides an adjunction 
\[F\dashv \TVect(L, -) : [\Vect, \Vect] \to \TVect\]
where the left adjoint is given by
\[FX \cong \varinjlim_{x: \Yon (H) \to X} LH \cong \int^{H\in \Vect^{\op}} XH \cdot LH\]
Here $\cdot$ is the following $\bfk$-additive version of the copower: $H\cdot S$ is $H \otimes S$ endowed with the finest linear topology making all the inclusions $h \otimes S \to H \otimes S$ continuous as $h$ ranges in $H$. 
To fully justify this we need to show that such large colimit exists in $\TVect$.
	
\begin{lemma}\label{lem:FXexists}
	If $X \in [\Vect, \Vect]$, then $FX$ exists and its supporting space is $X\bfk$.
	\begin{proof}
		Let us identify, as we did in the proof of Lemma~\ref{lem:smallness_crit}, $(\Yon \downarrow X)$ with the category of elements of $X$, and let $P : (\Yon \downarrow X) \to \Vect^\op$ denote the projection on the first factor $(H,x \in X H) \mapsto H$.
  
        Define $FX$ as $X\bfk$ together with the finest topology making all the maps $LH \to X\bfk$ given by $L\xi \mapsto X\xi \cdot x$, continuous as $H$ ranges through $\Vect$ and $x \in XH$. Note that $\eta(H, x): L\xi \mapsto X\xi \cdot x$ as $(H, x \in XH)$ range in $(\Yon \downarrow X)$, then defines a cone $LP \to FX$.

        Let $(\alpha (H,x)  \in \TVect(L H, S): H \in \Vect, x \in X H)$ be a cone for $L P : (\Yon \downarrow X) \to \TVect$.

        We claim $\alpha = f \cdot \eta$ for a unique continuous linear map $f: FX \to S$.
        Given a point $x \in X\bfk$ of $FX$ we necessarily have $f \cdot x = \alpha(\bfk, x)$, so the map is unique. Now we need to check that for such a defined $f$ we have $\alpha(H, x) \cdot L\xi = f \cdot (X\xi) \cdot x = \alpha(\bfk, (X\xi) \cdot x)$ for every $H$ and $x\in XH$, but this follows from the naturality of $\alpha$.

        Finally we have to show that such $f$ is continuous. Let $U\le S$ be an open subspace of $S$. By construction $f^{-1}(U)\le FX$ is open if and only $\eta(H,x)^{-1}f^{-1}(U)$ is open for every $H$ and  $x\in XH$ in $(\Yon \downarrow X)$.
        
        But by hypothesis for every $H$ and  $x\in XH$, $\alpha(H, x)$ is continuous and thus $\alpha(H,x)^{-1}(U)= \eta(H,x)^{-1} f^{-1} U$ is open in $LH$.
	\end{proof} 
\end{lemma}

\begin{lemma}\label{lem:counit_is_enrichment}
    Let $KS=F (\TVect(L, S))$ and $\varepsilon : K \to \id_{\TVect}$ be the counit of the adjunction $F \dashv \TVect(L,-)$. Then for every $S$ in $\TVect$, $\varepsilon S : KS \to S$ is mono and epi.
	\begin{proof}
		Recall that $\bfk$ is a separator for $\TVect$.
		So, given $f, g: S' \to KS$ with $f\neq g$, there is $h: \bfk \to S'$ with $f \cdot h \neq g \cdot h$. Now if $\varepsilon S \cdot f = \varepsilon S \cdot g$, then $\varepsilon S \cdot f \cdot h = \varepsilon S \cdot g \cdot h$, but since the domain of $h$ has the form $\bfk=L\bfk$, by the universal property of $\varepsilon$ it follows that $f \cdot h= g \cdot h$, contradiction.
		To see it is epi: let $f,g: S \to S'$ be s.t.\ $f\cdot \varepsilon S= g\cdot \varepsilon S$. Then for every $h: LH \to S$ we would have $f \cdot h = g \cdot h$, but again $\bfk=L\bfk$ is a separator, so it follows that $f=g$.
	\end{proof}
\end{lemma}

\begin{proposition}\label{prop:counit_description}
	For every $S$, $\varepsilon S : KS \to S$ is an enrichment of the topology on $S$, i.e.\ $KS$ has the same underlying vector space as $S$ but the topology is the finest topology for which all the continuous linear maps to $S$ from linearly compact spaces are still continuous in the new topology.
   \begin{proof}
        By Lemma~\ref{lem:counit_is_enrichment}, $\varepsilon S : KS \to S$ is an enrichment of the topology. The fact that the enriched topology is the one described in the statement is a consequence of the construction of $F$ (see the proof of Lemma~\ref{lem:FXexists}).
    \end{proof}
\end{proposition}
	
\begin{definition}
	A topological vector space $S$ is a \emph{K-space} if the map $\varepsilon S : F \TVect(L, S) \to S$ is an isomorphism, or, equivalently, if a subspace $U<S$ is open precisely when for every continuous linear $f: LH \to S$, $f^{-1}U$ is open in $LH$.
	Denote the full subcategory of $\TVect$ whose objects are the $K$-spaces, by $K\TVect$.
\end{definition}

\begin{lemma}\label{lem:F_factors}
    For $X \in [\Vect, \Vect]$, $FX$ is always a $K$-space, i.e.\ $\varepsilon F$ is a natural isomorphism.
    \begin{proof}
        Let $U\subseteq FX$ and suppose that $f^{-1}U$ is open for every continuous $f: LH \to FX$. Then in particular $f^{-1}U$ is open for every $f$ of the form $L\xi \mapsto X\xi \cdot x$ for some $x \in XH$. But then $f^{-1}U$ is open by definition of the topology on $FX$.
    \end{proof}
\end{lemma}
 
\begin{corollary}\label{cor:idempotent}
	The comonad $(K, \varepsilon, F \, \eta \, \TVect(L, -))$ on $\TVect$ is idempotent, in particular the K-spaces are objects for a full coreflective subcategory of $\TVect$.
    \begin{proof}
        Since $\varepsilon F$ is an isomorphism, by the unit-counit identities also $F \eta$ is an isomorphism, so $F \,\eta \, \TVect(L, -)$ is an isomorphism. The remainder of the statement follows form the relation between idempotent monads and reflective subcategories (see \cite[Cor.~4.2.4]{borceuxHCA2}).
    \end{proof}
\end{corollary}
		
\begin{corollary}
	$F$ factors through the inclusion $K\TVect \subseteq \TVect$, the factor 
	\[F': [\Vect, \Vect] \to K\TVect\]
	is a left adjoint to $K\TVect(L,-)$ and $K\TVect(L,-)$ is fully faithful.
	\begin{proof}
        The fact that $F$ factors through the inclusion $K\TVect \subseteq \TVect$ is due to Lemma~\ref{lem:F_factors}.
        The adjunction $F'\dashv K\TVect(L,-)$ is obtained restricting the adjunction $F\dashv \TVect(L,-)$.
        Finally $K\TVect(L,-)$ is fully faithful because the counit of the restricted adjunction is an isomorphism (see \cite[Ch.~IV, Sec.~3, Thm.~1]{maclane71categories}). To see this it suffices to observe that $\varepsilon K$ is an isomorphism because $\varepsilon F$ is (again by Lemma~\ref{lem:F_factors}).
	\end{proof}
\end{corollary}
	
\begin{remark}\label{rmk:fact_diagram}
	We thus gained a factorization of $F\dashv \TVect(L, -)$ as
	\[\begin{tikzcd}
		{[\Vect, \Vect]} \ar[rr, bend left = 10, "{F'}"{name=U1}] \ar[rrrr, bend left = 30, "{F=\int^H (-H) \cdot LH}"] 
		&&K\TVect \ar[ll, bend left = 10, "{K\TVect(L,-)}"{name=D1}] \ar[rr, bend left =10, hook, ""{name=U2}]
		&&\TVect \ar[ll, bend left = 10,"{K}"{name=D2}] \ar[llll, bend left = 30, "{\TVect(L,-)}"]
        \ar[from = U1, to = D1, phantom, "\bot" description]
        \ar[from = U2, to = D2, phantom, "\bot" description]
	\end{tikzcd}\]
	where $K\TVect(L,-)=K\TVect(L, -)$ is fully faithful and $F' \circ K\TVect(L,-)$ in an idempotent monad on $[\Vect, \Vect]$.
\end{remark}
	
\begin{proposition}\label{prop:KTVectchar}
	The category $K\TVect$ is the maximum dense extension of $L\Vect^\op$ in $\TVect$ containing all discrete subspaces, that is: if $\cC \subseteq \TVect$ is a full subcategory containing $L\Vect^\op$ and all discrete subspaces, and $L\Vect^\op \subseteq \cC$ is a dense extension, then $\cC \subseteq K\TVect$.
	\begin{proof}
		Assume $\cC \subseteq \TVect$ is a full subcategory containing $L\Vect^\op$ and all discrete subspaces and such that $L\Vect^\op\subseteq \cC$ is a dense extension, thus $\cC(L,-) = \TVect(L,-)|_\cC$ is fully faithful.
		We want to show that if $X \in \cC$, $\varepsilon X: KX \to X$ must be an isomorphism. By Lemma~\ref{lem:counit_is_enrichment}, we can assume $KX$ and $X$ are supported by the same vector space $V$ and we only need to prove that every $H$ which is open in the topology of $KX$, is open in the topology of $X$ as well.
		Since $H$ is open in $KX$, $KX/H$ is discrete, so given any $W \in \Vect$ and any $f \in \Vect(V/H, W)$, if $p\in \TVect(KX, (KX)/H)$ is the projection, $f \cdot p \in \TVect(KX,DW)$, where $DW$ is $W$ endowed with the discrete topology. Since by hypothesis $DW \in \cC$, it must be $\TVect(X,DW)= \TVect(KX,DW)$ so $f \cdot p \in \TVect(X,DW)$ for any $f \in \Vect(V/H, W)$. Choosing $W$ with $\dim W \ge \dim (V/H)$ we can pick $f$ to be injective, whence $\Kr(f\cdot p) = \Kr p = H$. Thus, since $DW$ is discrete and $f \cdot p: X \to DW$ is continuous, $H$ must be open in $X$.
	\end{proof}
\end{proposition}

\begin{remark}
	The category $K\TVect$ is not maximum among all dense extensions of $L\Vect^\op$. In other words Proposition~\ref{prop:KTVectchar} is not true without the assumption that $\cC$ contains all discrete spaces. A simple counterexample is given by $\cC=L\Vect^\op \cup \{X\}$ where $X$ is an infinite dimensional space whose open subspaces are precisely the finite-codimensional ones.
    The fact that $\TVect(L,-)$ is bijective on all $\TVect(LV, Y)$ is by Corollary~\ref{cor:idempotent} and by definition of $K$. Therefore one only needs to check that $\TVect(L,-)|_\cC$ is fully faithful on $\TVect(X, X)$ or equivalently that 
    $(-\cdot \varepsilon X): \TVect(X, X) \to \TVect(KX,X)$ is a bijection. But this follows from the fact that every linear endomorphism of $X$ is continuous and Proposition~\ref{prop:counit_description}.
    %
    %
\end{remark}

The following proposition provides a connection with the previous section.
	
\begin{proposition}\label{prop:TVect(L,S)_is_simple}
	For an object $X\in [\Vect, \Vect]$ among the statements
	\begin{enumerate}
		\item $X$ is in the essential image of $\TVect(L,-)$;
		\item $X$ is in the essential image of $K\TVect(L,-)$;
		\item $X$ is simple, i.e. $\cN_2 X=0$.
	\end{enumerate}
	the following implications hold: $1) \Leftrightarrow 2) \Rightarrow 3)$.
	\begin{proof}
		$1) \Leftrightarrow 2)$ is obvious at this stage.
		For $2) \Rightarrow 3)$ note that $\TVect(L, S)$ is $\bfk$-concrete for every $S$, because if $x \in \TVect(LH, S)$ is s.t. $\TVect(L\xi, S) \cdot x = 0$ for every $\xi \in H^*$, then $x=0$.
	\end{proof}
\end{proposition}
	
\begin{remark}
	The implication $3) \Rightarrow 2)$ does not hold. For a counterexample consider the left Kan-extension of some restriction of the double dualization functor to the category of vector spaces of dimension less than a given uncountable cardinal.
\end{remark}
	
Now we are going to relate separated spaces in $K\TVect$ with objects in $\Sigma\Vect$. First recall the definition of separated space.
	
\begin{definition}
	A (linearly topologized) topological vector space $S$ is \emph{separated} if the singleton $\{0\}$ is closed (this is equivalent to being T0 and to being T2). Denote the category of separated linearly topologized topological vector spaces by $\TVect_{s}$.
\end{definition}
	
In the remainder of this section, the main result will be Proposition~\ref{prop:KsepIffQ4free} from which Theorem~\ref{thm:KTVect_is_RCVS} follows. Toward Proposition~\ref{prop:KsepIffQ4free} the main step is to show that if $X$ is $\cQ_4$-free, then $FX$ is separated (Corollary~\ref{cor:Q4freeImpliesSeprated}).

\begin{remark}
    Recall that by Proposition~\ref{Q4-free_char}, if $X=\Vect(V,-)$, $Y\in \Ind\-(\Vect^\op)$ is $\cQ_4$-free, and $f: X \to Y$, then $\Kr(f)$ is representable, and thus also $\Img f \subseteq Y$ is representable.
\end{remark}

	
\begin{lemma}\label{lem:dense}
	If $X \in \Sigma\Vect$ and $(f_i)_{i \in \kappa}$ is a summable family in $X\bfk$, then the net of finite partial sums $(\sum_{i \in A} f_i)_{A\in [\kappa]^{<\omega}}$ converges to $\sum_{i \in \kappa} f_i$ with respect to the topology making the set $X\bfk$ into the space $FX$.
	\begin{proof}
		By definition $f_i = X\delta_i \cdot x$ for some $x \in X \bfk^{\oplus \kappa}$ and by definition $\sum_{i \in \kappa} f_i= X \sigma \cdot x$. Since $(\sum_{i \in A} L\delta_i^{\kappa})_{A \in [\kappa]^{<\omega}}$ converges to $L\sigma$ in $L\bfk^{\oplus I}$, and the topology must make the map $L\xi \mapsto X\xi \cdot x$ continuous, we have the thesis.
	\end{proof}
\end{lemma}
	
\begin{lemma}\label{lem:F(Q4-free_emb)_is_closed}
	Assume $f: X \to Y$ is mono in $\Ind\-(\Vect^\op)$ and assume $FY$ is separated and $f$ is a $\cQ_4$-free embedding. Then $Ff: FX \subseteq FY$ has a closed image.
	\begin{proof}
		A subspace $C$ of $Y\bfk$ is closed if and only if for every $V$ and $y \in YV$ the preimage of $C$ by $Y-\cdot y: \Vect(V, \bfk) \to Y\bfk$ is the annihilator of a subspace of $V$.
		Now, for $y \in YV$ consider the pullback
		\[\begin{tikzcd}
			X \ar[r] 		&Y\\
			Z \ar[r] \ar[u] &\Vect(V,-)\ar[u, "Y-\cdot y"]
		\end{tikzcd}\]
		We claim that $Z\cong\Vect(W,-)$ for some quotient $W=V/H$ for some $H<V$.
		Right orthogonal maps are closed under pullbacks, hence $Z \subseteq \Vect(V, -)$ should be a $\cQ_4$-free-embedding as well, hence it has a $\cQ_4$-free cokernel.
		But by Proposition~\ref{Q4-free_char} every map from a representable to a $\cQ_4$-free object has representable kernel, so $Z$ is representable.
	\end{proof}
\end{lemma}
	
\begin{corollary}\label{cor:Q4freeImpliesSeprated}
	If $X$ is $\cQ_4$-free then $FX$ is separated.
	\begin{proof}
		As a $\cQ_4$-free object, $X$ is in $\Ind\-(\Vect^{\op})$, so it is a quotient of a direct sum of representables. In particular, since it is $\cQ_4$-free, it is a cokernel of a $\cQ_4$-free embedding of some $Z$ into a direct sum of representables $Y= \bigoplus_{i \in I} \Yon V_i$.
		Now $F$, being a left adjoint, preserves colimits and $FX$ is a quotient of the separated space $FY$ by a closed subspace by Lemma~\ref{lem:F(Q4-free_emb)_is_closed}, hence it is separated.
	\end{proof}
\end{corollary}	
	
\begin{lemma}
	Let $V$ be a vector space and $\cF$ a filter of subspaces with $\bigcap \cF=0$ and let
	$Z=\bigcup_{H \in \cF} \Yon(V/H)\subseteq \Yon V$.
	Then the induced map $FZ \subseteq F \Yon V$ is a dense subspace of $F\Yon V = LV$.
	\begin{proof}
		The closed subspaces of $LV$ are precisely the $L(V/C)$ for some $C<V$.
		Note that the image of $F Z$ is precisely
		\[\bigcup_{H \in \cF} L(V/H) \subseteq LV\]
		so its points are the $\xi \in \TVect(L\bfk, LV)= V^*$ for which there is $H \in \cF$ with $\xi \cdot H=0$.\\
		Given any finite dimensional $A<V$, and $\xi \in V^*$, we must show that there is $H \in \cF$ with $(\xi + \Ann(A))\cap \Ann(H)\neq \emptyset$, or equivalently, that given any linear functional $\xi'$ on $A$ there is $\eta \in V^*$ with $\eta|A= \xi'$ and $\eta \cdot H=0$ for some $H \in \cF$.\\
		Now since $\bigcap \cF=0$ and $A$ is finite dimensional there has to be $H \in \cF$ such that $A \cap H=0$ and indeed any functional $\xi'$ on $A$ can be extended to one which is null on $H$.   
	\end{proof}
\end{lemma}
	
\begin{proposition}\label{prop:KsepIffQ4free}
	For an object $X\in [\Vect, \Vect]$ the following statements are equivalent:
	\begin{enumerate}
		\item $X$ is isomorphic to $\TVect(L,S)$ for some $S \in \TVect_{s}$;
		\item $X$ is isomorphic to $\TVect(L,S)$ for some $S \in K\TVect_{s}$;
		\item $X$ is $\cQ_4$-free.
	\end{enumerate}
    Moreover when they hold, the $X$-component $\eta X : X \to \TVect(L,FX)$ of the unit $\eta$ of the adjunction $F \dashv \TVect(L,-)$, is an isomorphism.
	\begin{proof}
		$(2) \Rightarrow (1)$ is obvious. $(1) \Rightarrow (2)$ follows from the fact that $KS$ is separated for a separated $S$ because $\varepsilon S : KS \to S$ is an enrichment.
  
		$(2) \Rightarrow (3)$. We know $X \cong \TVect(L, S)$ for a separated $K$-space $S$; let $V$ be a vector space, $\cF$ a filter of subspaces of $V$ with $\bigcap \cF=0$, $Z:= \bigcup_{H \in \cF} \Yon (V/H)$, and $\iota: Z\subseteq \Yon V$ the natural inclusion. Then by Lemma~\ref{lem:dense}, $F\iota$ has a dense image. Since $S$ is separated any $f: FZ \to S\cong\TVect(L,X)$ extends to at most one continuous $\bar{f}: LV \to S$ along the dense inclusion $F \iota$.
        
		Note that $\Yon V$ is $\bfk$-concrete because it is representable thus $Z\subseteq \Yon V$ is also simple. Moreover $X$ is $\bfk$-concrete by Proposition~\ref{prop:TVect(L,S)_is_simple}. Since $F$ is faithful when restricted to the category of $\bfk$-concrete objects, it follows that every map $h: Z \to X$ extends along $\iota$ to at most one map $\Yon V \to FX$. Thus $\Nat(\iota, X)$ is one-to-one and $X$ is $\cQ_4$-free.
  
		To prove $(3) \Rightarrow (2)$ and the ``moreover'', it suffices to show that if $X$ is $\cQ_4$-free, then $\eta X: X \cong \TVect(L,FX)$ is an isomorphism.
  
        Note that $\eta X$ is mono at every $\bfk$-concrete $X$, simply because by Lemma~\ref{lem:simple_prop}, $F$ is one-to-one on each $\Nat(Z,X)$ for every $Z$ (\cite[Ch.~IV, Sec.~3, Thm.~1]{maclane71categories}).

		Thus, since $[\Vect, \Vect]$ is Abelian we just need to show that $\eta X$ is epi. Recall that 
	    $(\eta X)H: XH \to \TVect(LH, FX)$ sends $x \in XH$ to the continuous function 
        \[LH \ni L\xi \mapsto X\xi \cdot x \in X\bfk.\]
		We have that if a function $f \in \TVect(LH, FX)\subseteq \Vect(\Vect(H, \bfk), X\bfk)$ is continuous, since $FX$ is separated, $f$ must preserve infinite sums (recall that $\cQ_4$-free objects are $\cQ_3$-free a fortiori). Therefore by Lemma~\ref{lem:sum-pres_characterization}, $f$ lifts to a map $\Vect(H, -) \to X$ i.e.\ via Yoneda, to the required element of $XH$.
	\end{proof}
\end{proposition}
	
\begin{remark}
	Let $V\in \Vect$ be infinite-dimensional and consider the quotient $X= \Vect(V,-)/(V \otimes_\bfk -)$ as in Example~\ref{ex:weirdo}.
	Then $FX$ is an indiscrete space and so $\TVect(L,FX)$ is not even small, in particular we have $\TVect(L, FX) \cong \Vect(\Vect(-,\bfk), X\bfk)$.
\end{remark}

\begin{theorem}\label{thm:KTVect_is_RCVS}
	The category $K\TVect_{s}$ is equivalent to the category of $\cQ_4$-free objects in $\Ind\-(\Vect^\op)$. Moreover $K\TVect_{s}$ is the largest dense extension of $L\Vect^\op$ in $\TVect_{s}$ containing all discrete subspaces; in particular $K\TVect_{s}$ is the largest category of reasonable strong vector spaces in $\TVect_{s}$ containing the discrete spaces.
	\begin{proof}
		As per Remark~\ref{rmk:fact_diagram}, $K\TVect(L,-): K\TVect \to [\Vect,\Vect]$, construes $K\TVect$ as a reflective subcategory of $[\Vect, \Vect]$. $K\TVect_{s}$ is in turn a reflective subcategory of $K\TVect$ and Proposition~\ref{prop:KsepIffQ4free} states precisely that the essential image of $K\TVect_{s}$ in $[\Vect,\Vect]$ consists of $\cQ_4$-free objects.
		The last assertion readily follows from Proposition~\ref{prop:KTVectchar}.
	\end{proof}
\end{theorem}
	
\begin{remark}\label{rmk:Based_are_topological}
    It is worth mentioning that all objects in $\Sigma\Vect$ of the form  $\bfk[\Gamma, \cF]$ are actually $\cQ_4$-free and so can be regarded as topological vector spaces.
\end{remark}

To conclude the comparison, we will show that all $\aleph_1$-ary strong vector spaces are of the topological kind. For this we will first need the following Proposition which is interesting in its own right.

\begin{proposition}\label{prop:countable_d_span_char}
	Let $\dim V = \aleph_0$ and $H < V^*$.
    \begin{enumerate}
        \item If $\Kr(H)=0$, then there is a basis $B$ of $V$ s.t.\ $\Span_\bfk \{\delta_b^B: b \in B\} \subseteq H$;
        \item If $\dim H \le \aleph_0$, then there is a basis $B$ of $V$ s.t.\ $\Span_\bfk \{\delta_b^B: b \in B\} \supseteq H$;
        \item If $\Kr(H)=0$ and $\dim H= \aleph_0$ then there is a basis $B$ of $V$ s.t. $\Span_\bfk \{\delta_b^B: b \in B\} = H$.
    \end{enumerate}
	\begin{proof}
        Note that (1) and (2) follow from (3). In fact, if $H$ is such that $\Kr(H)=0$, then $\dim H \ge \aleph_0$ and we can extract a subspace $H' \le H$ with $\dim(H')=\aleph_0$ and such that $\Kr(H')=0$. Since then the hypothesis of (3) holds for $H'$, if (3) holds, then we have that (1) holds for $H$. Similarly for (2): if $H < V^*$ is such that $\dim H \le \aleph_0$, then we can find $H'\ge H$ such that $\dim H'=\aleph_0$ and $\Kr(H')=0$, so by (3) there would be a basis $B$ of $V$ with $\Span_\bfk\{\delta_b^B: b \in B\}= H' \supseteq H$.

        Now let us prove (3). Consider a basis $\{\xi_n: n \in \omega\}$ of $H$ and a basis $\{e_n : n \in \omega\}$ of $V$. Define $V_n:= \bigcap_{i<n} \Kr (\xi_i)$ and $E_n:= \Span_\bfk \{e_i: i < n\}$. Since the $\xi_n$ are linearly independent we have $\dim V/V_n=n$, similarly $\dim(E_n)=n$. Consider the filter of subspaces $\cW:=\{W\le V: \exists n, V_n \subseteq W\}$ and let 
        \[\cF:=\{(U,W): U\le V,\; W \in \cW,\; V= U \oplus W\},\]
        be ordered by $(U,W) \le (U', W')$ if and only if $U \le U'$ and $W' \le W$. Now observe the following.
        \begin{itemize}
            \item For all $v \in V$ and $(U,W)\in \cF$, we can choose $(U,W) \le (U', W') \in \cF$ such that $v \in U'$: in fact by the hypothesis on $H$ we can find $n$ such that $v \notin V_n$ and $V_n \subseteq W$ and it suffices to pick $W'=V_n$ and choose $U'\ge U +  v \cdot \bfk$ so that $U'\oplus W'=V$;
            \item For all $\xi \in H$ and $(U,W)\in \cF$, we can choose $(U, W) \le (U', W') \in \cF$ such that $W' \subseteq \Kr(\xi)$: it suffices to choose $n$ so that $\xi \in \Span_\bfk\{\xi_i : i<n\}$ and such that $V_n \subseteq W$, then pick again $W':=V_n$ and let $U' \ge U$ be such that $V = U' \bigoplus W'$.
            \item Given $(U, W), (U', W') \in \cF$ with $(U, W) \le (U', W')$, there are $(U_i, W_i) \in \cF$ for $0 \le i\le n:=\dim(U'/U)$ such that $U_{i+1}/U_i \simeq \bfk$, $(U_0, W_0)=(U, W)$ and $(U_n, W_n)=(U', W')$. 
        \end{itemize}
        Thus it is possible to define an increasing sequence $\big((U_n, W_n): n \in \omega)$ in $\cF$ such that $U_{n+1}/U_n \simeq \bfk$ and for every $m \in \omega$, there is $N_m$ such that for all $n\ge N_m$, $E_m \le U_n$ and $V_m \ge W_n$.
        
        Choose $b_n \in U_{n+1}\cap (W_n \setminus W_{n+1})$ and note that then $B:=\{b_n: n \in \omega\}$ is a basis of $V$ because $e_m \in \Span_\bfk\{b_j :j<N_m\}$ for all $m \in \bN$, and moreover for all $n$, $W_{n} \cap \Span_\bfk\{b_j :j<n\}\subseteq W_n \cap U_{n} =0$, so $b_n \notin \Span_\bfk\{b_j :j<n\}$.
        
        Finally observe that by construction $\{W_n: n \in \omega\}$ generates $\cW$, i.e. $\cW=\{W: \exists n,\; W_n\le W\}$. It follows that $H = \bigcup_{n} \Ann(V_n)=\bigcup_{n} \Ann(W_n)= \Span\{\delta_{b_i}^N: i<n\}$.
    \end{proof}
\end{proposition}
	
\begin{corollary}
	If $X$ is $\aleph_1$-ary, then $X$ is $\cQ_4$-free if and only if it is $\cQ_3$-free.
    \begin{proof}
        Suppose $X$ is $\aleph_1$-ary and $\cQ_3$-free.
        Let $V$ be an infinite-dimensional vector space and let $x: \Yon V \to X$ be such that $K=\Kr(x) = \bigcup_{H \in \cF} \Yon (V/H)$ where $\cF$ is a filter of subspaces of $V$ such that $\bigcap \cF = \emptyset$. To show that $X$ is $\cQ_4$-free, we need to show that then $x=0$.
        
        Since $X$ is $\aleph_1$-ary, we can find $i: V'\le V$ with $\dim V'\le \aleph_0$ such that $x= x' \cdot \Yon i$. Then one has $K'=\Kr(x')=\bigcup_{H \in \cF'} \Yon (V'/H)$ where $\cF'=\cF|V':=\{H \cap V': H\in \cF\}$ and thus in particular $\bigcap \cF' = 0$. 
        
        Now $\bigcap \cF' = \emptyset$ entails that $T:=K'(\bfk) \subseteq (\Yon V)(\bfk)=V^*$ is such that $\Kr(T)=\{v \in V: \forall \xi \in T, \, \xi v=0\}=0$, whence by Proposition~\ref{prop:countable_d_span_char}, $T$ contains some $\{\delta_b^B: b \in B\}$ for some basis $B$ of $V'$.
        But then this means that $x'$ factors through the quotient $\Yon (V') / \bigoplus_{b \in B} \Yon(\delta_b^B)$ and thus since $X$ is $\cQ_3$-free, $x'=0$ and thus $x=x' \cdot \Yon i = 0$.
    \end{proof}
\end{corollary}
	
\subsection{Some questions}
	
In general the functor $\bfk[-;-]$ described in Subsection~\ref{ssec:Comp_B} and its essential image (equivalent to $B\Sigma\Vect$) have a mysterious nature in the following regard: it seems ``hard'' to determine whether an $X \in \Sigma\Vect$ has the form $\bfk[-;-]$ or not. A necessary condition would be that it is $\cQ_4$-free, but it seems unlikely it is sufficient.
		
\begin{question}
	Are there non-based $\cQ_4$-free strong vector spaces (i.e. that are topological in nature)?
\end{question}
	
\begin{question}
	Does the functor $\bfk[-;-]$ reflect isomorphisms?
\end{question}
	
\section{Monoidal Closed Structure}\label{sec:ClMon}

Toward the goal of considering ``strongly linear $\bfk$-algebras'' we will now be concerned with defining a suitable monoidal closed structure on reasonable categories of strong vector spaces and study some of their basic properties.

The plan will be the following:
\begin{itemize}
    \item In Subsection~\ref{ssec:Tensor_Products} we will make some basic observations on the natural monoidal structure on $\Ind(\Vect^\op)$ given by the indization $-\hotimes-$ of the natural tensor product on $\Vect^\op$. In particular we will show that $B\Sigma\Vect$ and the category of $\cQ_4$ objects are closed under $-\hotimes$ (Lemma~\ref{lem:hotimes_on_BSigmaVect} and Proposition~\ref{prop:hotimes_on_Q4_free}) and that $\Sigma\Vect$ is not (Example~\ref{ex:hotimes_on_SigmaVect});
    \item Afterwards (Subsection~\ref{ssec:Internal_Hom}) we will show that $-\hotimes-$ makes $\Ind\-(\Vect^\op)$ into a closed monoidal category, i.e. that each $Y \hotimes -$ has a right adjoint $\Hom(Y,-)$ (Proposition~\ref{prop:hotimes_is_closed}). We will also show that the r.c.s.v.s.\ considered so far, $B\Sigma\Vect$, $K\TVect_{s}$ and $\Sigma\Vect$ are all closed under $\Hom(-,-)$ (Propositions~\ref{prop:Hom_on_BSigma} and \ref{prop:Hom_restricts}). We will conclude the subsection with some (mainly notational) remarks on the monoidal closed structure on $\Ind\-(\Vect^\op)$. 
    \item Finally in the Subsections~\ref{ssec:SLAlg} and \ref{ssec:SLDifferential} we will consider strongly linear algebras and their modules, observing that these admit a notion of strongly linear Kähler differential.  
\end{itemize}

\subsection{Tensor products}\label{ssec:Tensor_Products}
		
The category $\Ind\-(\Vect^{\op})$ has a natural monoidal structure induced by the one on $\Vect^{\op}$ given by $-\otimes_\bfk -$. Indeed composing the natural equivalence
\[\Ind\-(\Vect^{\op} \times \Vect^\op) \cong \Ind\-(\Vect^{\op}) \times \Ind\-(\Vect^{\op})\]
with $\Ind(-\otimes_\bfk-): \Ind\-(\Vect^{\op} \times \Vect^\op) \to \Ind\-(\Vect^{\op})$ (see \cite[Prop.~6.1.9]{kashiwara2005categories}) gives a functor
\[-\hotimes - : \Ind\-(\Vect^{\op}) \times \Ind\-(\Vect^{\op}) \to \Ind\-(\Vect^{\op})\]
given by
\[\begin{aligned}
	(X \hotimes  Y) V \cong \varinjlim_{x \in (\Yon \downarrow X)} \varinjlim_{y \in (\Yon \downarrow Y)}
	\Vect(Px \otimes Py, V) \cong\\
	\cong \int^{H_0}\int^{H_1 } X(H_0) \otimes Y(H_1) \otimes \Vect(H_0 \otimes H_1, V),
\end{aligned}\]
where we denote with the same letter $P$, the canonical $\Yon$-diagrams $P: (\Yon \downarrow X) \to \Vect^\op$ and $P: (\Yon \downarrow Y) \to \Vect^\op$ of $X$ and $Y$ respectively.

\begin{remark}[Elements of $X \hotimes Y$]
	Recall that by Yoneda if $x \in XV$ then $\dot{x}=(X-) \cdot x : \Yon V \to X$ uniquely represents $x$. Now if $y \in YW$, we have a natural transformation $\dot{x} \hotimes \dot{y} : \Yon V \hotimes \Yon W \to X \hotimes Y$.
		
	Recalling that $\Yon (V \hotimes W) \cong \Yon V \hotimes \Yon W$, there has to be a unique element $x \hotimes y \in (X\hotimes Y)(V \otimes W)$ such that $((X\hotimes Y) -)\cdot (x \hotimes y)$ is the composite
	\[\begin{tikzcd}
        \Yon (V \hotimes W) \cong \Yon V \hotimes \Yon W \ar[r, "{\dot{x}\hotimes \dot{y}}"] &X \hotimes Y.
    \end{tikzcd}\]
	Slightly more concretely $x \hotimes y$ is the image of
	\[x \otimes y \otimes \id_{V\otimes W} \in XV \otimes YW \otimes \Vect(V \otimes W, V \otimes W) \]
	for the $(V,W)$-component of the universal wedge to the coend 
    \[X \hotimes Y\cong \int^{V}\int^{W} X(V) \otimes Y(W) \otimes \Vect(V \otimes W, -).\]
\end{remark}
	
\begin{remark}[$-\hotimes-$ represents strong bilinear maps in $\Sigma\Vect$]\label{rmk:hotimes_and_strongbilinear}
	Given $X, Y \in \Ind\-(\Vect^{\op})$, denote by $X \times Y$ the (non-additive) functor
	\[X \times Y : \Vect \times \Vect \to \Vect, \quad (V, W) \mapsto XV\times XW\]
	then there is a natural transformation 
	\[\eta: X \times Y \to (X \hotimes  Y)(- \otimes -)\]
	such that $\Nat(X \times Y, Z (-\otimes-)) \cong \Nat( (X \hotimes  Y) (- \otimes -), Z (- \otimes -))$ via $\Nat(\eta, -)$.
	More precisely for every $V, W \in \Vect$, $x \in XV$ and $y \in YW$, $\eta$ sends $(x,y)$ to 
	\[x \hotimes y \in (X\hotimes Y) (V \hotimes W).\]
	In particular if $X, Y, Z\in \Sigma\Vect$ the isomorphism $\Nat(\eta, -)$ gives a natural bijection between maps $f:X \times Y \to Z$ such that
    \[(f(\bfk, \bfk))\left(\sum_{i\in I} x_i, \sum _{j \in J} y_j\right) = \sum_{(i,j) \in I \times J} (f\bfk)(x_i, y_j)\]
	for every $I$, $J$, $(x_i)_{i \in I} \in (X\bfk)^I$, $(y_j)_{j \in J} \in (Y\bfk)^J$, and maps $X \hotimes  Y \to Z$, or equivalently maps $\cR^\infty_3(X \hotimes Y)\to Z$ (recall Notation~\ref{not:reflections}).
\end{remark}
 
\begin{remark}
	The functor $- \hotimes  -$ is exact in both arguments.
	In fact every exact sequence $\begin{tikzcd}[cramped]0 \ar[r] & Y_0 \ar[r] &Y_1 \ar[r] &Y_2 \ar[r] &0\end{tikzcd}$ in $\Ind\-(\Vect^\op)$ is isomorphic to a sequence of the form
	\[\begin{tikzcd}
	    0 \ar[r] & \displaystyle{\varinjlim_{d \in D} \Yon F_0d} \ar[r] &\displaystyle{\varinjlim_{d \in D} \Yon F_1 d} \ar[r] &\displaystyle{\varinjlim_{d \in D} \Yon F_2 d} \ar[r] &0.
	\end{tikzcd}\]
	Here $D$ is a directed poset, $F_0, F_1, F_2 : D \to \Vect^\op$ are diagrams, and the arrows are of the form $\varinjlim \alpha$ and $\varinjlim \beta$ for $\alpha : F_0 \to F_1$ and $\beta : F_1 \to F_2$ such that for every $d \in D$ the short sequence given by $\alpha d$ and $\beta d$ is exact (\cite[Prop.~8.6.6~(a)]{kashiwara2005categories}).
	If $X \cong \varinjlim_{b \in B} \Yon Gb$ for a directed diagram $G: B \to \Vect^\op$, then applying $X \hotimes  -$ yields
	\[\begin{tikzcd}[column sep = 1.2em]
        0 \ar[r] & \displaystyle{\varinjlim_{\substack{d \in D\\b \in B}} \Yon (F_0d \otimes Gb)} \ar[r] &\displaystyle{\varinjlim_{\substack{d \in D\\b \in B}} \Yon (F_1d \otimes Gb)} \ar[r] &\displaystyle{\varinjlim_{\substack{d \in D\\b \in B}} \Yon (F_2d \otimes Gb)} \ar[r] &0
    \end{tikzcd}\]
	which is exact because each 
	\[\begin{tikzcd}
		0 \ar[r]
		&F_0d \otimes Gb \ar[rr, "{\alpha d \otimes Gb}"]
		&&F_1d \otimes Gb \ar[rr, "{\beta d \otimes Gb}"]
		&&F_2d \otimes Gb \ar[r]
		&0	 		
	\end{tikzcd}\]
	is exact.
\end{remark}

\begin{lemma}\label{lem:hotimes_on_BSigmaVect}
	For $\Gamma, \Delta$ sets and $\cF, \cG$ bornologies on $\Gamma$ and $\Delta$ respectively, we have a canonical isomorphism
	\[\bfk[\Gamma; \cF] \hotimes  \bfk[\Delta; \cG] \cong \bfk[\Gamma \times \Delta; \cF \otimes \cG]\]
	where $\cF \otimes \cG$ is the smallest bornology on $\Gamma\times \Delta$ containing all the $F \times G$ with $F \in \cF$ and $G \in \cG$.
	\begin{proof}
		The canonical isomorphism is given by the following chain of natural isomorphisms
		\[\bfk[\Gamma; \cF] \hotimes  \bfk[\Delta; \cG] \cong \varinjlim_{\substack{F \in \cF\\G \in \cG}} \Vect\left(\frac{\bfk^{\oplus\Gamma}}{\bfk^{\oplus (\Gamma\setminus F)}} \otimes \frac{\bfk^{\oplus\Delta}}{\bfk^{\oplus (\Delta \setminus G)}}, -\right)\cong\]
		\[\cong \varinjlim_{\substack{F \in \cF\\G \in \cG}} \Vect\left(\frac{\bfk^{\oplus(\Gamma\times \Delta)}}{\bfk^{\oplus (\Gamma \times \Delta \setminus F \times G)}}, -\right)\cong \bfk[\Gamma \times \Delta; \cF \otimes \cG].\]
	\end{proof}
\end{lemma}
	
\begin{lemma}\label{lem:toname2}
	If $V \in \Vect$, then for every $X$ in $\Ind\-(\Vect^\op)$
	\[X \hotimes  (\Yon V) \cong X \circ (\Yon V) \quad \text{and} \quad X \hotimes  (\Fin V) \cong (\Fin V) \circ X.\]
	In particular $\Fin \bfk \cong \Yon \bfk$ is a unit for the monoidal structure on $\Ind\-(\Vect^\op)$ with product $- \hotimes  -$.
	\begin{proof}
		For the first assertion observe that
		\[X \hotimes \Yon V \cong \varinjlim_{x : \Yon W \to X} \Yon (W \otimes V) \cong \varinjlim_{x : \Yon W \to X} \Yon W \circ \Yon V\cong X \circ (\Yon V)\]
		For the second observe that 
		\[\Fin V \hotimes X \cong \varinjlim_{\substack{H \le V\\\dim H< \aleph_0}} (\Fin H) \hotimes X \cong \varinjlim_{\substack{H \le V\\\dim H< \aleph_0}}(\Fin H) \circ X \cong (\Fin V) \circ X.\]
		Here the second isomorphism is because when $H$ is finite dimensional, $(\Fin H) \cong (\Yon H^*)$ and $(\Yon H^*) \circ X \cong X \circ (\Yon H^*)$. This in turn is because $X: \Vect \to \Vect$ preserves biproducts and $H^* \simeq \bfk^{\oplus n}$ for some $n$ so $\Yon (H^*) \cong (-)^{\oplus n}\cong (-)^n$.
	\end{proof}
\end{lemma}
	
\begin{proposition}\label{prop:hotimes_on_Q4_free}
	If $X_0, X_1 \in \Ind-(\Vect^\op)$ are $\cQ_4$-free, then $X_0 \hotimes X_1$ is $\cQ_4$-free.
	\begin{proof}
		By Proposition~\ref{Q4-free_char}, we can write $X_0$ and $X_1$ as directed unions of representables:
		\[X_i\cong \varinjlim_{d \in D_i} \Yon F_id\]
		with $F_i(d<d')$ an epimorphism in $\Vect$ for every $d,d' \in D_i$, $d<d'$. Now the tensor product is
		\[X_0\hotimes X_1\cong \varinjlim_{\substack{d_0 \in D_0\\d_1\in D_1}} \Yon (F_0 d_0 \otimes F_1 d_1).\]
		So, it suffices to observe that if $(d_0, d_1)<(d_0', d_1')$ then 
		\[F_0(d_0<d_0') \otimes F_1(d_1<d_1') : F_0(d_0) \otimes F_1(d_1) \to F_0(d_0') \otimes F_1(d_1')\]
		is an epimorphism in $\Vect$ and conclude again by Proposition~\ref{Q4-free_char}.
	\end{proof}
\end{proposition}
	
\begin{remark}
	Contrary to the case of $\cQ_4$-free objects, it is not true that if $X_0, X_1$ are $\cQ_3$-free then $X_0 \hotimes X_1$ is necessarily $\cQ_3$-free, as shown by the example below. 
\end{remark}
	
\begin{example}\label{ex:hotimes_on_SigmaVect}
	Recall Example~\ref{ex:weirdo}: if $V\in \Vect$ is infinite dimensional, then the cokernel $C$ of the canonical embedding $\iota :V \otimes - \subseteq \Vect(V^*, -)$ is $\cQ_3$-free but not $\cQ_4$-free.
		
	We are going to show that $C \hotimes \Yon V$ is not $\cQ_3$-free. By exactness of $- \hotimes -$, this amounts to show that $\iota \hotimes \Yon V$ is not a $\Sigma$-closed embedding. Observe that $(\iota \hotimes \Yon V)(W)$ is given, for $W \in \Vect$, by
	\begin{equation}\label{eq:iota}
		\begin{aligned}
			V \otimes \Vect(V, W) &\to \Vect(V^* \otimes V, W)\\
			v\otimes h &\mapsto \big(\xi \otimes v' \mapsto (\xi \cdot v) \cdot (h\cdot v')\big)
		\end{aligned}
	\end{equation}
	Showing that $i \hotimes \Yon V$ is not a $\Sigma$-closed embedding amounts to showing that for some $\lambda \in \Card$ the following square diagram is not a pullback
	\[\begin{tikzcd}[row sep = 0.05cm, column sep = 0.00cm]
		{}&v \otimes h \ar[r, mapsto] &\big(\xi \otimes v' \mapsto  (h \cdot v') \cdot (\xi \cdot v) \big)\\
		v \otimes h \ar[ddd, mapsto]&V \otimes \Vect(V, \bfk^{\oplus \lambda}) \ar[r] \ar[ddd] 	&\Vect(V^* \otimes V, \bfk^{\oplus \lambda})\ar[ddd] & f \ar[ddd, mapsto]\\ &&&{}\\  &&&{}\\
		\big(v \otimes (\delta_i \cdot h)\big)_{i \in \lambda} &\big(V \otimes \Vect(V, \bfk)\big)^\lambda \ar[r]			&\Vect(V^* \otimes V, \bfk)^\lambda&  (\delta_i \cdot f)_{i \in \lambda}\\
		&\bigg(v_i \otimes \eta_i \ar[r, mapsto] &\big(\xi \otimes v' \mapsto (\eta_i \cdot v')\cdot (\xi \cdot v_i)\big) \bigg)_{i \in \lambda}
    \end{tikzcd}\]
	Let $\lambda= \dim V$. Without loss of generality assume $V=\bfk^{\oplus \lambda}$. Now observe that $f \in \Vect\big((\bfk^{\oplus \lambda})^* \otimes \bfk^{\oplus \lambda}, \bfk^{\oplus \lambda}\big)$ given by $f( \xi \otimes e_i ) = e_i \cdot (\xi \cdot e_i)$ and $(e_i \otimes \delta_i)_{i \in \lambda} \in \big(\bfk^{\oplus \lambda} \otimes \Vect(\bfk^{\oplus \lambda}, \bfk)\big)^\lambda$ have the same image in $\Vect\big((\bfk^{\oplus \lambda})^* \otimes \bfk^{\oplus \lambda}, \bfk \big)^\lambda$ because $\delta_i \cdot f (\xi \otimes e_i) = \xi \cdot e_i = (\delta_i \cdot e_i) \cdot (\xi \cdot e_i)$. On the other hand, the $\lambda$-uple $(e_i \otimes \delta_i)_{i \in \lambda}$ is not in the image of the vertical arrow, so the square diagram is not a pullback.
			
	From a more concrete perspective, we showed that if we set 
	\[Y:=\Vect((\bfk^{\oplus \lambda})^* \otimes \bfk^{\oplus \lambda}, -), \qquad X:=\bfk ^{\oplus \lambda} \otimes \Vect(\bfk^{\oplus \lambda}, -),\]
	and let $\iota : X \to Y$ be given as in \ref{eq:iota} with $V=\bfk^{\oplus \lambda}$, then the family $(x_i)_{i \in \lambda} \in \Vect((\bfk^{\oplus \lambda})^* \otimes \bfk^{\oplus \lambda}, \bfk)$ defined by $x_i(\xi \otimes e_i)= \xi \cdot e_i$ is summable in $Y$ and its sums is given by
	\[(\sum_{i \in \lambda} x_i \cdot k_i)(\xi \otimes e_i)= (\xi \cdot e_i) \cdot k_i.\]
	Moreover for every $i \in \lambda$, we have $x_i \in (\Img \iota)(\bfk)$, however $\sum_{i \in \lambda} x_i \notin (\Img \iota)(\bfk)$.
\end{example}

\subsection{Internal Hom}\label{ssec:Internal_Hom} We now deal with the internal Hom.
	
\begin{proposition}\label{prop:hotimes_is_closed}
	For every $X\in \Ind\-(\Vect^{\op})$, the functor 
	\[X \hotimes  - : \Ind\-(\Vect^{\op}) \to \Ind\-(\Vect^{\op})\]
	has a right adjoint $\Hom(X, -)$ computed as
	\[\begin{aligned}
		\Hom(X,Y)(V)&\cong
		\Ind\-(\Vect^{\op}) (X \hotimes  \Yon V, Y) \cong\\
		&\cong \varprojlim_{x \in (\Yon \downarrow X)^{\op}} \varinjlim_{y \in (\Yon \downarrow Y)} \Vect(Py, Px \otimes V) \cong\\
		&\cong \int_{H} \Vect\big(XH, Y(H \otimes V)\big) \cong \\
		&\cong \int_{H_0}\int^{H_1} \Vect\big(XH_0, YH_1 \otimes \Vect(H_1 , H_0 \otimes V) \big).
	\end{aligned}\]
		
	\begin{proof}
		We need to check that such a defined $\Hom(X,Y)$ is small.
		Since $\Vect^{\op}$ has small limits, so do $\Ind\-(\Vect^{\op})$ and $[\Vect, \Vect]$. Moreover the inclusion $\Ind\-(\Vect^{\op}) \subseteq [\Vect, \Vect]$ preserves them, so they are computed pointwise: it means that the question whether $\Hom(X,Y)\in \Ind\-(\Vect^{\op})$ reduces to finding out whether for every $V, W \in \Vect$ one has $\Vect(V, W \otimes -) \in \Ind\-(\Vect^{\op})$. We see this is the case: indeed $\Vect(\bfk^{\oplus S}, W \otimes -) \cong (W \otimes -)^S$ and $\Ind\-(\Vect^{\op})$ has small products.
	\end{proof}
\end{proposition}

\begin{remark}\label{rmk:Hom_is_continuous}
    Note that $\Hom: \Ind\-(\Vect^\op)^\op \times \Ind\-(\Vect^\op) \to \Ind\-(\Vect^\op)$ is continuous. It is continuous in the second argument because it is a right adjoint, and it is continuous in the first argument by construction (because it is explicitly given as a limit in the first argument).
\end{remark}
		
\begin{proposition}\label{prop:Hom_on_BSigma}
	For $\Gamma, \Delta$ sets and $\cF, \cG$ bornologies on $\Gamma$ and $\Delta$ respectively, we have
	\[\Hom\big(\bfk[\Gamma, \cF]; \bfk[\Delta, \cG]\big) \cong \bfk[\Gamma \times \Delta; (\cF \to \cG)].\]
	Here $\cF \to \cG$ is the bornology on $\Gamma \times \Delta$ consisting of those $S$ such that for every $F \in \cF$,  there is a $G \in \cG$ such that $S \cap (F \times \Delta)\subseteq F \times G$ and $S \cap (F \times \Delta)$ has finite fibers over $G$. 
	\begin{proof}
        For $\Gamma' \subseteq \Gamma$, let $\cF\|\Gamma'$ be the bornology on $\Gamma$, given by
        \[\cF\|\Gamma':=\{(A\cap \Gamma') \cup B: A \in \cF, |B|<\aleph_0\}.\]
        Note that for each $F\subseteq F' \in \cF$, naturally $\bfk[\Gamma, \cF\|F] \subseteq \bfk[\Gamma, \cF\|F']\subseteq \bfk[\Gamma, \cF]$ and that 
        \[\bfk[\Gamma, \cF] = \varinjlim_{F \in \cF} \bfk[\Gamma, \cF\|F] = \bigcup_{F \in \cF} \bfk[\Gamma, \cF\|F].\]
        It is not hard to check that
        \[\Hom(\bfk[\Gamma, \cF\|F], Y)(V) \cong Y(V^{\oplus F}) \oplus (Y(V))^{\Gamma \setminus F} \cong \bfk[\Gamma \times \Delta, \cH_F](V)\]
        where $S \in \cH_F$ if and only if for all $\gamma \in \Gamma$, $\{\delta: (\gamma, \delta) \in S\}\in \cG$ and
        \begin{itemize}
            \item there is $G \in \cG$ with $S \cap (F \times \Delta)\subseteq F \times G$ and $S \cap (F \times \Delta)$ has finite fibers over $G$.
        \end{itemize}
        Note that if $F \subseteq F'$, then the natural inclusion $\bfk[\Gamma \times \Delta, \cH_{F'}]\subseteq \bfk[\Gamma \times \Delta, \cH_{F}]$ is identified through the isomorphisms $\Hom(\bfk[\Gamma, \cF\|F], Y) \cong \bfk[\Gamma \times \Delta, \cH_F]$ with the image under the functor $\Hom(-, Y)$ of the inclusion $\bfk[\Gamma, \cF\|F] \subseteq \bfk[\Gamma,\cF\|F']$.

        Thus, since $\Hom(-,Y)$ is continuous (Remark~\ref{rmk:Hom_is_continuous}), we get 
        \[\Hom(X, Y)\cong \varprojlim_{F \in \cF} \Hom(\bfk[\Gamma, \cF\|F], Y) \cong \bigcap_{F \in \cF} \bfk[\Gamma\times \Delta, \cH_F]\]
        which yields the bornology in the statement.
	\end{proof}
\end{proposition}
	
Now we show the classes of $\cQ_4$-free objects and of $\cQ_3$-free objects (i.e. objects in $\Sigma\Vect$) are closed under $\Hom(-, -)$.
	
\begin{proposition}\label{prop:Hom_restricts}
	Let $\cQ=\cQ_i=\cQ(\bE_i)$ for $i \in \{3,4\}$ as in Example~\ref{ex:mainex}. Assume $Y\in \Ind\-(\Vect^\op)$ is $\cQ$-free. Then $\Hom(X, Y)$ is $\cQ$-free for every $X \in \Ind\-(\Vect^\op)$.
	\begin{proof}
		Let $Q \in \cQ$ and observe that by definition of $\Hom$,
		\[\Nat(Q, \Hom (X,Y))\cong \Nat(Q \hotimes X , Y).\]
		Thus it suffices to showing that $Q\hotimes X$ is a torsion object for the torsion theory $({}^\bot(\cQ^\bot),\cQ^\bot)$, that is $Q \hotimes X \in \Tors({}^\bot(\cQ^\bot))$. Since the torsion objects are closed under colimits and $-\hotimes -$ preserves colimits, this reduces to show that $Q\hotimes \Yon V\in \Tors({}^\bot(\cQ^\bot))$ for every representable $\Yon V$.
		Now recall that by Lemma~\ref{lem:toname2}, $Q \hotimes (\Yon (\bfk^{\oplus \lambda})) \cong Q \circ (\Yon (\bfk^{\oplus \lambda})) \cong Q(-^\lambda)$.
			
		For $i=3$, for every $\mu \in \Card$ we have
		\[\Vect(\bfk,-)^{\oplus (\lambda \times \mu)} \hookrightarrow (\Vect(\bfk, -)^\lambda)^{\oplus \mu} \hookrightarrow \Vect(\bfk, -)^{\lambda \times \mu},\]
		thus $Q(-^\lambda)=\Ckr\big( (\Vect(\bfk, -)^\lambda)^{\oplus \mu} \hookrightarrow \Vect(\bfk, -)^{\lambda \times \mu}\big)$ is a quotient of the cokernel $\Ckr \big(\Vect(\bfk,-)^{\oplus (\lambda \times \mu)} \hookrightarrow \Vect(\bfk, -)^{\lambda \times \mu}\big)$ and thus is in $\Tors({}^\bot(\cQ^\bot))$.
		For $i=4$ it is enough to observe that if $(W, \cF) \in \bE_4$ then the filter $\cF_\lambda$ of subspaces of $W^\lambda$ generated by $\{H^\lambda: H \in \cF\}$ is still such that $\bigcap\cF_\lambda=0$ so $(W^\lambda, \cF_\lambda)$ is still in $\bE_4$.
	\end{proof}
\end{proposition}

\begin{remark}[Summable families in $\Hom(X, Y)$]
    Let $X$ and $Y$ in $\Sigma\Vect$, $(f_i)_{i \in I}$ be an $I$-indexed family in $\Hom(X, Y)(\bfk)\cong \Sigma\Vect(X,Y)$. Then $(f_i)_{i \in I}$ is summable if and only if for every set $J$ and every $J$-indexed summable family in $X\bfk$, the family $(f_i(x_j))_{I \times J}$ is a summable family in $Y\bfk$. To see this note that by the $\Hom$-$\hotimes$ adjunction and by Remark~\ref{rmk:hotimes_and_strongbilinear} an $I$-indexed family is in the image of $\Hom(X, Y)(\delta_I):\Hom(X,Y)(\Yon (\bfk^{\oplus I})) \to \Hom(X,Y)(\bfk)^I$ if and only if it induces a strongly bilinear map $\bfk^I \times X\bfk \to Y\bfk$.
\end{remark}
	
\begin{corollary}\label{cor:Induced_tensor}
	Let $\cQ=\cQ_i=\cQ(\bE_i)$ for $i \in \{3,4\}$ as in Example~\ref{ex:mainex}.
	The monoidal closed structure given by $- \hotimes -$ and $\Hom$ on $\Ind\-(\Vect^\op)$ induces a monoidal closed structure on the reflective subcategory of $\cQ$-free objects given by $\cR^\infty_\cQ(- \hotimes -)$ and $\Hom$.
	\begin{proof}
		Set $\cR^\infty:=\cR^\infty_\cQ$ (recall Notation~\ref{not:reflections}). Observe that for every $\cQ$-free $X$, $\cR^\infty_\cQ(-\hotimes X)$ is the left adjoint of the restriction of $\Hom(X,-)$ to the category of $\cQ$-free objects which is well defined by Proposition~\ref{prop:Hom_restricts}.
			
		The only non-trivial thing left to prove is that $\cR^\infty_\cQ(- \hotimes -)$ is associative (i.e.\ it has a natural associator).
		Observe that there are bijections, natural in $X_0, X_1, X_2, Z$ ranging over the $\cQ$-free objects
		\[\begin{aligned}
			\Nat(\cR^\infty(\cR^\infty(X_0 \hotimes X_1) \hotimes X_2), Z) \cong \Nat(\cR^\infty(X_0 \hotimes X_1) \hotimes X_2, Z) \cong\\
			\cong \Nat(\cR^\infty(X_0 \hotimes X_1), \Hom(X_2,Z))\cong \Nat(X_0\hotimes X_1, \Hom (X_2, Z))\cong \\
			\cong \Nat((X_0 \hotimes X_1) \hotimes X_2, Z)  \cong \Nat(\cR^\infty((X_0 \hotimes X_1) \hotimes X_2), Z)
		\end{aligned}\]
		and that their composition is $\Nat(\cR^\infty(\rho^\infty(X_0 \hotimes X_1) \hotimes X_2), Z)$. This, since $Z$ ranges over all $\cQ$-free objects, implies that
		\[\cR^\infty(\rho^\infty(X_0 \hotimes X_1) \hotimes X_2): \cR^\infty((X_0 \hotimes X_1) \hotimes X_2) \to \cR^\infty(\cR^\infty(X_0 \hotimes X_1) \hotimes X_2)\]
		is an isomorphism. Similarly
		\[\cR^\infty(X_0 \hotimes \rho^\infty(X_1 \hotimes X_2)): \cR^\infty(X_0 \hotimes X_1 \hotimes X_2) \to \cR^\infty(X_0 \hotimes \cR^\infty(X_1 \hotimes X_2))\]
		is an isomorphism as well. The associator can then be obtained from the associator of $-\hotimes-$ from the composition of isomorphisms:
		\[\begin{aligned}
			\cR^\infty(\cR^\infty(X_0 \hotimes X_1) \hotimes X_2)\cong \cR^\infty((X_0 \hotimes X_1) \hotimes X_2) \cong\\
			\cR^\infty(X_0 \hotimes (X_1 \hotimes X_2))\cong \cR^\infty(X_0 \hotimes \cR^\infty(X_1 \hotimes X_2)),
		\end{aligned}\]
		and its required properties can be checked likewise from the analogous ones for the associator of $-\hotimes-$.
	\end{proof}
\end{corollary}

\begin{corollary}
	The monoidal closed structure given by $- \hotimes -$ and $\Hom$ on $\Ind\-(\Vect^\op)$ induces a monoidal closed structure on $\Sigma\Vect$ given by $\cR^\infty_3(- \hotimes -)$ and $\Hom$.
	\begin{proof}
		This is the previous Corollary with $\cQ=\cQ_3$.
	\end{proof}
\end{corollary}

\begin{remark}
	The category $\Ind\-(\Vect^\op)$ with the above defined $\Hom(-,-)$ and $-\hotimes-$ is a commutative monoidal closed category and we saw it restricts to some subcategories of $\cQ$-free objects. For the remainder of this Chapter we fix the following symbols: 
	\begin{itemize}
		\item $\alpha: (X \hotimes Y) \hotimes Z \to X \hotimes (Y\hotimes Z)$ for its associator;
		\item $\nu^l : (\Yon\bfk) \hotimes X \to X$ and $\nu^r: X\hotimes (\Yon \bfk) \to X$ for its left and right unitors;
		\item $\tau: X \hotimes Y \to Y \hotimes X$ for its commutor;
		\item we denote the $\otimes\dashv \Hom$ adjunction by
		\[-^\sharp: \Nat(X \hotimes Y, Z)\cong \Nat(X, \Hom(Y, Z)) : -^{\flat}\]
		and its unit and counit by $\ad(X) : X \to \Hom(Y, X\hotimes Y)$ and $\ev(Z) : \Hom(Y,Z) \hotimes Y \to Z$ respectively.
	\end{itemize}  
\end{remark}

\begin{remark}
	As for any monoidal closed category, the natural isomorphism $\Nat(X \hotimes Y, Z)\cong \Nat(X, \Hom(Y, Z))$, lifts to an internalized version \[\Hom(X \hotimes Y, Z)\cong \Hom(X, \Hom(Y, Z)),\] for which we use the same symbols. We will also make use of the internalized version of the components of the functors $\Hom(-,-)$ and $-\hotimes-$. Given $X, Y$, there are arrows in $\Ind\-(\Vect^\op)$, natural in $X$, $Y$, and $S$,
	\[S \hotimes \bullet : \Hom (X, Y) \to \Hom (S \hotimes X, S \hotimes Y),\]
	given by the composition
	\[\Hom (X, Y) \longrightarrow \Hom (X, \Hom (S, S\hotimes Y)) \cong \Hom (S \hotimes X, S \hotimes Y),\]
	where the first arrow is $\Hom (X,\, \ad(S) )$ and the second is $-^\sharp$.
	Similarly there are natural arrows
	\[\bullet \hotimes S : \Hom (X, Y) \to \Hom (X \hotimes S, Y \hotimes S).\]
\end{remark}
	
\subsection{Strongly linear algebras and modules}\label{ssec:SLAlg}

Having a monoidal structure available on $\Sigma\Vect$ one can talk about strongly linear algebras as monoids in $\Sigma\Vect$, and of modules for those algebras.
	
\begin{definition}\label{defn:slalg}
	In $(\Ind\-(\Vect^\op), \Yon \bfk, -\hotimes- )$ a monoid is a triple $(X, \mu, u)$ with $X \in \Ind\-(\Vect^\op)$, $\mu  \in \Nat(X \hotimes X, X)$, $u \in \Nat(\Yon \bfk, X)\cong X\bfk$ such that
	\[\mu \cdot (\mu \hotimes  1_X) = \mu \cdot (1_X \hotimes  \mu)\]
	and $\mu \cdot (u \hotimes 1_X)$ and $\mu \cdot (u \hotimes 1_X)$ are respectively the left and right unitors of the monoidal structure $(- \hotimes  -, \Yon \bfk)$. 		
	A monoid  $(X,\mu,u)$ in $\Ind\-(\Vect^\op)$ with $X \in \Sigma\Vect$ will be called an associative \emph{strongly linear $\bfk$-algebra}.
	We will usually denote by $1\in X\bfk$ the image of $1 \in \bfk$ by $(u\bfk): \bfk \to X\bfk$.
\end{definition}

\begin{remark}
    A strongly linear $\bfk$-algebra is ``the same'' as a monoid in $\Sigma\Vect$ with the induced tensor product $\cR_3^\infty(-\hotimes -)$. This is because for $Z\in \Sigma\Vect$ there is a bijective correspondence between maps $\cR_3^\infty(X\hotimes Y) \to Z$ and maps $X \hotimes Y \to Z$.
\end{remark}
 
\begin{remark}[Relation with \emph{summability algebras} in \cite{bagayoko2024automorphisms}]
	If $(X, \mu, u)$ is a strongly linear $\bfk$-algebra and $\eta: X \times X \to X \hotimes X$ is as in Remark~\ref{rmk:hotimes_and_strongbilinear} then
    \[(\mu \cdot \eta)(\bfk) : X\bfk \times X\bfk \to X \bfk \]
    defines an associative $\bfk$-algebra structure on $X\bfk$ which preserves infinite sums in both arguments and conversely any associative $\bfk$-algebra structure on $X\bfk$ which preserves infinite sums in both arguments comes from a  unique such $\mu$. Hence a strong $\bfk$-algebra can be equivalently defined as a strong vector space $X$ together with a multiplication $*: X \bfk \times X\bfk \to X\bfk$ and a distinguished element $1 \in X\bfk$ such that
	\begin{enumerate}
		\item if $(f_i)_{i \in I}$ and $(g_j)_{j \in J}$ are summable families in $X\bfk$ then the family $(f_i * g_j)_{(i,j) \in I\times J}$ is summable and for all $(k_i)_{i \in I}\in \bfk^I$ and $(l_j)_{j \in J}\in \bfk^J$,
		\[\left(\sum_{i} f_i k_i\right) * \left(\sum_{j} g_j l_j\right)=\sum_{i, j} (f_i * g_j) \cdot (k_i\cdot l_j)\]
		\item $1*f = f*1$ for every $f \in X\bfk$;
		\item $f * (g * h) = (f * g) * h$ for every $f,g,h \in X\bfk$.
	\end{enumerate}   
    In other words, a strongly linear $\bfk$-algebra is essentially the same as a \emph{summability algebra} in \cite[Sec.~3]{bagayoko2024automorphisms}.
\end{remark}
	
\begin{example}
	If $X \in \Sigma\Vect$, then $\Hom (X, X)$ is a strongly linear algebra with composition (or rather with the map induced by composition). 
\end{example}
	
\begin{definition}
	A strongly linear algebra is \emph{commutative} (or \emph{Abelian}) if $\mu \cdot \tau = \mu$.
\end{definition}
	
\begin{example}\label{ex:algebra_of_series}
	Let $(G, \cdot, 1)$ be a monoid and $\cF$ be a bornology on $G$ such that $F_0 \cdot F_1  \in \cF$ whenever $F_0$ and $F_1$ are in $\cF$ and for every $g \in G$, $\{(g_0, g_1) \in F_0 \times F_1: g_0 g_1=g\}$ is finite. Call a pair $(G, \cF)$ as just described a \emph{bornological monoid}. On $\bfk[G, \cF]$ there is a unique strongly linear algebra structure making the inclusion $G \subseteq \bfk[G, \cF](\bfk)$ a monoid embedding. We can call it an \emph{algebra of series over} $(G, \cF)$. Observe that if $G$ is Abelian the s.l.\ algebra $\bfk[G, \cF]$ is Abelian.
\end{example}
	
\begin{definition}
	If $(R, \mu, u)$ is a strongly linear algebra, a \emph{left strongly linear $R$-module} is a pair $(M, \lambda)$ where $M \in \Sigma\Vect$ and $\lambda : R \hotimes  M \to M$ is such that $\lambda \cdot (\mu \hotimes  M) = \lambda \cdot (R \hotimes  \lambda)$. A morphism of s.l.\ left $R$-modules $f: M \to N$ is an $f \in \Sigma\Vect (M, N)$ which preserves the left action of $R$ i.e. $\lambda_N \cdot (R \hotimes f) = f \cdot \lambda_M$. 
\end{definition}

\begin{remark}
	In general morphisms of left strongly linear $R$-modules are naturally points for a $\Sigma$-closed sub-object of $\Hom (M, N)$ given by the kernel of the following map $\Hom (M, N) \to \Hom (R \hotimes M, N)$
	\[\Hom (\lambda_M, N) - \Hom (R \hotimes M, \lambda_N) \cdot (R \hotimes \bullet)\]
	We denote such a subspace by $R\Mod(M, N).$
\end{remark}
	
\begin{remark}
	If $(R, \mu, u)$ is a commutative s.l.\ algebra and $M, N$ are s.l.\ $R$-modules then $R\Mod(M, N)$ is naturally a s.l. $R$-module. 
\end{remark}
	
\begin{example}
	Let $\Gamma$ be a set, $\cF$ be a bornology on $\Gamma$. We can define 
	\[R[\Gamma; \cF]=\bigcup_{F \in \cF} R^F.\]
	Note that a family $(f_i)_{i \in I}$ in $R[\Gamma; \cF]$ is summable if and only if for every $\gamma \in \Gamma$, $\{f_{i, \gamma}: i \in I\}$ is summable in $R$ and $\{\gamma \in \Gamma: \exists i \in I , f_{i, \gamma}\neq 0\}\in \cF$.
\end{example}
	
\begin{example}
	Assume $(G, \cF)$ is a bornological monoid as in Example~\ref{ex:algebra_of_series}, and $X$ is a free left $G$-set with $\cH$ a bornology on $X$ such that $F \cdot H \in \cH$ for every $F \in \cF$ and every $H \in \cH$. Then $\bfk[X, \cH]$ has a natural structure of left $R= \bfk[G, \cF]$-module.
\end{example}

\subsection{Strongly linear Kähler differentials}\label{ssec:SLDifferential}
	
In the following assume $S$ is a given commutative strongly linear algebra.
	
\begin{definition}\label{defn:sl_diff}
	If $(M, \lambda)$ is a s.l.\ $S$-module, a \emph{$M$-valued strongly linear derivation} is a strongly linear map $\partial: S \to M$, such that: $\partial \cdot \mu= \lambda \cdot (S\hotimes \partial) + \lambda \cdot \tau \cdot (\partial \hotimes S)$.
	Strongly linear derivations form the points of the s.l.\ vector space given by the kernel of the map $\Hom (S, M) \to \Hom (S\hotimes S, M)$
	\[\Hom (\mu, M) - \Hom (S \hotimes S, \lambda) \cdot (S \hotimes \bullet) - \Hom (S \hotimes S, \lambda\cdot \tau) \cdot (\bullet \hotimes S).\]
	We denote the s.l.\ vector space of derivations by $\Sigma \Deriv(S, M)$. It is also easy to see that $\Sigma \Deriv(S, M)$ is naturally a s.l.\ $S$-module.
\end{definition}

A consequence of the completeness and cocompleteness of $\Sigma\Vect$ is that strongly linear derivations are representable in the sense that there is a s.l.\ $S$-module $\Omega_{S/\bfk}^{\Sigma}$ such that $S \Mod (\Omega_{S/\bfk}^{\Sigma}, M)\cong \Sigma \Deriv(S, M)$. Let $dS \in \Sigma\Vect$ denote an isomorphic copy of $S$ just as an object in $\Sigma\Vect$ and $\bar{d}\in \Sigma\Vect(S, \cR^\infty_3(S \hotimes dS))$ denote the composition 
\[\begin{tikzcd}
	S \cong \bfk \hotimes S \ar[r, "{1 \hotimes S}"]& S\hotimes S \ar[r, "{\rho^\infty}"] &\cR^\infty_3(S \hotimes S).
\end{tikzcd}\]
Then $f \mapsto f \cdot \bar{d}$ is an isomorphism $S\Mod (\cR^\infty_3(S \hotimes dS), M) \cong \Sigma\Vect (S, M)$, where $S \hotimes_\bfk dS$ is regarded as a strong $S$-module with action on the left factor. It follows that it suffices to define $\Omega_{S/\bfk}^\Sigma$ as
\[\Omega_{S/\bfk}^\Sigma = \frac{\cR^\infty_\cQ(S \hotimes_\bfk dS)}{H}\]
where $H= \bigcap \{\Kr(f): M \in S\Mod,\; f \cdot \bar{d} \in \Sigma \Deriv(S, M)\}$ provided that the composition $d$ of $\bar{d}$ with the quotient map is a strongly linear derivation. In fact it turns out to be the case: it is strongly linear by construction, and is a derivation because for every $s_1, s_2 \in S$, $\bar{d}(s_1  s_2) - s_1 \bar{d} s_2 -s_2 \bar{d}s_1$ is annihilated by all derivations and so is in $H$.
	
\begin{question}
	Is there a nice description of the strongly linear Kähler differential of the Hahn field $\bR[\bfx^\bR; \WO(\bfx^\bR)]$? Is its underlying strong vector space not in the essential image of $\bfk[-;-]$?
\end{question}

\bibliography{Res}

@book {maclane71categories,
    AUTHOR = {MacLane, Saunders},
     TITLE = {Categories for the working mathematician},
    SERIES = {Graduate Texts in Mathematics},
    VOLUME = {Vol. 5},
 PUBLISHER = {Springer-Verlag, New York-Berlin},
      YEAR = {1971},
     PAGES = {ix+262},
   MRCLASS = {18-02},
  MRNUMBER = {354798},
MRREVIEWER = {H.-B.\ Brinkmann},
}

@article{kelly1980unified,
    AUTHOR = {Kelly, G. M.},
     TITLE = {A unified treatment of transfinite constructions for free
              algebras, free monoids, colimits, associated sheaves, and so
              on},
   JOURNAL = {Bull. Austral. Math. Soc.},
  FJOURNAL = {Bulletin of the Australian Mathematical Society},
    VOLUME = {22},
      YEAR = {1980},
    NUMBER = {1},
     PAGES = {1--83},
      ISSN = {0004-9727},
   MRCLASS = {18C15 (18A40)},
  MRNUMBER = {589937},
MRREVIEWER = {V\'aclav\ Koubek},
       DOI = {10.1017/S0004972700006353},
       URL = {https://doi.org/10.1017/S0004972700006353},
}

@article{neumann1949ordered,
    AUTHOR = {Neumann, B. H.},
     TITLE = {On ordered division rings},
   JOURNAL = {Trans. Amer. Math. Soc.},
  FJOURNAL = {Transactions of the American Mathematical Society},
    VOLUME = {66},
      YEAR = {1949},
     PAGES = {202--252},
      ISSN = {0002-9947,1088-6850},
   MRCLASS = {09.1X},
  MRNUMBER = {32593},
MRREVIEWER = {R.\ Moufang},
       DOI = {10.2307/1990552},
       URL = {https://doi.org/10.2307/1990552},
}

@article{kaplansky1942maximal,
    AUTHOR = {Kaplansky, Irving},
     TITLE = {Maximal fields with valuations},
   JOURNAL = {Duke Math. J.},
  FJOURNAL = {Duke Mathematical Journal},
    VOLUME = {9},
      YEAR = {1942},
     PAGES = {303--321},
      ISSN = {0012-7094,1547-7398},
   MRCLASS = {09.1X},
  MRNUMBER = {6161},
MRREVIEWER = {Saunders\ Mac Lane},
       URL = {http://projecteuclid.org/euclid.dmj/1077493226},
}

@article {kuhlmann2005kappa,
    AUTHOR = {Kuhlmann, Salma and Shelah, Saharon},
     TITLE = {{$\kappa$}-bounded exponential-logarithmic power series
              fields},
   JOURNAL = {Ann. Pure Appl. Logic},
  FJOURNAL = {Annals of Pure and Applied Logic},
    VOLUME = {136},
      YEAR = {2005},
    NUMBER = {3},
     PAGES = {284--296},
      ISSN = {0168-0072,1873-2461},
   MRCLASS = {03C60 (06A05)},
  MRNUMBER = {2169687},
MRREVIEWER = {Martin\ Weese},
       DOI = {10.1016/j.apal.2005.04.001},
       URL = {https://doi.org/10.1016/j.apal.2005.04.001},
}

@article {berarducci2015surreal,
    AUTHOR = {Berarducci, Alessandro and Mantova, Vincenzo},
     TITLE = {Surreal numbers, derivations and transseries},
   JOURNAL = {J. Eur. Math. Soc. (JEMS)},
  FJOURNAL = {Journal of the European Mathematical Society (JEMS)},
    VOLUME = {20},
      YEAR = {2018},
    NUMBER = {2},
     PAGES = {339--390},
      ISSN = {1435-9855,1435-9863},
   MRCLASS = {03C64 (12H05 13N15 16W60 26A12)},
  MRNUMBER = {3760298},
       DOI = {10.4171/JEMS/769},
       URL = {https://doi.org/10.4171/JEMS/769},
}

@book{kashiwara2005categories,
  title={Categories and Sheaves},
  author={Kashiwara, M. and Schapira, P.},
  isbn={9783540279495},
  lccn={2005930329},
  series={Grundlehren der mathematischen Wissenschaften},
  url={https://books.google.it/books?id=K-SjOw\_2gXwC},
  year={2005},
  publisher={Springer Berlin Heidelberg}
}

@book {lefschetz1942algebraic,
    AUTHOR = {Lefschetz, Solomon},
     TITLE = {Algebraic {T}opology},
    SERIES = {American Mathematical Society Colloquium Publications},
    VOLUME = {Vol. 27},
 PUBLISHER = {American Mathematical Society, New York},
      YEAR = {1942},
     PAGES = {vi+389},
   MRCLASS = {56.0X},
  MRNUMBER = {7093},
MRREVIEWER = {H.\ Whitney},
}

@article{hoeven2001operators,
AUTHOR = {{\noopsort{hoeven}van der Hoeven}, Joris},
     TITLE = {Operators on generalized power series},
   JOURNAL = {Illinois J. Math.},
  FJOURNAL = {Illinois Journal of Mathematics},
    VOLUME = {45},
      YEAR = {2001},
    NUMBER = {4},
     PAGES = {1161--1190},
      ISSN = {0019-2082,1945-6581},
   MRCLASS = {16W60},
  MRNUMBER = {1894891},
MRREVIEWER = {N.\ Sankaran},
       URL = {http://projecteuclid.org/euclid.ijm/1258138061},
}

@article{dries2001logarithmic,
 AUTHOR = {{\noopsort{dries}van den Dries}, Lou and Macintyre, Angus and Marker, David},
     TITLE = {Logarithmic-exponential series},
 BOOKTITLE = {Proceedings of the {I}nternational {C}onference ``{A}nalyse \&
              {L}ogique'' ({M}ons, 1997)},
   JOURNAL = {Ann. Pure Appl. Logic},
  FJOURNAL = {Annals of Pure and Applied Logic},
    VOLUME = {111},
      YEAR = {2001},
    NUMBER = {1-2},
     PAGES = {61--113},
      ISSN = {0168-0072,1873-2461},
   MRCLASS = {12J15 (03H05 12H05 26E05)},
  MRNUMBER = {1848569},
MRREVIEWER = {Niels\ Schwartz},
       DOI = {10.1016/S0168-0072(01)00035-5},
       URL = {https://doi.org/10.1016/S0168-0072(01)00035-5},
}

@book{adamek1994locally,
    AUTHOR = {Ad\'{a}mek, Ji\v{r}\'{\i} and Rosick\'{y}, Ji\v{r}\'{\i}},
     TITLE = {Locally presentable and accessible categories},
    SERIES = {London Mathematical Society Lecture Note Series},
    VOLUME = {189},
 PUBLISHER = {Cambridge University Press, Cambridge},
      YEAR = {1994},
     PAGES = {xiv+316},
      ISBN = {0-521-42261-2},
   MRCLASS = {18Axx (18-02)},
  MRNUMBER = {1294136},
       DOI = {10.1017/CBO9780511600579},
       URL = {https://doi.org/10.1017/CBO9780511600579},
}

@article {diliberti2020codensity,
    AUTHOR = {{\noopsort{liberti}Di Liberti}, Ivan},
     TITLE = {Codensity: {I}sbell duality, pro-objects, compactness and
              accessibility},
   JOURNAL = {J. Pure Appl. Algebra},
  FJOURNAL = {Journal of Pure and Applied Algebra},
    VOLUME = {224},
      YEAR = {2020},
    NUMBER = {10},
     PAGES = {106379, 18},
      ISSN = {0022-4049,1873-1376},
   MRCLASS = {18C15 (18A40 18A99 18C35 20E18 54H99)},
  MRNUMBER = {4093066},
MRREVIEWER = {Partha\ Pratim\ Ghosh},
       DOI = {10.1016/j.jpaa.2020.106379},
       URL = {https://doi.org/10.1016/j.jpaa.2020.106379},
}

@article {higman1952ordering,
    AUTHOR = {Higman, Graham},
     TITLE = {Ordering by divisibility in abstract algebras},
   JOURNAL = {Proc. London Math. Soc. (3)},
  FJOURNAL = {Proceedings of the London Mathematical Society. Third Series},
    VOLUME = {2},
      YEAR = {1952},
     PAGES = {326--336},
      ISSN = {0024-6115,1460-244X},
   MRCLASS = {09.1X},
  MRNUMBER = {49867},
MRREVIEWER = {D.\ Zelinsky},
       DOI = {10.1112/plms/s3-2.1.326},
       URL = {https://doi.org/10.1112/plms/s3-2.1.326},
}

@article{berarducci2021value,
    AUTHOR = {Berarducci, Alessandro and Freni, Pietro},
     TITLE = {On the value group of the transseries},
   JOURNAL = {Pacific J. Math.},
  FJOURNAL = {Pacific Journal of Mathematics},
    VOLUME = {312},
      YEAR = {2021},
    NUMBER = {2},
     PAGES = {335--354},
      ISSN = {0030-8730,1945-5844},
   MRCLASS = {16W60 (03C64 12L12)},
  MRNUMBER = {4305776},
MRREVIEWER = {Lothar\ Sebastian\ Krapp},
       DOI = {10.2140/pjm.2021.312.335},
       URL = {https://doi.org/10.2140/pjm.2021.312.335},
}

@article {freyd1972categories,
    AUTHOR = {Freyd, P. J. and Kelly, G. M.},
     TITLE = {Categories of continuous functors. {I}},
   JOURNAL = {J. Pure Appl. Algebra},
  FJOURNAL = {Journal of Pure and Applied Algebra},
    VOLUME = {2},
      YEAR = {1972},
     PAGES = {169--191},
      ISSN = {0022-4049,1873-1376},
   MRCLASS = {18A25},
  MRNUMBER = {322004},
MRREVIEWER = {M.\ Tierney},
       DOI = {10.1016/0022-4049(72)90001-1},
       URL = {https://doi.org/10.1016/0022-4049(72)90001-1},
}

@article{krapp2022rayner,
 AUTHOR = {Krapp, Lothar Sebastian and Kuhlmann, Salma and Serra,
              Michele},
     TITLE = {On {R}ayner structures},
   JOURNAL = {Comm. Algebra},
  FJOURNAL = {Communications in Algebra},
    VOLUME = {50},
      YEAR = {2022},
    NUMBER = {3},
     PAGES = {940--948},
      ISSN = {0092-7872,1532-4125},
   MRCLASS = {13J05 (06F20 12J20 16W60 20K01)},
  MRNUMBER = {4379648},
MRREVIEWER = {G\'{e}rard\ Leloup},
       DOI = {10.1080/00927872.2021.1976789},
       URL = {https://doi.org/10.1080/00927872.2021.1976789},
}

@article {leinster2012codensity,
    AUTHOR = {Leinster, Tom},
     TITLE = {Codensity and the ultrafilter monad},
   JOURNAL = {Theory Appl. Categ.},
  FJOURNAL = {Theory and Applications of Categories},
    VOLUME = {28},
      YEAR = {2013},
     PAGES = {No. 13, 332--370},
      ISSN = {1201-561X},
   MRCLASS = {18C15 (03C20 18A30 18A40)},
  MRNUMBER = {3080612},
MRREVIEWER = {Lutz\ Schr\"oder},
}

@article {positselski2015contramodules,
    AUTHOR = {Positselski, Leonid},
     TITLE = {Contramodules},
   JOURNAL = {Confluentes Math.},
  FJOURNAL = {Confluentes Mathematici},
    VOLUME = {13},
      YEAR = {2021},
    NUMBER = {2},
     PAGES = {93--182},
      ISSN = {1793-7434},
   MRCLASS = {16T15 (16W60 16W80 17B65 18E10 22D12)},
  MRNUMBER = {4400900},
MRREVIEWER = {Paolo\ Saracco},
}

@book {loregian2021coend,
    AUTHOR = {Loregian, Fosco},
     TITLE = {({C}o)end calculus},
    SERIES = {London Mathematical Society Lecture Note Series},
    VOLUME = {468},
 PUBLISHER = {Cambridge University Press, Cambridge},
      YEAR = {2021},
     PAGES = {xxi+308},
      ISBN = {978-1-108-74612-0},
   MRCLASS = {18-02 (18A40 18D60 18D70 18M60)},
  MRNUMBER = {4274071},
MRREVIEWER = {Nicola\ Gambino},
       DOI = {10.1017/9781108778657},
       URL = {https://doi.org/10.1017/9781108778657},
}

@misc{bagayoko2024automorphisms,
      title={Automorphisms and derivations on algebras endowed with formal infinite sums}, 
      author={Vincent Bagayoko and Lothar Sebastian Krapp and Salma Kuhlmann and Daniel Panazzolo and Michele Serra},
      year={2024},
      eprint={2403.05827},
      archivePrefix={arXiv},
      primaryClass={math.RA},
      url={https://arxiv.org/abs/2403.05827}, 
}

@article {maclane1948groups,
    AUTHOR = {MacLane, Saunders},
     TITLE = {Groups, categories and duality},
   JOURNAL = {Proc. Nat. Acad. Sci. U.S.A.},
  FJOURNAL = {Proceedings of the National Academy of Sciences of the United States of America},
    VOLUME = {34},
      YEAR = {1948},
     PAGES = {263--267},
      ISSN = {0027-8424},
   MRCLASS = {20.0X},
  MRNUMBER = {25464},
MRREVIEWER = {P.\ Hall},
       DOI = {10.1073/pnas.34.6.263},
       URL = {https://doi.org/10.1073/pnas.34.6.263},
}

@article {isbell1957remarks,
    AUTHOR = {Isbell, J. R.},
     TITLE = {Some remarks concerning categories and subspaces},
   JOURNAL = {Canadian J. Math.},
  FJOURNAL = {Canadian Journal of Mathematics. Journal Canadien de
              Math\'ematiques},
    VOLUME = {9},
      YEAR = {1957},
     PAGES = {563--577},
      ISSN = {0008-414X,1496-4279},
   MRCLASS = {20.00 (08.00)},
  MRNUMBER = {94405},
MRREVIEWER = {S.\ Eilenberg},
       DOI = {10.4153/CJM-1957-064-6},
       URL = {https://doi.org/10.4153/CJM-1957-064-6},
}

@article {dickson1966torsion,
    AUTHOR = {Dickson, Spencer E.},
     TITLE = {A torsion theory for {A}belian categories},
   JOURNAL = {Trans. Amer. Math. Soc.},
  FJOURNAL = {Transactions of the American Mathematical Society},
    VOLUME = {121},
      YEAR = {1966},
     PAGES = {223--235},
      ISSN = {0002-9947,1088-6850},
   MRCLASS = {18.15 (18.20)},
  MRNUMBER = {191935},
MRREVIEWER = {E.\ A.\ Walker},
       DOI = {10.2307/1994341},
       URL = {https://doi.org/10.2307/1994341},
}

@article {grandis2006natural,
    AUTHOR = {Grandis, Marco and Tholen, Walter},
     TITLE = {Natural weak factorization systems},
   JOURNAL = {Arch. Math. (Brno)},
  FJOURNAL = {Universitatis Masarykianae Brunensis. Facultas Scientiarum
              Naturalium. Archivum Mathematicum},
    VOLUME = {42},
      YEAR = {2006},
    NUMBER = {4},
     PAGES = {397--408},
      ISSN = {0044-8753,1212-5059},
   MRCLASS = {18C15},
  MRNUMBER = {2283020},
MRREVIEWER = {Thomas\ Tradler},
}

@article {rosicky2008factorization,
    AUTHOR = {Rosick\'y, Ji\v{r}\'{\i} and Tholen, Walter},
     TITLE = {Factorization, fibration and torsion},
   JOURNAL = {J. Homotopy Relat. Struct.},
  FJOURNAL = {Journal of Homotopy and Related Structures},
    VOLUME = {2},
      YEAR = {2007},
    NUMBER = {2},
     PAGES = {295--314},
      ISSN = {1512-2891},
   MRCLASS = {18E40 (18D30)},
  MRNUMBER = {2369170},
MRREVIEWER = {R.\ H.\ Street},
}

@article {cassidy1985reflective,
    AUTHOR = {Cassidy, C. and H\'ebert, M. and Kelly, G. M.},
     TITLE = {Reflective subcategories, localizations and factorization
              systems},
   JOURNAL = {J. Austral. Math. Soc. Ser. A},
  FJOURNAL = {Australian Mathematical Society. Journal. Series A. Pure
              Mathematics and Statistics},
    VOLUME = {38},
      YEAR = {1985},
    NUMBER = {3},
     PAGES = {287--329},
      ISSN = {0263-6115},
   MRCLASS = {18A20 (18E35)},
  MRNUMBER = {779198},
MRREVIEWER = {Ji\v r\'i\ Ad\'amek},
}

@book {borceuxHCA2,
    AUTHOR = {Borceux, Francis},
     TITLE = {Handbook of categorical algebra. 2},
    SERIES = {Encyclopedia of Mathematics and its Applications},
    VOLUME = {51},
      NOTE = {Categories and structures},
 PUBLISHER = {Cambridge University Press, Cambridge},
      YEAR = {1994},
     PAGES = {xviii+443},
      ISBN = {0-521-44179-X},
   MRCLASS = {18-02 (18Exx)},
  MRNUMBER = {1313497},
MRREVIEWER = {Martin\ Hyland},
}

@unpublished{freni2024structural,
            year = {2024},
           title = {Structural Investigations in some Classes of o-minimal Fields},
            note = {Unpublished},
          school = {University of Leeds},
           month = {September},
          author = {Freni, Pietro},
             url = {https://etheses.whiterose.ac.uk/id/eprint/36168/}
}

@article{hahn1907nichtarchimedeischen,
    author = {Hahn, Hans},
    year = {1907},
    title = {{\"U}ber die nichtarchimedeischen {G}rossensysteme},
    journal = {Sitzungberichte Keiserlichen der Akademie Wissenschaften Wien, Mathematisch - Naturwissenschaflichte Klasse },
}

@article {stone1936theory,
    AUTHOR = {Stone, M. H.},
     TITLE = {The theory of representations for {B}oolean algebras},
   JOURNAL = {Trans. Amer. Math. Soc.},
  FJOURNAL = {Transactions of the American Mathematical Society},
    VOLUME = {40},
      YEAR = {1936},
    NUMBER = {1},
     PAGES = {37--111},
      ISSN = {0002-9947,1088-6850},
   MRCLASS = {06E20},
  MRNUMBER = {1501865},
       DOI = {10.2307/1989664},
       URL = {https://doi.org/10.2307/1989664},
}

@book {johnstone1982stone,
    AUTHOR = {Johnstone, Peter T.},
     TITLE = {Stone spaces},
    SERIES = {Cambridge Studies in Advanced Mathematics},
    VOLUME = {3},
 PUBLISHER = {Cambridge University Press, Cambridge},
      YEAR = {1982},
     PAGES = {xxi+370},
      ISBN = {0-521-23893-5},
   MRCLASS = {54-02 (06Dxx 18B30)},
  MRNUMBER = {698074},
MRREVIEWER = {Andreas\ Blass},
}

@article {ringel1970diagonalisierungspaare,
    AUTHOR = {Ringel, Claus Michael},
     TITLE = {Diagonalisierungspaare. {I}},
   JOURNAL = {Math. Z.},
  FJOURNAL = {Mathematische Zeitschrift},
    VOLUME = {117},
      YEAR = {1970},
     PAGES = {249--266},
      ISSN = {0025-5874,1432-1823},
   MRCLASS = {18.10},
  MRNUMBER = {272864},
MRREVIEWER = {J.\ R.\ Isbell},
       DOI = {10.1007/BF01109846},
       URL = {https://doi.org/10.1007/BF01109846},
}

@misc{riehl2008factorization,
  title={Factorization systems},
  author={Riehl, Emily},
  note={Notes available at {\texttt{http://www.math.jhu.edu/\~{}eriehl/factorization.pdf}}},
  year={2008}
}

@article {cech1937bicompact,
    AUTHOR = {{\v C}ech, Eduard},
     TITLE = {On bicompact spaces},
   JOURNAL = {Ann. of Math. (2)},
  FJOURNAL = {Annals of Mathematics. Second Series},
    VOLUME = {38},
      YEAR = {1937},
    NUMBER = {4},
     PAGES = {823--844},
      ISSN = {0003-486X,1939-8980},
   MRCLASS = {99-04},
  MRNUMBER = {1503374},
       DOI = {10.2307/1968839},
       URL = {https://doi.org/10.2307/1968839},
}

@article {ribenboim1997semisimple,
    AUTHOR = {Ribenboim, Paulo},
     TITLE = {Semisimple rings and von {N}eumann regular rings of
              generalized power series},
   JOURNAL = {J. Algebra},
  FJOURNAL = {Journal of Algebra},
    VOLUME = {198},
      YEAR = {1997},
    NUMBER = {2},
     PAGES = {327--338},
      ISSN = {0021-8693,1090-266X},
   MRCLASS = {16E50 (16S36)},
  MRNUMBER = {1489900},
MRREVIEWER = {Giuseppe\ Baccella},
       DOI = {10.1006/jabr.1997.7063},
       URL = {https://doi.org/10.1006/jabr.1997.7063},
}
\bibliographystyle{plain}
	
\end{document}